\documentclass[12pt, leqno]{amsart}

\usepackage{shuffle} 

\usepackage{enumerate}

\usepackage{upgreek}

\usepackage{lscape}
\usepackage{graphicx}
\usepackage{amssymb,amsmath,amsthm,amscd}
\usepackage{mathrsfs}
\usepackage[usenames,dvipsnames]{color}
\usepackage[colorlinks=true, pdfstartview=FitV,
 linkcolor=blue,citecolor=blue,urlcolor=blue]{hyperref}

\usepackage[all]{xy}
\usepackage[normalem]{ulem}  

\allowdisplaybreaks[4]
\usepackage{oldgerm}
\usepackage{xspace}

\usepackage{marginnote}

\newcommand{\nc}{\newcommand}
\numberwithin{equation}{section}

\newenvironment{red}{\relax\color{red}}{\relax}
\newenvironment{blue}{\relax\color{Dandelion}}{\hspace*{.5ex}\relax}

\newcommand{\beb}{\begin{blue}}
\newcommand{\eb}{\end{blue}}

\newcommand{\ber}{\begin{red}}
\newcommand{\er}{\end{red}}

\newcommand{\berm}[1]{\begin{red}{}\marginnote{\fbox{\scshape\lowercase{M}}}%
#1}  

\newcommand{\berE}[1]{\begin{red}{}\marginnote{\fbox{\scshape\lowercase{E}}}%
#1}  

\newcommand{\berMH}[1]{\begin{red}{}\marginnote{\fbox{\scshape\lowercase{MH}}}%
#1}  

\newcommand{\berS}[1]{\begin{red}{}\marginnote{\fbox{\scshape\lowercase{S}}}%
#1}  

\theoremstyle{plain}
\newtheorem{lemma}{Lemma}[section]
\newtheorem{prop}[lemma]{Proposition}
\newtheorem{theorem}[lemma]{Theorem}

\newcommand{\Prop}{\begin{prop}}
\newcommand{\enprop}{\end{prop}}
\newcommand{\Lemma}{\begin{lemma}}
\newcommand{\enlemma}{\end{lemma}}
\newcommand{\Th}{\begin{theorem}}
\newcommand{\enth}{\end{theorem}}
\newtheorem{corollary}[lemma]{Corollary}
\newcommand{\Cor}{\begin{corollary}}
\newcommand{\encor}{\end{corollary}}
\newtheorem{definition}[lemma]{Definition}
\newtheorem*{conjecture}{Conjecture}
\newcommand{\Def}{\begin{definition}}
\newcommand{\edf}{\end{definition}}
\newtheorem{sublemma}[lemma]{Sublemma}
\newcommand{\Sublemma}{\begin{sublemma}}
\newcommand{\ensub}{\end{sublemma}}
\newtheorem*{convention}{Convention}

\theoremstyle{definition}
\newtheorem{remark}[lemma]{Remark}
\newtheorem{example}[lemma]{Example}
\newtheorem{Convention}[lemma]{Convention}
\newcommand{\Conv}{\begin{Convention}}
\newcommand{\enconv}{\end{Convention}}
\nc{\Conj}{\begin{conjecture}}
\nc{\enconj}{\end{conjecture}}
\nc{\Rem}{\begin{remark}}
\nc{\enrem}{\end{remark}}
\newcommand{\C}{{\mathbb C}}
\newcommand{\Q}{\mathbb {Q}}

\newcommand{\Z}{{\mathbb Z}}
\newcommand{\B}{{\mathbf{B}}}
\newcommand{\A}{{\mathbf A}}

\newcommand{\D}{\mathscr{D}}

\newcommand{\q}{{\mathfrak{q}}}

\newcommand{\R}{{\rm R}}

\newcommand{\one}{{\bf{1}}}
\newcommand{\trivial}{\one}
\newcommand{\seteq}{\mathbin{:=}}

\newcommand{\hd}{{\operatorname{hd}}}
\newcommand{\soc}{{\operatorname{soc}}}

\newcommand{\g}{{\mathfrak{g}}}

\newcommand{\Hom}{\operatorname{Hom}}
\newcommand{\End}{\operatorname{End}}

\newcommand{\isoto}[1][]{\mathop{\xrightarrow%
[{\raisebox{.3ex}[0ex][.3ex]{$\scriptstyle{#1}$}}]%
{{\raisebox{-.6ex}[0ex][-.6ex]{$\mspace{2mu}\sim\mspace{2mu}$}}}}}

\newcommand{\M}{{\mathscr M}}

\newenvironment{myequation}
{\relax\setlength{\arraycolsep}{1pt}\begin{eqnarray}}
{\end{eqnarray}}
\newenvironment{myequationn}
{\relax\setlength{\arraycolsep}{1pt}\begin{eqnarray*}}
{\end{eqnarray*}}

\nc{\eq}{\begin{myequation}}
\nc{\eneq}{\end{myequation}}
\nc{\eqn}{\begin{myequationn}}
\nc{\eneqn}{\end{myequationn}}

\newcommand{\hs}{\hspace*}
\newcommand{\ms}{\mspace}

\newcommand{\To}[1][{\hs{2ex}}]{\xrightarrow{\,#1\,}}

\newcommand{\on}{\operatorname}
\newcommand{\Ker}{\on{Ker}}

\newcommand{\bnum}{\be[{\rm(i)}]}
\nc{\bni}{\bnum}
\nc{\bna}{\be[{\rm(a)}]}

\newcommand{\soplus}{\mathop{\mbox{\normalsize$\bigoplus$}}\limits}

\newcommand{\id}{\on{id}}
\newcommand{\ba}{\begin{array}}
\newcommand{\ea}{\end{array}}

\newcommand{\Coker}{{\operatorname{Coker}}}

\newcommand{\monoto}{\rightarrowtail}

\newcommand{\indlim}{\varinjlim\limits}
\newcommand{\set}[2]{\left\{#1 \mid #2 \right\}}

\newcommand{\eqsub}{\begin{subequations}\begin{eqnarray}}
\newcommand{\eneqsub}{\end{eqnarray}\end{subequations}}

\newcommand{\ol}{\overline}

\newcommand{\Pro}{{{\on{Pro}}}}
\newcommand{\ko}{{{\mathbf{k}}}}

\nc{\la}{\lambda}
\nc{\lam}{\lambda}
\nc{\U}[1][\g]{U_q(#1)}
\nc{\te}{\tilde{e}}
\nc{\tei}{\tilde{e}_i}
\nc{\tf}{\tilde{f}}
\nc{\tfi}{\tilde{f}_i}
\nc{\tU}{\widetilde U_q(\g)}
\nc{\tE}{\tilde{E}}
\nc{\tF}{\widetilde{\F}}
\nc{\tK}{\widetilde{K}}

\nc{\tk}{\tilde{k}}
\nc{\tkone}{\tk_{\ol{1}}}
\nc{\teone}{\tilde{e}_{\ol{1}}}
\nc{\tfone}{\tilde{f}_{\ol{1}}}

\nc{\teibar}{\tilde{e}_{\ol{i}}} \nc{\tfibar}{\tilde{f}_{\ol{i}}}
\nc{\tki}{{\tk}_{\ol {i}}}

\nc{\BZ}{{\mathbb{Z}}}
\nc{\al}{\alpha}
\nc{\qs}{{q}}
\nc{\lan}{\langle}
\nc{\ran}{\rangle}
\nc{\re}{{\mathrm{re}}}
\nc{\wt}{\operatorname{wt}}
\nc{\ch}{\operatorname{ch}}
\nc{\Um}[1][\g]{U^-_q(#1)}
\nc{\Ue}{U^+_q(\g)}
\nc{\eps}{\varepsilon}
\nc{\vphi}{\varphi}
\nc{\sphi}{\varphi^*}
\nc{\seps}{\varepsilon^*}

\nc{\nn}{\nonumber}
\def\max{{\mathop{\mathrm{max}}}}
\nc{\vp}{\varpi}
\nc{\cls}{{\operatorname{cl}}}
\nc{\Wt}{{\operatorname{Wt}}}
\nc{\Us}{U'_q(\g)}
\nc{\La}{\Lambda}
\nc{\tLa}{\widetilde\Lambda}
\nc{\ro}{{\rm(}}
\nc{\rf}{{\rm)}}
\nc{\norm}{{\mathrm{norm}}}
\nc{\qbox}{\quad\mbox}
\nc{\braid}{{\mathfrak{B}}}
\nc{\Ad}{\operatorname{Ad}}
\nc{\Aut}{\operatorname{Aut}}
\nc{\dt}[1]{\tilde{\tilde #1}}
\nc{\Sn}{S^{{\mathrm{norm}}}}
\nc{\aff}{{\rm{aff}}}
\nc{\rk}{{\mathrm{rk}}}
\nc{\tP}{\widetilde{P}}
\nc{\tW}{\widetilde{W}}
\nc{\Dyn}{\mathrm{Dyn}}
\nc{\tD}{\widetilde{\Delta}}
\nc{\height}[1]{{\operatorname{ht}}(#1)}
\nc{\bl}{\bigl(}
\nc{\br}{\bigr)}
\nc{\Hecke}{\mathrm{H}}
\nc{\HA}{\Hecke^{\mathrm{A}}}
\nc{\HB}{\Hecke^{\mathrm{B}}}
\newcommand{\scbul}{{\,\raise1pt\hbox{$\scriptscriptstyle\bullet$}\,}}
\nc{\vac}{{\phi}}
\nc{\Bt}{\B_\theta(\g)}
\nc{\be}{\begin{enumerate}}
\nc{\ee}{\end{enumerate}}
\nc{\low}{{\mathrm{low}}}
\nc{\upper}{{\mathrm{up}}}
\nc{\Zodd}{\Z_{\mathrm{odd}}}
\nc{\Ft}[1][n]{\mathbb{P}\mathrm{ol}_{#1}}
\nc{\Ftf}[1][n]{\widetilde{\mathbb{P}\mathrm{ol}}_{#1}}
\nc{\KA}{\on{K}^{\mathrm{A}}}
\nc{\KB}{\on{K}^{\mathrm{B}}}
\nc{\Res}{\on{Res}}
\nc{\Fc}[1][{n,m}]{\mathbf{F}_{#1}}
\nc{\tphi}{\tilde{\varphi}}
\nc{\CO}{\mathscr{O}}
\nc{\inte}{\mathrm{int}}
\nc{\Oint}{\mathcal{O}^{\ge0}_{\inte}}
\nc{\vs}{\vspace*}
\nc{\tLt}{\widetilde{L}}
\nc{\tL}{\widetilde{\Lambda}}
\nc{\tu}{\tilde{u}}
\nc{\noi}{\noindent}
\nc{\heigh}{\mathfrak{t}}
\nc{\lowest}{\mathfrak{l}}
\nc{\rootl}{\mathsf{Q}}
\nc{\rlQ}{\rootl}
\nc{\cl}{\colon}

\nc{\uqpg}{U'_q(\mathfrak g)}
\nc{\uq}{\uqpg}
\nc{\Oh}{\widehat{\mathcal{O}}}
\nc{\pn}{p_{\mathfrak{n}}}

\nc{\KLR}{KLR algebra}
\nc{\KLRs}{KLR algebras}
\nc{\cor}{\mathbf{k}}
\nc{\cora}{{\cor(A)}}
\nc{\haut}{\mathrm{ht}}
\nc{\tens}{\mathop\otimes}
\nc{\gmod}{\mbox{-$\mathrm{gmod}$}}
\nc{\gMod}{\mbox{-$\mathrm{gMod}$}}
\nc{\proj}{\mbox{-$\mathrm{proj}$}}
\nc{\gproj}{\mbox{-$\mathrm{gproj}$}}
\nc{\smod}{\mbox{-$\mathrm{mod}$}}
\nc{\Mod}{\mbox{-$\mathrm{gMod}$}}
\nc{\h}{\mathfrak h}
\nc{\Rnorm}{R^{\rm{norm}}}
\nc{\Runiv}{R^{\rm{univ}}}
\nc{\Rren}{R^{\rm{ren}}}

\nc{\Vhat}{\widehat{V}}
\nc{\F}{\mathcal{F}}

\def\T{{\mathcal T}}

\nc{\fd}[1][A]{\on{\mathrm{flat.dim}_{#1}}}
\nc{\bP}{{\mathbb{P}}}
\nc{\bPh}{\widehat{\mathbb{P}}}
\nc{\bK}[1][{n}]{\widehat{\mathbb{K}}_{#1}}
\nc{\bV}[1][{n}]{\widehat{V}^{\otimes{#1}}}
\nc{\bVK}[1][{n}]{\widehat{V}^{\otimes{#1}}_{\widehat{\mathbb{K}}}}
\nc{\hV}{\widehat{V}}
\nc{\opp}{\mathrm{opp}}
\nc{\col}{\colon}
\nc{\oep}{\epsilon}
\nc{\qtext}{\quad\text}
\nc{\qtextq}[1]{\quad\text{#1}\quad}
\nc{\longtwoheadrightarrow}[1][]{\xymatrix{\ar@{->>}[r]^-{{#1}}&}}
\nc{\epiTo}[1][]{\longtwoheadrightarrow[{#1}]}
\nc{\epito}{\twoheadrightarrow}
\nc{\monoTo}[1][]{\xymatrix{\ar@{>->}[r]^-{{#1}}&}}
\nc{\sym}{\mathfrak{S}}
\nc{\inp}[1]{{({#1})_{\mathrm{n}}}}
\nc{\rtl}{\rootl}
\nc{\wtd}{\widetilde}
\nc{\etens}{\boxtimes}
\nc{\ds}[1]{\mathrm{d}(#1)}
\nc{\rmat}[1]{{\mathbf{r}}_%
{\mspace{-2mu}\raisebox{-.6ex}{${\scriptstyle{#1}}$}}}
\nc{\rmats}[1]{{\mathbf{r}}_%
{\mspace{-2mu}\raisebox{-.6ex}{${\scriptscriptstyle{#1}}$}}}
\nc{\shc}{\mathcal{C}}
\nc{\shs}{\mathcal{S}}
\nc{\Fct}{{\on{Fct}}}
\nc{\tC}{\widetilde{\shc}}
\nc{\Zp}{\Z_{\ge0}}
\nc{\tPhi}{\widetilde{\Phi}}
\nc{\tT}{{\widetilde{\T}}}
\nc{\Ob}{\on{Ob}}
\nc{\bwr}{\mbox{\large$\wr$}}
\nc{\Img}{\on{Im}}
\nc{\Ab}{\mathcal{A}^{\mathrm{big}}}
\nc{\Sb}{\mathcal{S}^{\mathrm{big}}}
\nc{\As}{\mathcal{A}}
\nc{\Ss}{\mathcal{S}}
\nc{\ntens}{\widetilde{\otimes}}
\nc{\hR}{\widehat{R}}
\nc{\nconv}{\mathop{\mbox{\large $\odot$}}}
\nc{\snconv}{\mbox{\scriptsize$\odot$}}
\nc{\ts}{\tilde{s}}
\nc{\sho}{\mathcal{O}}
\nc{\bc}{\begin{cases}}
\nc{\ec}{\end{cases}}
\nc{\slnh}{{\widehat{\mathfrak{sl}}_N}}
\nc{\UA}{U_q'(\slnh)}
\nc{\KR}{R_K}
\nc{\cQ}{\mathcal{Q}}
\nc{\Irr}{\mathcal{I}rr}
\nc{\tQ}{\widetilde{\cQ}}
\nc{\bs}{\mathbf{s}}
\nc{\bL}{\mathbb{L}}
\nc{\tg}{\tilde{g}}
\nc{\lie}{\mathsf{g}}
\nc{\conv}{\mathbin{\mbox{\large $\circ$}}}
\nc{\shconv}{\mathbin{\large\diamond}}
\nc{\sconv}{\mathbin{\large\Delta}}
\nc{\stens}{\mathbin{\large\Delta}}
\nc{\hconv}{\mathbin{\nabla}}
\nc{\htens}{\mathbin{\nabla}}

\nc{\Rm}{R^{\mathrm{ren}}}

\nc{\bQ}{\ol{Q}}
\renewcommand{\Im}{\on{Im}}

\nc{\de}{\on{\textfrak{d}\ms{1mu}}}

\nc{\xmono}{\ar@{>->}}
\nc{\xepi}{\ar@{->>}}
\nc{\db}[1]{\raisebox{-.5ex}[2ex][1.8ex]{$#1$}}
\nc{\wb}[1]{\mbox{$\rule[-1.1ex]{0ex}{2ex}#1$}}
\nc{\univ}{\mathrm{univ}}
\nc{\rM}{{}^*\mspace{-2mu}M}
\nc{\lM}{M^*}
\nc{\uqm}{\uq\smod}
\nc{\tR}{\widetilde{R}_{\gamma,\beta}}
\nc{\tx}{\tilde{x}}
\nc{\bi}{\mathbf{i}}
\nc{\ttau}{\widetilde{\tau}}

\nc{\tEnd}{\on{\widetilde{E}nd}}
\nc{\tHom}{\on{\widetilde{H}om}}

\nc{\K}{{J}}
\nc{\Kex}{{\K}_{\mathrm{ex}}}
\nc{\Kfr}{{\K}_{\mathrm{f\mspace{.01mu}r}}}
\nc{\coro}{\cor}
\nc{\tB}{\widetilde{B}}
\nc{\seed}{\mathscr{S}}

\nc{\up}{\mathrm{up}}
\nc{\bfa}{\mathbf{a}}
\nc{\bfb}{\mathbf{b}}
\nc{\bfc}{\mathbf{c}}

\newlength{\mylength}
\nc{\ov}[1]{\overline{#1}}
\nc{\Wlmj}[3]{\W_{#2,#3}^{(#1)}}
\nc{\Mkl}[2]{\M_\ttww(#1,#2)}
\nc{\mqs}{(-q^2)}
\nc{\Cquiver}{\upsigma}

\nc{\mut}[1]{{\mu}_{\mspace{-2mu}\raisebox{-.5ex}{${\scriptstyle{#1}}$}}}

\nc{\Kt}{\mathcal K_t}
\nc{\KT}{\mathbb{K}_t}
\nc{\yim}{y_{i,m}}
\nc{\yjm}{y_{j,m}}
\nc{\yjp}{y_{j,p}}
\nc{\yimp}{y_{i,m+1}}
\nc{\yjmp}{y_{j,m+1}}

\nc{\Refl}{\mathscr{S}}
\nc{\Reflinv}{{\Refl}^{-1}}

\nc{\catC}{\mathscr C}
\nc{\catA}{\mathcal A}
\nc{\shift}{{\mathrm T}}
\nc{\rE}{ \mathsf{E} }
\nc{\rW}{ \mathcal{W} }
\nc{\rES}{ \mathcal{E} }

\nc{\brd}{\sigma} 
\nc{\into}{\xymatrix@C=3ex{{}\ar@{^{(}->}[r]&{}}}
\nc{\dual}{\D}
\nc{\cat}[1][{\g}]{\catC_{#1}}
\nc{\cato}[1][{\g}]{\catC_{#1}^0}

\nc{\catCO}{{\catC_\g^0}}
\nc{\catCQ}{{\catC_{\qQ}}}

\nc{\catCQd}{{\catC_{\widetilde{\qQ}}}}
\nc{\catCD}{{\catC_{\ddD}}}
\nc{\Li}{{\La^\infty}}
\nc{\sig}{{\sigma(\g)}}
\nc{\sigZ}{{\sigma_0(\g)}}
\nc{\sigQ}{{\sigma_\qQ(\g)}}
\nc{\sigQd}{{\sigma_{\widetilde{\qQ}}(\g)}}
\nc{\phiQd}{\phi_{\widetilde{\qQ}}}
\nc{\sigD}{{\sigma_\ddD(\g)}}

\nc{\ZZ}{{\mathbf{Z}}}
\nc{\sP}{{\mathsf{P}}}
\nc{\sV}{{\mathsf{V}}}
\nc{\rxw}{{\underline{w_0}}}
\nc{\boten}[1]{\overrightarrow{\bigotimes_{#1}}}
\nc{\cmA}{{\mathsf{A}}}
\nc{\cmC}{{\mathsf{C}}}
\nc{\ddD}{{\mathcal{D}}}
\nc{\ddDQ}{{\ddD_\qQ}}
\nc{\ddDQd}{{\ddD_{\widetilde{\qQ}}}}
\nc{\qQ}{{\mathscr{Q}}}
\nc{\gf}{{\g_{\rm fin}}}
\nc{\Df}{{\Delta_{\rm fin}}}
\nc{\prDf}{{\Phi^+_{\rm fin}}}
\nc{\If}{{I_{\rm fin}}}
\nc{\cmAf}{{\cmA_{\rm fin}}}
\nc{\weyl}{{\mathsf{W}}}
\nc{\sg}{{\mathfrak{S}}}
\nc{\weylA}{{\mathsf{W}_\cmA}}
\nc{\weylC}{{\mathsf{W}_\cmC}}
\nc{\weylf}{{\weyl_{\rm fin}}}
\nc{\Deg}{\mathrm{Deg}}
\nc{\Di}{\Deg^\infty}
\nc{\KRc}{{K_{q=1}(R_\cmC\gmod)}}

\nc{\prD}{{\Phi^+}}		
\nc{\nrD}{{\Phi^-}}
\nc{\prDA}{{\Phi^+_\cmA}}
\nc{\prDC}{{\Phi^+_\cmC}}
\nc{\nrDC}{{\Phi^-_\cmC}}
\nc{\rDC}{{\Phi_\cmC}}

\nc{\n}{{\mathfrak{n}}}
\nc{\Rt}{\ms{1.5mu}\mathsf{L}} 
\nc{\Cp}{\mathsf{V}} 
\nc{\cuspS}{\ms{1.5mu}{\mathsf{S}}}
\nc{\st}[1]{\{{#1}\}}
\nc{\bst}[1]{\bigl\{{#1}\bigr\}}
\nc{\WS}{quantum affine Schur-Weyl duality\xspace}
\nc{\CWS}{Quantum affine Schur-Weyl duality}
\nc{\zz}{{{\mathbf{z}}}}
\nc{\wlP}{\mathsf{P}}
\nc{\wl}{\wlP}
\nc{\clp}{{\mathrm{cl}}}
\nc{\wlPc}{{\wlP_\clp}}
\nc{\awlP}{\widehat{\mathsf{P}}}

\nc{\dM}{\mathsf{D}}	

\nc{\dC}{\mathsf{C}}
\nc{\cC}{\mathcal{C}}
\nc{\lR}{\widetilde{{R}}}
\nc{\zero}{\ms{3mu}\mathrm{zero}}
\nc{\PD}{principal }

\nc{\prtl}[1][J]{\rootl_{#1}^+}
\nc{\hL}{\widehat{\Rt}}
\nc{\hF}{\widehat{\F}}
\nc{\Proof}{\begin{proof}}
\nc{\QED}{\end{proof}}
\nc{\e}{\mathrm{e}}
\renewcommand{\P}{\mathrm{P}}
\nc{\Aff}{\mathrm{Aff}}
\nc{\rT}{\mathcal{T}}		
\nc{\rr}{rationally renormalizable\xspace}
\nc{\RA}{{R_\cmA}}		
\nc{\RC}{{R_\cmC}}		
\nc{\proolim}[1][]{\mathop{``{}\ms{1mu}{\varprojlim}{\mbox{''}}}\limits_{#1}}
\nc{\qtq}[1][\text{and}]{\quad\text{#1}\quad}
\newcommand{\gW}{\mathsf{W}}
\newcommand{\sgW}{\mathsf{W}^*}
\newcommand{\Sp}{\mathrm{span}_{\mathbb{R}_{\ge0}}}  	
\newcommand{\rF}{\mathcal{T}} 		

\nc{\corh}{\widehat{\cor}}
\nc{\ang}[1]{\langle{#1}\rangle}
\nc{\rc}{renormalizing coefficient\xspace}
\nc{\cz}{{\cor[z^{\pm1}]}}
\nc{\tp}{\ms{1.5mu}{\widetilde{p}}\ms{2mu}}
\nc{\G}{\mathcal{G}}
\nc{\cc}{\mathfrak{c}}
\nc{\rsP}{{\Upsilon_\g}}  		 
\nc{\rsX}{{X_\g}}
\nc{\rs}{ \mathsf{s} }
\nc{\Dynkin}{\triangle}		
\nc{\Dat}{\sigma}
\nc{\hf}{\xi}
\nc{\ake}[1][1ex]{\rule[-#1]{0ex}{1ex}}
\nc{\akew}[1][1ex]{\rule[-1ex]{#1}{0ex}}
\nc{\akeu}[1][1ex]{\rule[#1]{0ex}{1ex}}
\nc{\snoi}{\smallskip\noi}
\nc{\principal}{complete\xspace}
\nc{\Principal}{Complete\xspace}
\nc{\mnoi}{\medskip\noindent}
\nc{\vpi}{\varpi}
\nc{\qt}[1]{\quad\text{#1}}
\nc{\Fd}{\F_\ddD}
\nc{\hFd}{\hF_\ddD}
\nc{\zconv}{\mathop{\conv}\limits_z}
\nc{\phiQ}{{\upphi_\qQ}}
\nc{\catCab}[1]{\catC_\g^{#1}}
\nc{\Rgf}{{R^\gf}}

\title
[{PBW theory for quantum affine algebras}]
{PBW theory for quantum affine algebras}

\author[M. Kashiwara]{Masaki Kashiwara}
\thanks{The research of M.\ Kashiwara
was supported by Grant-in-Aid for Scientific Research (B)
20H01795, Japan Society for the Promotion of Science.}
\address[M. Kashiwara]{
Kyoto University Institute for Advanced Study,
Research Institute for Mathematical Sciences, Kyoto University,
Kyoto 606-8502, Japan \& Korea Institute for Advanced Study, Seoul 02455, Korea }
\email[M. Kashiwara]{masaki@kurims.kyoto-u.ac.jp}

\author[M. Kim]{Myungho Kim}
\address[M. Kim]{Department of Mathematics, Kyung Hee University, Seoul 02447, Korea}
\email[M. Kim]{mkim@khu.ac.kr}
\thanks{The research of M.\ Kim was supported by the National Research Foundation of
Korea(NRF) Grant funded by the Korea government(MSIP) (NRF-2017R1C1B2007824 and NRF-2020R1A5A1016126).}

\author[S.-j. Oh]{Se-jin Oh}
\thanks{ The research of S.-j.\ Oh was supported by the Ministry of Education of the Republic of Korea and the National Research Foundation of Korea (NRF-2019R1A2C4069647).}
\address[S.-j. Oh]{Department of Mathematics, Ewha Womans University, Seoul 03760, Korea}
\email[S.-j. Oh]{sejin092@gmail.com}

\author[E. Park]{Euiyong Park}
\thanks{The research of E.\ Park was supported by the National Research Foundation of Korea(NRF) Grant funded by the Korea Government(MSIP)(NRF-2020R1F1A1A01065992 and NRF-2020R1A5A1016126).}
\address[E. Park]{Department of Mathematics, University of Seoul, Seoul 02504, Korea}
\email[E. Park]{epark@uos.ac.kr}


\date{May 20, 2021}

\begin{document}

\maketitle

\begin{abstract}

Let $U_q'(\g)$ be a quantum affine algebra of \emph{arbitrary} type and let $\catCO$ be Hernandez-Leclerc's category. 
We can associate the \WS functor $\F_\ddD$ 
to a duality datum $\ddD$ in $\catCO$. 
In this paper, we introduce the notion of a strong (complete) duality datum $\ddD$ and prove that, when $\ddD$ is strong, the induced duality functor $\F_\ddD$ 
sends simple modules to simple modules and preserves the invariants $\La$, $\tLa$ and $\Li$ introduced by the authors. 
We next define the reflections $ \Refl_k$ and $ \Refl^{-1}_k$ acting on strong duality data $\ddD$.
We prove that if $\ddD$ is a strong (resp.\ complete) duality datum, then $\Refl_k(\ddD)$ and $\Refl_k^{-1}(\ddD)$ are also strong (resp.\ complete) duality data.
This allows us to make new strong (resp.\ complete) duality data by applying the reflections $\Refl_k$ and $\Refl^{-1}_k$ from known strong (resp.\ complete) duality data.
We finally introduce the notion of affine cuspidal modules in $\catCO$ by
using the duality functor $\F_\ddD$, and develop
the cuspidal module theory
for quantum affine algebras similar to
the quiver Hecke algebra case.
When $\ddD$ is complete,  we
show that all simple modules in $\catCO$ can be constructed as the heads of ordered tensor products of affine cuspidal modules.
We further prove that the ordered tensor products of affine cuspidal modules have the unitriangularity property.
This generalizes the classical simple module construction using ordered tensor products of fundamental modules.

\end{abstract}

\tableofcontents

\section{Introduction}

Let $q$ be an indeterminate and let $\catC_\g$ be the category of finite-dimensional integrable modules over  a quantum affine algebra $U_q'(\g)$.
The category $\catC_\g$ occupies an important position in the  study of quantum affine algebras because of its rich structure. 
The simple modules in $\catC_\g$ are indexed by using $n$-tuples of polynomials with constant term 1 (called \emph{Drinfeld polynomials}) (\cite{CP91,CP94,CP95} for the 
 untwisted cases and \cite{CP98} for the twisted cases). 
The simple modules can be obtained as the head of
ordered tensor product of \emph{fundamental representations} (\cite{AK97, Kas02, VV02}), and a geometric approach to simple modules was also studied in \cite{Nak01, Nak04, VV02}.

Let $\g_0$ be a finite-dimensional simple Lie algebra of $ADE$ type and let $U_q'(\g)$ be a quantum affine algebra of untwisted affine $ADE$ type.  
Hernandez and Leclerc introduced the monoidal full subcategory $\catCO$ of $\catC_\g$,  which consists of objects whose all simple subquotients are obtained from the heads of tensor products of certain fundamental representations (\cite{HL10}).
Any simple module in $\catC_\g$ can be obtained as a tensor product of suitable parameter shifts of simple modules in $\catCO$.
For each Dynkin quiver $ Q$ of $\g_0$ with a height function,  Hernandez and Leclerc introduced a monoidal subcategory $\catC_Q$ of $\catCO$ (\cite{HL15}).
 The category $\catC_Q$  is defined by using certain fundamental representations parameterized by vertices of the Auslander-Reiten quiver of $Q$.
It turns out that the complexified Grothendieck ring $\C \otimes _\Z K(\catC_Q)$ is isomorphic to the coordinate ring $\C[N]$ of the unipotent group $N$ associated with $\g_0$ and, 
under this isomorphism, the set of isomorphism classes of simple modules in $\catC_Q$ corresponds to the \emph{upper global basis} (or \emph{dual canonical basis}) of $\C[N]$ (\cite{HL15}).

In \cite{KKKO16D, KO18, OS19, OhSuh19}, the notion of the categories $\catC_\g^0$ and $\catC_Q$ is extended to all untwisted and twisted quantum affine algebras. 
Suppose that $U_q'(\g)$ is of an arbitrary affine type.  
We consider the set  $ \sig \seteq I_0 \times \cor^\times / \sim $, where the equivalence relation is given by $\eqref{eq:sig}$, with the arrows 
determined by the pole of R-matrices between tensor products of fundamental representations $V(\varpi_i)_x$ ($(i,x)\in \sig$).
Let $\sigZ$ be a connected component of $\sig$. The category $\catCO$ is defined to be the full subcategory of $\catC_\g$ determined by $\sigZ$ (see Section \ref{Sec: HL cat}).
Let $\gf$ be the simple Lie algebra of type $\rsX$ defined in $\eqref {Table: root system}$.
Note that, when $\g$ is of untwisted affine type $ADE$, $\gf$ coincides with $\g_0$. 
A $Q$-datum is a triple $\qQ \seteq (\Dynkin, \Dat, \hf)$  consisting of the Dynkin diagram $\Dynkin$ of $\gf$, an automorphism $\Dat$ on $\Dynkin$ and a height function $\hf$ (see Section \ref{Sec: Q-data}).
When $\g$ is of untwisted affine type $ADE$, $\Dat$ is the identity and $\qQ$ is equal to a Dynkin quiver with a height function.
To a $Q$-datum $\qQ$, the monoidal subcategory $\catCQ$ of $\catCO$ was introduced in \cite{HL15} for untwisted affine type $ADE$, in \cite{KKKO16D} for twisted affine type $A^{(2)}$ and $D^{(2)}$, 
in  \cite{KO18, OhSuh19} for untwisted affine type $B^{(1)}$ and $C^{(1)}$, and in \cite{OS19} for exceptional affine type.
Similarly to the untwisted affine $ADE$ case,  the category $\catCQ$ categorifies the coordinate ring $\C[N]$ of the maximal unipotent group $N$ associated with $\gf$. 
The simple Lie algebra $\gf$ is more deeply related to the structure of the category $\catC_\g$.
It is proved in \cite{KKOP20} that the simply-laced root system $\rsP$ of $\gf$ arises form $\catC_\g$ in a natural way 
and the  block decompositions of $\catC_\g$ and $\catCO$ are parameterized by the lattice associated with the root system $\rsP$.
In the course of the proof, the new invariants $\La $ and $\Li $ for $\catC_\g $ introduced in \cite{KKOP19C} are used in a crucial way.
These invariants are quantum affine algebra analogues of the invariants (with the same notations) for the quiver Hecke algebras (see \cite{KKK18A, KKKO18}).

Let $\RC$ be a \emph{quiver Hecke algebra}  (or \emph{Khovanov-Lauda-Rouquier algebras}) corresponding to a generalized Cartan matrix $\cmC$ and denote by $\RC\gmod$ its finite-dimensional graded module category. 
The algebra $\RC$ categorifies the half of the quantum group $U_q(\g)$ associated with $\cmC$ (\cite{KL09, KL11, R08}).
The simple $\RC$-modules were studied and classified by using the structure of $U_q^-(\g)$ via the categorification (see \cite{BKOP14, HMM09, Kato14, KR09, LV09,   Mc15, TW16}).
When $\RC$ is symmetric and the base field is of characteristic 0, the set of isomorphism classes of simple $\RC$-modules correspond to the upper global basis of $U_q^-(\g)$ (\cite{R11, VV09}).   
Suppose that $\cmC$ is of finite type.
One of the most successful construction for simple $\RC$-modules is the construction using \emph{cuspidal modules} via the (dual) PBW theory for $U^-_q(\g)$. 
For a reduced expression $\rxw$ of the longest element $w_0$ of the Weyl group $\weylC$, one can define the associated  cuspidal modules $\{ \Cp_k \}_{k=1, \ldots, \ell}$, which correspond to the \emph{dual PBW vectors}, 
and all simple $\RC$-modules are obtained as the simple heads of ordered tensor products of cuspidal modules.
The construction using Lyndon words was introduced in \cite{KR09} (see also \cite{HMM09}) and the construction in a general setting with a convex order was studied in \cite{Kato14, Mc15}. 
It was also studied in \cite{Kle14, Mc17} for an affine case, and in \cite{TW16} for a symmetrizable case with the viewpoint of MV polytopes.

The \emph{\WS} \cite{KKK18A} gives a connection between quiver Hecke algebras and quantum affine algebras.  
The \WS says that, for each \emph{duality datum} $\ddD = \{ \Rt_i \}_{i\in J} \subset \catC_\g$ associated with a generalized symmetric Cartan matrix $\cmC$, 
there exists a monoidal functor $\F_\ddD$, called a \emph{duality functor} shortly, from the category $R_\cmC\gmod$  to the category $\catC_\g$.
The duality functor is very interesting and useful, but it is difficult to handle it because the functor does not enjoy good properties in general. 
When $\ddD$ arises from a Q-datum, the duality functor $\F_\ddD$ enjoys good properties. 
It was shown in \cite{KKK15B, KKKO16D,KO18,OS19, Fu18} that, for each choice of Q-data $\qQ$, the \WS functor
$$
\F_\qQ \col \Rgf \gmod \longrightarrow \catCQ \subset \catC_\g^0
$$
is exact and sends simple modules to simple modules, thus it induces an isomorphism at the Grothendieck ring level.
Here $\Rgf$ is the symmetric quiver Hecke algebra associated with $\gf$. 
In this viewpoint, it is natural and important to ask which conditions for $\ddD$ provide the duality functor $\F_\ddD$ with such good properties, 
and what properties are preserved from $\RC\gmod$ to $\catC_\g$ under the duality functor $\F_\ddD$.

\medskip

This paper is a complete version of the announcement \cite{KKOP20A2}. The main results of this paper can be summarized as follows:
\bni
\item Let $U_q'(\g)$ be a quantum affine algebra of \emph{arbitrary type}.
 We find a sufficient condition for a duality datum $\ddD = \{ \Rt_i \}_{i\in J}$ to provide the functor $\F_\ddD$ with good properties. We introduce the notion of \emph{strong duality datum} by investigating \emph{root modules}.
We prove that the associated duality functor $\F_\ddD$ 
sends simple modules to simple modules and preserves the invariants $\La$, $\tLa$ and $\Li$. 
We also introduce  the notion of \emph{complete duality datum}, which can be understood as a generalization of the duality datum arising from a Q-datum.
It turns out that the Cartan matrix $\cmC$ associated with a compete duality datum $\ddD$ is equal to the one of $\gf$.

\item We introduce the \emph{reflections} $ \Refl_{  i  }$ and $ \Refl^{-1}_{  i }$ $(i\in J)$ acting on strong duality data $\ddD$.
We prove that if $\ddD$ is a strong (resp.\ complete) duality datum, then $\Refl_{  i }(\ddD)$ and $\Refl_{  i }^{-1}(\ddD)$ are also strong (resp.\ complete) duality data.
This allows us to create new strong (resp.\ complete) duality data from known strong (resp.\ complete) duality data by applying a finite sequence of the reflections $\Refl_{  i }$ and $\Refl^{-1}_{ i}$.
Indeed, the family $\st{\Refl_{  i }}_{  i  \in J}$ satisfies the \emph{braid relations}, 
etc.\ (see \cite{KKOP20A1}). It will be discussed in a forthcoming paper. 

\item We introduce the notion of \emph{affine cuspidal modules} for the category $\catCO$. Let $\ddD$ be a complete duality datum associated with a Cartan matrix $\cmC$. 
For a reduced expression $\rxw$ of the longest element of the Weyl group $\weylC$, 
we define the affine cuspidal modules $\{ \cuspS_{k} \}_{k\in \Z}$ for $\catCO$ by using the duality functor $\F_\ddD$, the right and left duals $\dual$, $\dual^{-1}$, and the cuspidal modules $\{ \Cp_k \}_{k=1, \ldots, \ell}$ of the quiver Hecke algebra $\RC$ associated with $\rxw$.
If $\ddD$ arises from a $Q$-datum, then the affine cuspidal modules $\{ \cuspS_{k} \}_{k\in \Z}$ consist of fundamental modules. But, in general, affine cuspidal modules are not fundamental.
We prove that all simple modules in $\catCO$ can be obtained uniquely as the simple heads of the ordered tensor products $\sP_{\ddD, \rxw} (\bfa)$, called \emph{standard modules}, of cuspidal modules.
We then show that the standard module  $\sP_{\ddD, \rxw} (\bfa)$ has the \emph{unitriangularity property}. 
This generalizes the classical simple module construction taking the head of ordered tensor products of fundamental representations (\cite{AK97, Kas02, Nak01, Nak04, VV02}).
The unitriangularity property allows us to define a monoidal subcategory  $ \catCab{[a,b], \ddD, \rxw}$ of $\catCO$ for an interval $[a,b]$, which is a generalization of the subcategory $\catC_l$ ($l\in \Z_{>0}$) introduced in  \cite{HL10}.
This approach can be understood as a counterpart of the PBW theory for quiver Hecke algebras via the duality functor $\F_\ddD$.
Hence we establish a base to answer the monoidal categorification conjecture for various monoidal subcategories of $\catCO$ in the same spirit of \cite{KKKO18, KKOP19C}. The monoidal categorification conjecture will be discussed in a forthcoming paper (\cite{KKOP20A3, KKOP20B}).

\ee 

We remark that, when $ U_q'(\g)$ is of untwisted affine ADE type, it has been established by Hernandez-Leclerc that the complexified Grothendieck ring $\C \otimes_\Z K(\catCO)$ can be written as a product of copies of  $\C \otimes_\Z K( \catC_Q) \simeq \C[N]$, where $N$ is the unipotent group associated with $\g_0$ (see the proof of \cite[Theorem 7.3]{HL15}).   
When the orientation of the quiver $Q$ varies, one gets various copies of $\C[N]$ in $\C \otimes_\Z K(\catCO)$  and the basis of standard modules correspond to various PBW basis.
The PBW theory developed in this paper explains this story transparently at the level of the module category.

\medskip 

Let us explain our results more precisely. Let $U_q'(\g)$ be a quantum affine algebra of an arbitrary type.
We first investigate several properties of root modules about the new invariants $\La$, $\de$, etc in Section \ref{Sec: root modules}. 
A root module is a real simple module $L$ such that
\begin{align*}
& \de(L, \dual^k(L))\ms{2mu} =\ms{3mu}\delta(k=\pm 1)
\qtext{for any $k\in\Z$.}
\end{align*}
Note that the name ``root module'' comes from Lemma~\ref{lem:root}. 
We prove several lemmas and propositions on root modules, which are used crucially in the proofs of the main results.

We next deal with the \WS. Let $\ddD $ be a duality datum associated with a generalized Cartan matrix $\cmC = (c_{i,j})_{i,j\in J}$ of symmetric type.
We study the \emph{affinizations} of modules appeared in both of two categories $\RC\gmod$ and $\catC_\g$ as pro-objects and
modify slightly the definition of \WS in order that the duality functor $\F_\ddD$ preserves the affinizations (see Theorem \ref{Thm: affinization}).
This allow us to compare the invariants $\La$, $\de$, etc between quiver Hecke algebras and quantum affine algebras via the duality functor $\F_\ddD$.
When $ \ddD= \{ \Rt_i \}_{i\in J}$ is a strong duality datum of a Cartan matrix $\cmC = (c_{i,j})_{i,j\in J}$ of simply-laced finite type (see Definition \ref{Def: SDD}), 
we prove that $\F_\ddD$ sends simple modules to simple modules (Theorem \ref{Thm: simple to simple}), i.e., $\F_\ddD$ is faithful (Corollary \ref{Cor:faithful}), and $\F_\ddD$ preserves the invariants: 
for any simple modules $M$, $N$ in $\RC\gmod$,
\bnum
\item $\La(M,N) = \La( \F_\ddD(M), \F_\ddD(N) )$,
\item $\de(M,N) = \de( \F_\ddD(M), \F_\ddD(N) )$,
\item $(\wt M,\wt N) =-\Li( \F_\ddD(M), \F_\ddD(N) )$,\label{it:wt}
\item $\de\bl \dual^k\F_\ddD(M), \F_\ddD(N)\br=0$ for any $k\not=0,\pm1$,
\item \raisebox{-.7ex}{$\ba[t]{rl}\tL(M,N) =\de\bl \dual\F_\ddD(M),\F_\ddD(N)\br =\de\bl \F_\ddD(M), \dual^{-1}\F_\ddD(N)\br\ea$}
\ee
(see Theorem \ref{Thm: invariants preserve}). The key part for the proof is to show that the invariants for \emph{determinantial modules} $\dM(w\La, \La)$ (see Section \ref{Sec: QHA})
are preserved under $\F_\ddD$, which is given in Theorem \ref{Thm: d=d}. 
Corollary \ref{Cor: mono [F]} says that the duality functor $\F_\ddD$ induces an injective ring homomorphism 
$$\KRc\monoto K(\catC_\g  ),$$
where $ \KRc$ is the specialization of the $K( \RC \gmod)$ at $q=1$.
Interestingly, the $\eps_i$ and $\eps_i^*$ in the crystal theory for $\RC\gmod$ can be interpreted in terms of the invariants $\de$ for $\catC_\g$ (see Corollary \ref{Cor: ep under F}).

Let $ \ddD= \{ \Rt_i \}_{i\in J}$ be a strong duality datum 
associated with a Cartan matrix $\cmC = (c_{i,j})_{i,j\in J}$ of simply-laced finite type, and define 
 $\catCD$ to be  the smallest full subcategory of $\catCO$ such that
\bna
\item   it contains $\F_\ddD( L )$ for any simple  $R_{\cmC}$-module $L$, 
\item  it is stable by taking subquotients, extensions, and tensor products. 
\ee
The induced map $[\F_\ddD]$ gives an isomorphism between  $K(\catCD)$ and $ \KRc$ as a ring.
We introduce the notion of \emph{unmixed pairs} of modules in $\catC_\g$ (Definition \ref{Def: unmixed}) and 
investigate several properties.  Lemma \ref{Lem: unmixed to unmixed} says that if $(M,N)$ are an unmixed pair of simple modules in $\RC\gmod$, then $(\F_\ddD(M), \F_\ddD(N))$ is strongly unmixed. 
Let $w_0$ be  the longest element
of the Weyl group $\weylC$ of $\g_\cmC$, and $\ell$ the length of $w_0$.
We define the affine cuspidal modules  $\{ \cuspS_k \}_{ k\in \Z } \subset \catCO $ to be the simple $U_q'(\g)$-modules given by 
\bna
\item $\cuspS_k = \F_\ddD(\Cp_k)$ for any $k=1, \ldots, \ell$, 
\item $\cuspS_{k+\ell} = \dual( \cuspS_k )$ for any $k\in \Z$,
\ee
where $\{\Cp_k \}_{k=1, \ldots, \ell} \subset \RC\gmod$ are  the cuspidal modules associated with $\rxw$.
Note that the cuspidal module $\Cp_k$ corresponds to the dual PBW vectors associated with $\rxw$ under the categorification using quiver Hecke algebras.
We then prove that $\cuspS_a$ is a root module for any $a\in \Z$, and $(\cuspS_a, \cuspS_b)$ is strongly unmixed for any $a > b$, which tells us that the ordered tensor product $ \cuspS_{k_1}^{\tens a_1} \tens \cdots \tens \cuspS_{k_t}^{\tens a_t} $ has a simple head 
for any decreasing integers $ k_1 > \ldots > k_t$ and $a_1, \ldots, a_t\in \Z_{\ge0}$ (see Proposition \ref{Prop: cusp}).
We next define the reflections $\Refl_k$ and $\Refl_k^{-1}$ on duality data (see $\eqref{Def: refl}$) and prove that the reflections preserve strong duality data with the same Cartan matrix (Proposition \ref{Prop: strong to strong}). 
Furthermore, we characterize simple modules in the intersections $\catC_{\Refl_i(\ddD)} \cap \catC_{\ddD}  $ and $\catC_{\Refl_i^{-1}(\ddD)} \cap \catC_{\ddD}  $ by using the cuspidal modules $\{\Cp_k \}_{k=1, \ldots, \ell}$ (see Proposition \ref{Prop: CsIDCD}).

We finally introduce the notion of a complete duality datum (see Definition \ref{Def: principal}).
We prove that, if $ \ddD= \{ \Rt_i \}_{i\in J}$ is a complete duality datum, then the associated Cartan matrix $\cmC$ has the same type as that of $\gf$.
Note that the root system $\rsP$ of $\gf$ provides the block decomposition of $\catC_\g$ (\cite{KKOP20}). 
The reflections $\Refl_k$ and $\Refl_k^{-1}$ preserve complete duality data with the same Cartan matrix (Theorem \ref{Thm: refl D}) and the duality datum $\ddDQ$ arising from a Q-datum $\qQ = (\Dynkin, \Dat, \hf)$ is complete (Proposition \ref{Prop: DQ complete}).  
By the definition, $\catC_{\ddDQ}$ is equal to $ \catC_{\qQ}$. Since a new complete duality datum can be constructed by applying the reflections to $\ddDQ$,
when $\ddD$ is complete, the category $\catCD$ can be viewed as a generalization of $\catCQ$.
We now assume that $\ddD$ is complete. Let $\{ \cuspS_k\}_{k\in \Z}$ be the affine cuspidal modules corresponding to $\ddD$ and a reduced expression $\rxw$, and set 
$
 \ZZ \seteq \Z_{\ge0}^{\oplus \Z}
$.
We denote by $\prec$ the bi-lexicographic order on $\ZZ$.
 For any $\bfa = (a_k)_{k\in \Z} \in \ZZ$, we define 
the standard module by
$$
\sP_{Q, \rxw} (\bfa) \seteq 
\cdots\tens\cuspS_1^{\otimes a_1} \otimes  \cuspS_{0}^{\otimes a_{0}}
 \otimes \cuspS_{-1}^{\otimes a_{-1}}\otimes \cdots, 
$$
and set $ \sV _{Q, \rxw} (\bfa) \seteq \hd \bl\sP_{Q, \rxw} (\bfa)\br$. 
We prove that $\sV _{Q, \rxw} (\bfa)$ is simple for any $\bfa \in \ZZ$ and the set $\{   \sV _{Q, \rxw} (\bfa) \mid \bfa \in \ZZ \}$ is a complete and irredundant set of simple modules of $\catCO$ up to isomorphisms (see Theorem \ref{Thm: main pbw}). Furthermore, Theorem \ref{Thm: main tri} says that, 
if $V$ is a simple subquotient of $\sP_{Q, \rxw} (\bfa)$ which is not isomorphic to $ \sV _{Q, \rxw} (\bfa) $, then  
$$
\bfa_{Q, \rxw}(V) \prec \bfa,
$$
which means that the module $ \sP_{Q, \rxw} (\bfa) $ has the unitriangularity property with respect to $\prec$.
For an interval $[a,b]$, we define 
$ \catCab{[a,b], \ddD, \rxw}$ to be the full subcategory of $\catC_\g$ whose objects have all their composition factors $V$ satisfying 
$$
b \ge l( \bfa_{\ddD, \rxw}(V) ) \quad \text{ and }\quad r( \bfa_{\ddD, \rxw}(V) )\ge a,
$$
where $l$ and $r$ are defined in $\eqref{Eq: def lr}$.
By the unitriangularity, the category $ \catCab{[a,b], \ddD, \rxw}$ is stable by taking tensor products, and 
it also enjoys the same properties (Theorem \ref{Thm: unitry for C[a,b]}).

\smallskip

This paper is organized as follows.
In Section \ref{Sec: preliminaries},  we give the necessary background on quiver Hecke algebras, quantum affine algebras, and the invariants related to R-matrices.
In Section \ref{Sec: root modules}, we introduce the notion of root modules and investigate several properties.
In Section \ref{Sec: SW}, we study affinizations and the duality functor $\F_\ddD$, and prove that, when $\ddD$ is strong,   
$\F_\ddD$ sends simple modules to simple modules and preserves the new invariants.   
In Section \ref{Sec: SD and ACM}, we introduce the notions of affine cuspidal modules and reflections, and prove that the reflections preserve the strong duality data. 
In Section \ref{Sec: PBW}, we study the PBW theoretic approach to $\catCO$ using a complete duality datum and affine cuspidal modules.

\medskip

\vskip 1em 

{\bf Acknowledgments}

The second, third and fourth authors gratefully acknowledge for the hospitality of RIMS (Kyoto University) during their visit in 2020.
The authors would like to thank the anonymous referee for valuable comments and suggestions.

\vskip 2em

\section{Preliminaries} \label{Sec: preliminaries}

\begin{convention}
\ 
\bnum
\item
For a statement $P$, $\delta(P)$ is $1$ or $0$ according that
$P$ is true or not.
\item For a field $ \cor$, $a\in\cor$ and $f(z)\in\cor(z)$, 
we denote
by $\zero_{z=a}f(z)$ the order of zero of $f(z)$ at $z=a$.\label{conv:zero}
\item For a ring $A$, $A^\times$ is the set of invertible elements of $A$. 
\ee
\end{convention}

\subsection{Quantum groups}\

Let $I $ be an index set. A quintuple $ (\cmA,\wlP,\Pi,\wlP^\vee,\Pi^\vee) $ is called a  (symmetrizable) {\it Cartan datum} if it
consists of
\begin{enumerate}
\item[(a)] a generalized \emph{Cartan matrix} $\cmA=(a_{ij})_{i,j\in I}$, 
\item[(b)] a free abelian group $\wlP$, called the {\em weight lattice},
\item[(c)] $\Pi = \{ \alpha_i \mid i\in I \} \subset \wlP$,
called the set of {\em simple roots},
\item[(d)] $\wlP^{\vee}=
\Hom_{\Z}( \wlP, \Z )$, called the \emph{coweight lattice},
\item[(e)] $\Pi^{\vee} =\{ h_i \in \wlP^\vee \mid i\in I\}$, called the set of {\em simple coroots} 
\end{enumerate}
satisfying the following:
\bnum
\item $\lan h_i, \alpha_j \ran = a_{ij}$ for $i,j \in I$,
\item $\Pi$ is linearly independent over $\Q$,
\item for each $i\in I$, there exists $\Lambda_i \in \wlP$, called the \emph{fundamental weight}, such that $\lan h_j,\Lambda_i \ran =\delta_{j,i}$ for all $j \in I$.
\item there is a symmetric bilinear 
form $( \cdot \, , \cdot )$ on $\wlP$ satisfying 
\begin{equation*}
(\alpha_i,\alpha_i)\in 2\Z_{>0}\qtq \lan h_i,  \lambda\ran = \dfrac{2 (\alpha_i,\lambda)}{(\alpha_i,\alpha_i)}.
\end{equation*} 
\end{enumerate}

We set $ \rlQ \seteq \bigoplus_{i \in I} \Z \alpha_i$ and  $\rlQ^+ \seteq  \sum_{i\in I} \Z_{\ge 0} \alpha_i$ and define $\height{\beta}=\sum_{i \in I} k_i $  for $\beta = \sum_{i \in I} k_i \alpha_i \in \rlQ^+$.
We define $\wlP^+ \seteq  \{ \La \in \wlP \mid \langle h_i , \La\rangle \in \Z_{\ge0} \text{ for any } i\in I \}$.

We write $\prD$ for the set of  positive roots associated with $\cmA$ and set $\nrD \seteq  - \prD$.
Denote by $\weyl$  the \emph{Weyl group}, which is the subgroup of $\mathrm{Aut}(\wlP)$ generated by  
$ s_i(\lambda) \seteq  \lambda - \langle h_i, \lambda \rangle \alpha_i  $ for $i\in I$. 

We denote  by $U_q(\g)$ the \emph{quantum group} associated with $(\cmA, \wlP,\wlP^\vee, \Pi, \Pi^{\vee})$, which
is a $\Q(q)$-algebra generated by $f_i$, $e_i$ $(i\in I)$ and $q^h$ $(h\in \wlP^\vee)$ with certain defining relations (see \cite[Chapter 3]{HK02} for details).
We denote by $U_q^+(\g)$ (resp.\ $U_q^-(\g)$) the subalgebra of $U_q(\g)$ generated by $e_i$'s (resp.\ $f_i$'s).
Set $\A \seteq\Z[q,q^{-1}]$ and write $U^{\pm}_\A(\g)$ for  the $\A$-lattice  of $U_q^\pm(\g)$, which is the $\A$-subalgebra generated by $e_i^{(n)}$ (resp.\ $f_i^{(n)}$) for $i\in I$ and $n \in \Z_{\ge 0}$.  
We write the \emph{unipotent quantum coordinate ring}
$$
A_q(\n) \seteq  \bigoplus_{\beta \in \rlQ^-} A_q(\n)_\beta\qt{where 
$A_q(\n)_\beta \seteq  \Hom_{\Q(\q)}( U_q^+(\g)_{ -\beta }, \Q(q))$,}
$$
and denote by $ A_q(\n)_\A $ the $\A$-lattice of $A_q(\n)$.
Note that $A_q(\n)$ is isomorphic to $U_q^-(\g)$ as a $\Q(q)$-algebra \cite[Lemma 8.2.2]{KKKO18}.

\subsection{Quiver Hecke algebras} \label{Sec: QHA} \ 

Let $\mathbf k$ be a field and let $(\cmA,\wlP, \Pi,\wlP^{\vee},\Pi^{\vee})$ be a Cartan datum.
 Choose polynomials 
 $$
Q_{i,j}(u,v) =
\delta (i\not=j)  \sum\limits_{ \substack{ (p,q)\in \Z^2_{\ge0} 
\\ (\alpha_i , \alpha_i)p+(\alpha_j , \alpha_j)q=-2(\alpha_i , \alpha_j)}}
t_{i,j;p,q} u^p v^q \in \cor [u,v]
$$
with $t_{i,j;p,q}\in{\mathbf k}$, $t_{i,j;p,q}=t_{j,i;q,p}$ and $t_{i,j;-a_{ij},0} \in {\mathbf k}^{\times}$.
Note that $Q_{i,j}(u,v)=Q_{j,i}(v,u)$ for $i,j\in I$. 
Let
${\sg}_{n} = \langle s_1, \ldots, s_{n-1} \rangle$ be the symmetric group
on $n$ letters with the action of  ${\mathfrak{S}}_n$ on $I^n$ by place permutation.
For  $\beta \in \rlQ^+$ with $\height{\beta} = n$, we set
$$I^{\beta} \seteq  \{\nu = (\nu_1, \ldots, \nu_n) \in I^{n}\mid \alpha_{\nu_1} + \cdots + \alpha_{\nu_n} = \beta\}.$$

\begin{definition}
Let $\beta \in \rlQ^+$ with $\height{\beta}=n$. The \emph{quiver Hecke algebra}  
$R(\beta)$  associated
with the parameters $\{ Q_{i,j} \}_{i,j \in I}$ is the ${\mathbf k}$-algebra generated by $\{ e(\nu) \}_{\nu \in  I^{\beta}}$, $ \{x_k \}_{1 \le
k \le n}$, $\{ \tau_m \}_{1 \le m \le n-1}$ satisfying the following defining relations:
\begin{align*} 
& e(\nu) e(\nu') = \delta_{\nu, \nu'} e(\nu), \ \
\sum_{\nu \in  I^{\beta} } e(\nu) = 1, \allowdisplaybreaks\\
& x_{k} x_{m} = x_{m} x_{k}, \ \ x_{k} e(\nu) = e(\nu) x_{k}, \allowdisplaybreaks\\
& \tau_{m} e(\nu) = e(s_{m}(\nu)) \tau_{m}, \ \ \tau_{k} \tau_{m} =
\tau_{m} \tau_{k} \ \ \text{if} \ |k-m|>1, \allowdisplaybreaks\\
& \tau_{k}^2 e(\nu) = Q_{\nu_{k}, \nu_{k+1}} (x_{k}, x_{k+1})
e(\nu), \allowdisplaybreaks\\
& (\tau_{k} x_{m} - x_{s_k(m)} \tau_{k}) e(\nu) = \begin{cases}
-e(\nu) \ \ & \text{if $m=k$, $\nu_{k} = \nu_{k+1}$,} \\
e(\nu) \ \ & \text{if $m=k+1$, $\nu_{k}=\nu_{k+1}$,} \\
0 \ \ & \text{otherwise,}
\end{cases} \allowdisplaybreaks\\
& (\tau_{k+1} \tau_{k} \tau_{k+1}-\tau_{k} \tau_{k+1} \tau_{k}) e(\nu)\\
& =\begin{cases} \dfrac{Q_{\nu_{k}, \nu_{k+1}}(x_{k},
x_{k+1}) - Q_{\nu_{k}, \nu_{k+1}}(x_{k+2}, x_{k+1})} {x_{k} -
x_{k+2}}e(\nu) \ \ & \text{if} \
\nu_{k} = \nu_{k+2}, \\
0 \ \ & \text{otherwise}.
\end{cases}
\end{align*}
\end{definition}

The algebra $R(\beta)$ has the  $\Z$-grading defined by 
\begin{equation*} \label{eq:Z-grading}
\deg e(\nu) =0, \quad \deg\, x_{k} e(\nu) = (\alpha_{\nu_k}
, \alpha_{\nu_k}), \quad\deg\, \tau_{l} e(\nu) = -
(\alpha_{\nu_l} , \alpha_{\nu_{l+1}}).
\end{equation*}
For a $\Z$-graded $\cor$-algebra $A$, we denote by $A \Mod$ the category of  graded left $A$-modules, 
 and write $A \gproj$ (resp.\ $A \gmod$) for the full subcategory of $A \Mod$ consisting of
finitely generated projective (resp.\ finite-dimensional) graded $A$-modules. 
We set $R\gproj \seteq  \bigoplus_{\beta \in \rlQ^+}R(\beta) \gproj$ and 
$R\gmod \seteq  \bigoplus_{\beta \in \rlQ^+}R(\beta) \gmod$.

 For $M\in R(\beta)\Mod$ and $N\in R(\gamma)\Mod$, we define their convolution product by 
$$
M\conv N \seteq R(\beta+\gamma)e(\beta,\gamma) \otimes_{R(\beta) \otimes R(\gamma)} (M  \otimes  N),
$$
where $e(\beta,\gamma) =\displaystyle\sum_{ \nu_1 \in I^{\beta},  \nu_2 \in I^{\gamma} } e(\nu_1*\nu_2)$. Here $\nu_1*\nu_2$ is 
the concatenation of $\nu_1$ and $\nu_2$. 
We denote by $M \hconv N$ the head of $M\conv N$ and by $M\sconv N$ the socle of $M \conv N$. We say that simple $R$-modules $M$ and $N$ \emph{strongly commute} if $M \conv N$ is simple. 
A simple $R$-module $L$ is \emph{real} if $L \conv L $ is simple. For $i\in I$ and an $R(\beta)$-module $M$, we define 
$$
E_i(M) \seteq  e(\al_i, \beta-\al_i) M, \quad F_i(M) \seteq  R(\al_i)\conv M,
$$
and 
\begin{align*}
\wt(M) &\seteq -\beta, \\ 
\eps_i (M) &\seteq  \max\{ k \ge0 \mid E_i^k(M) \ne 0 \},  \\
 \vphi_i(M) &\seteq  \eps_i(M) + \langle h_i, \wt(M) \rangle.
\end{align*}
For $i\in I$, we denote by $L(i)$ the self-dual 1-dimensional  simple $R(\al_i)$-module.
For a simple module $M$, $\tf_i(M)$ (resp.\ $\te_i(M)$) is the self-dual simple $R$-module being isomorphic to $L(i) \hconv M$ (resp.\ $\mathrm{soc}(E_i M)$). 
One also define $E_i^*$, $F_i^*$, $\eps_i^*$, etc in the same manner as above if replacing the role of $e(\al_i, \beta-\al_i)$ and $R(\al_i)\conv -$ with the ones of $e( \beta-\al_i, \al_i)$ and $-\conv R(\al_i)$.

\begin{theorem}[\cite{KL09, KL11, R08}]\label{Thm: categorification}
There exist  $\A$-bialgebra isomorphisms 
\begin{eqnarray*}
 U_\A^-(\g) &\buildrel \sim \over \longrightarrow& K(R \gproj) \quad \text{ and } \quad A_q(\n)_\A \buildrel \sim \over \longrightarrow K(R \gmod),
\end{eqnarray*}
where $K(R\gproj)$ and $K(R\gmod)$ are the Grothendieck groups of $R \gproj$ and $R\gmod$.
\end{theorem}

\smallskip

\begin{definition}
The quiver Hecke algebra $R(\beta)$ is said to be \emph{symmetric} if $ Q_{i,j}(u,v)$ is a polynomial in $u-v$ for any
$i,j \in I$.
\end{definition}

When $R$ is symmetric, the Cartan matrix $\cmA$ is of symmetric type. In this case we assume that $(\al_i,\al_i)=2$ for all $i\in I$.

{\em In the sequel, we assume that $R$ is symmetric.} 

Let $z$ be an indeterminate with homogeneous degree $2$.  
For an  $R(\beta)$-module $M$, 
we denote by $M^\aff $  the \emph{affinization} of $M$ (see \cite{KKK18A, KP18}).
When $R(\beta)$ is symmetric,  $M^\aff = \cor[z]\otimes_\cor M$ and the $R(\beta)$-module structure of $M^\aff$
is defined by
\begin{align*}
e(\nu) (f \otimes m) &= f \otimes e(\nu)m, \\
 x_j (f \otimes m) &= (zf) \otimes m + f \otimes x_j m, \\
  \tau_k(f \otimes m) &= f \otimes (\tau_k m)
\end{align*}
for $f \in \cor[z]$, $m \in M$, $\nu \in I^\beta$ and admissible $j,k$.
We sometimes write $M_z$ instead of $M^\aff$ to emphasize $z$.

Let $\beta\in\rlQ^+$ and $m=\height{\beta}$.
For $k = 1, \ldots, m-1$ and $\nu\in I^\beta$, the {\em intertwiner} $\varphi_k \in R(\beta)$ is defined by
$$ \varphi_k e(\nu) \seteq \left\{
              \begin{array}{ll}
                (\tau_kx_k - x_k\tau_k) e(\nu) & \hbox{ if } \nu_k = \nu_{k+1}, \\
                \tau_k e(\nu) & \hbox{ otherwise.}
              \end{array}
            \right.
$$
Note that $\{\vphi_k\}_{1\le k\le m-1}$ satisfies the braid relation.
Hence, we can define $\vphi_w$ for any $w\in\sym_m$. 
Let $M$ be an $R(\beta)$-module with $\height{\beta}=m$ and $N$ an $R(\beta')$-module with $ \height{\beta'}=n$.
Let $w[n,m]$ be the element of $\sym_{m+n}$ which sends
$k\mapsto k+m$ for $1\le k\le n$ and $k\mapsto k-n$ if $n<k\le m+n$. 
Then the $R(\beta)\otimes R(\beta')$-linear map $M\otimes N \longrightarrow N \conv M$ defined by $u\otimes v \mapsto \varphi_{w[n,m]}(v \otimes u)$ can be extended to the $R(\beta+\beta')$-module
homomorphism (up to a grading shift)
\begin{align*}
R_{M,N}\cl M\conv N \longrightarrow N \conv M.
\end{align*}
For non-zero $R$-modules $M$ and $N$, we set
$$
\Rm_{M_{z},N_{z'}} \seteq (z' - z)^{-s} R_{M_{z},N_{z'}}:  M_z \conv N_{z'} \to N_{z'}\conv M_z,
$$
where $s$ is the largest integer such that $R_{M_z,N_{z'}}(M_z \conv N_{z'})\subset (z'-z)^s N_{z'}\conv M_z$.
We call it the {\em renormalized $R$-matrix}. 
Then, we define 
$$
\rmat{M,N} \cl M\conv N\to N\conv M
$$
as the specialization of $\Rm_{M_z,N_{z'}}$ at $z=z'=0$ 
(up to a constant multiple), 
which never vanishes by the definition (see \cite[Section 1]{KKK18A} and \cite[Section 2]{KP18} for details).

\begin{definition}
Let $M$ and $N$ be simple $R$-modules. We set
\begin{align*}
\La(M,N) &\seteq \deg (\rmat{M,N}), \\
\tLa(M,N)&\seteq\dfrac{1}{2}\bl \La(M,N)+  \bl \wt(M), \wt(N)   \br  \br, \\
\de(M,N) &\seteq \dfrac{1}{2}\bl\La(M,N) + \La(N,M)\br.
\end{align*}
\end{definition}
Many properties of $\La$, $\tLa$, and $\de$ were obtained in \cite{KKKO18, KKOP18, KKOP19A}. 

\smallskip

We now define the monoidal subcategories $\cC_w$, $\cC_{*,v}$ and $\cC_{w,v}$ of $R\gmod$ for $w,v \in \weyl$. 
For $M \in  R(\beta)\Mod$, we define
\begin{align*}
\gW(M) &\seteq  \{  \gamma \in  \rlQ^+ \cap (\beta - \rlQ^+)  \mid  e(\gamma, \beta-\gamma) M \ne 0  \}, \\
\sgW(M) &\seteq  \{  \gamma \in  \rlQ^+ \cap (\beta - \rlQ^+)  \mid  e(\beta-\gamma, \gamma) M \ne 0  \}.
\end{align*}
For $w\in \weyl$, we denote by $\cC_{w}$  the  full subcategory of $R\gmod$ whose objects $M$ satisfy
\begin{align*}
\gW(M) \subset \Sp( \prD \cap w \nrD ).
\end{align*}
Similarly, for $v\in \weyl$, we define $\cC_{*,v}$ to be the  full subcategory of $R\gmod$ whose objects $N$ satisfy
\begin{align*}
\sgW(N) \subset \Sp ( \prD \cap v \prD ).
\end{align*}
Finally, we define $\cC_{w,v}\seteq\cC_w\cap\cC_{*,v}$.

 When $\g$ is of finite type, we have $\cC_{w_0} = R\gmod$ and
\begin{align*}
&\text{ $M \in \cC_{s_i w_0} $ if and only if $\eps_i(M)=0$,}  \\
&\text{ $M \in \cC_{*, s_i} $ \  if and only if $\eps_i^*(M)=0$,} 
\end{align*}
for any  $R$-module $M$ in $\R\gmod$ and $i\in I$.  Here $w_0$ denotes the longest element of $\weyl$ (see \cite{KKOP18} for details).

Let $\underline{w} \seteq  s_{i_1}\cdots s_{i_l}$ be a reduced expression of $w\in \weyl$ and define 
\begin{align} \label{Eq: beta_k}
\beta_k \seteq  s_{i_1} \cdots s_{i_{k-1}} ( \al_{i_k} )\quad \text{ for $k=1, \ldots, l$.}
\end{align}
Then we have  $ \prD \cap w \nrD = \{ \beta_1, \ldots, \beta_l \}$ with the \emph{convex order} $\prec$ on $\prD \cap w \nrD$, i.e.,  $\beta_a \prec\beta_b$ for any $a < b$.
For $\beta \in \prD \cap w \nrD$, a pair $(\al, \gamma)$ is called a \emph{minimal pair} of $\beta$ if $\beta = \al+\gamma$, $ \al \prec\gamma$ and there exists no pair $(\al', \gamma')$ such that $ \beta = \al' + \gamma'$ and $\al \prec \al' \prec \gamma' \prec \gamma$.
The convex order provides the \emph{PBW vectors} $\{ E(\beta_{k}) \}_{k=1, \ldots, l} $ in $U_\A^-(\g)$ and the \emph{dual PBW vectors} $\{ E^*(\beta_{k}) \}_{k=1, \ldots, l} $ in $A_{q}(\n)_{\A}$.
We set $A_q(\n(w))$ to be the subalgebra of $A_q(\n)$ generated by $  E^*(\beta_{k})$ for $k=1, \ldots, l$.
The category $\cC_w$ categorifies the algebra $A_q(\n(w))$ (\cite{KKKO18, KKOP18}).

For $k=1, \ldots,l $,  let $\Cp_k$ be the \emph{cuspidal module} corresponding to $\beta_k$ with respect to $\underline{w}$ (see \cite[Section 2]{KKOP18} for precise definition). 
Under the categorification, the cuspidal module $\Cp_k$ corresponds to the dual PBW vector $E^*(\beta_k)$.
It is known that the set
$$
\{  \hd ( \Cp_{l}^{\circ a_l} \conv \cdots \conv \Cp_1^{\circ a_1} ) \mid (a_1, \ldots, a_l) \in \Z_{\ge0}^l  \}
$$
gives a complete set of pairwise non-isomorphic simple graded modules in $\cC_w$, up to a grading shift (see \cite{BKM14, Kato14, Mc15, TW16}).
Note that, for a minimal pair $(\beta_a, \beta_b)$ of $\beta_k$, there exists an isomorphism
\begin{align} \label{Eq: minimal for QHA}
\Cp_a \hconv \Cp_b \simeq \Cp_k
\end{align}
(see \cite[Lemma 4.2]{Mc15} and \cite[Section 4.3]{BKM14}.)

For $\La \in \wlP^+$ and $w,v \in \weyl $ with $w \ge v$, we denote by 
$\dM(w\La, v\La)$ the \emph{determinantial module} in $R\gmod$ corresponding to the pair $(w\La, v\La)$ (see \cite[Section 10.2]{KKKO18} and \cite[Section 4]{KKOP18} for precise definition).
Under the categorification, the determinantial module $\dM(w\La, v\La)$ corresponds to the \emph{unipotent quantum minor} $D(w\La, v\La)$ in $A_q(\n)$ (\cite[Proposition 4.1]{KKOP18}).

{}From now on, {\em we assume that $\cor$ is a field of characteristic $0$ and that $R$ is a symmetric quiver Hecke algebra of finite {\rm ADE}\;type.}

Note that, under the categorification by $R\gmod$, 
the \emph{upper global basis} (or \emph{dual canonical basis}) of $A_q(\n)$ corresponds to the set of isomorphism classes of simple $R$-modules \cite{R11, VV09}. 
Then the reflection functor $\rT_i$ constructed in \cite{Kato14} gives an equivalence of categories:
\begin{align*}
\rT_i \cl \cC_{s_i w_0} \buildrel \sim \over \longrightarrow  \cC_{*, s_i}. 
\end{align*}
Note that $\rT_i$ is  denoted by $\rT_i^*$ in \cite{KKOP18}. 
Since, at the crystal level, this functor corresponds to the Saito crystal reflection (\cite{Saito94}), we have 
\begin{align} \label{Eq: Saito refl}
 \rT_i(M) \simeq \tf_i^{ \vphi_i^*(M)}  \te_i^{*\hskip 0.1em \eps_i^*(M)} (M)
\end{align}
for a simple module $ M$ with $\eps_i(M)=0$. 
For a reduced expression $\underline{w} \seteq  s_{i_1}\cdots s_{i_l}$, the cuspidal module $\Cp_k$ can be computed as follows (see \cite[Section 5]{KKOP18}): 
\begin{align} \label{Eq: cusp and refl}
\Cp_k \simeq  \rF_{i_1} \rF_{i_2} \cdots \rF_{i_{k-1}} (L( i_k)  )
\end{align}
for $k=1, \ldots, l$.

\subsection{Quantum affine algebras}\

We assume that $\cmA=(a_{i,j})_{i,j\in I}$
is an affine Cartan matrix. Note that the rank of $\wlP$ is $|I|+1$.
We denote by $\delta \in \rlQ$  the \emph{imaginary root} and by $c $ 
the  \emph{central element} in $\wlP^\vee$. Note that 
the positive  imaginary root $\Delta_+^{\rm im}$ is equal to $\Z_{>0} \delta$  and the center of  $\g$ is generated by $c$.
We write $\wlPc  \seteq  \wlP / (\wlP \cap \Q \delta) $, called the \emph{classical weight lattice}, and take $\rho \in \wlP$ (resp.\ $\rho^\vee \in \wlP^\vee$) such that $\lan h_i,\rho \ran=1$ (resp.\ $\lan \rho^\vee,\al_i\ran =1$) for any $i \in I$.
We choose a  $\Q$-valued non-degenerate  symmetric bilinear form $( \ , \ )$ on $\wlP$ satisfying 
$$
\lan h_i,\la \ran= \dfrac{2(\al_i,\la)}{(\al_i,\al_i)} \quad \text{ and} \quad \lan c,\la \ran = (\delta,\la)
$$
for any $i \in I$ and $\la \in \wlP$.
We define $\g$ to be the \emph{affine Kac-Moody algebra} associated with $\cmA$. 
We shall use the standard convention in~\cite{Kac}
to choose $0\in I$ except $A^{(2)}_{2n}$ type, in which we take the longest simple root as $\al_0$, and 
 $B_2^{(1)}$ and $A_3^{(2)}$ types, in which we take the following Dynkin diagrams:
\begin{equation*} 
\begin{aligned} 
& A^{(2)}_{2n} :   \xymatrix@C=4ex@R=3ex{
  *{\circ}<3pt> \ar@{<=}[r]_<{n \ } & *{ \circ }<3pt> \ar@{-}[r]_<{n-1}  & *{ \circ }<3pt> \ar@{-}[r]_<{n-2} & \cdots \ar@{-}[r]_<{ }   &*{\circ}<3pt> \ar@{-}[l]^<{1} &*{ \circ }<3pt>  \ar@{=>}[l]^<{ \ \ 0}  } \hs{3ex}
B^{(1)}_{2} :   \xymatrix@C=4ex@R=3ex{
  *{\circ}<3pt> \ar@{=>}[r]_<{0 \ } & *{ \circ }<3pt>  \ar@{<=}[l]^<{2} &*{ \circ }<3pt>  \ar@{=>}[l]^<{ \ \ 1}  } \hs{3ex}
A^{(2)}_{3} :   \xymatrix@C=4ex@R=3ex{
  *{\circ}<3pt> \ar@{<=}[r]_<{0 \ } & *{ \circ }<3pt>  \ar@{=>}[l]^<{2} &*{ \circ }<3pt>  \ar@{<=}[l]^<{ \ \ 1}  }
\end{aligned}
\end{equation*}
Note that $B^{(1)}_2$ and $A^{(2)}_3$ in the above diagram 
are denoted by $C_2^{(1)}$ and $D_3^{(2)}$ respectively in \cite{Kac}.

Set $ I_0 \seteq I \setminus \{ 0 \}$.

Let $q$ be an indeterminate and $\ko$ the algebraic closure of the subfield $\C(q)$
in the algebraically closed field $\corh\seteq\bigcup_{m >0}\C((q^{1/m}))$. 
 For $m,n \in \Z_{\ge 0}$ and $i\in I$, we define
$q_i = q^{(\alpha_i,\alpha_i)/2}$ and
\begin{equation*}
 \begin{aligned}
 \ &[n]_i =\frac{ q^n_{i} - q^{-n}_{i} }{ q_{i} - q^{-1}_{i} },
 \ &[n]_i! = \prod^{n}_{k=1} [k]_i ,
 \ &\left[\begin{matrix}m \\ n\\ \end{matrix} \right]_i=  \frac{ [m]_i! }{[m-n]_i! [n]_i! }.
 \end{aligned}
\end{equation*}

Let $d$ be the smallest positive integer such that 
$d \frac{(\al_i, \al_i)}{2}\in\Z$  for all $i\in I$.

\begin{definition} \label{Def: GKM}
The {\em quantum affine algebra} $U_q(\g)$ associated with an affine Cartan datum $(\cmA,\wl,\Pi,\wl^\vee,\Pi^\vee)$ is the associative algebra over $\ko$ with $1$ generated by $e_i,f_i$ $(i \in I)$ and
$q^{h}$ $(h \in  d^{-1} \wl^{\vee})$ satisfying following relations:
\bnum
\item  $q^0=1, q^{h} q^{h'}=q^{h+h'} $\quad for $ h,h' \in d^{-1} \wl^{\vee},$
\item  $q^{h}e_i q^{-h}= q^{\langle h, \alpha_i \rangle} e_i$,
$q^{h}f_i q^{-h} = q^{-\langle h, \alpha_i \rangle }f_i$\quad for $h \in d^{-1}\wl^{\vee}, i \in I$,
\item  $e_if_j - f_je_i =  \delta_{ij} \dfrac{K_i -K^{-1}_i}{q_i- q^{-1}_i }, \ \ \text{ where } K_i=q_i^{ h_i},$
\item  $\displaystyle \sum^{1-a_{ij}}_{k=0}
(-1)^ke^{(1-a_{ij}-k)}_i e_j e^{(k)}_i =  \sum^{1-a_{ij}}_{k=0} (-1)^k
f^{(1-a_{ij}-k)}_i f_jf^{(k)}_i=0 \quad \text{ for }  i \ne j, $
\ee
where $e_i^{(k)}=e_i^k/[k]_i!$ and $f_i^{(k)}=f_i^k/[k]_i!$.
\end{definition}

Let us denote by $U_q'(\g)$ the $\cor$-subalgebra of $U_q(\g)$ generated by $e_i,f_i,K^{\pm 1}_i$ $(i \in I)$.
Let $\catC_\g$ be the category of finite-dimensional integrable
$\uqpg$-modules, 
i.e., finite-dimensional modules $M$ with a weight decomposition
$$
M=\soplus_{\la\in\wlPc}M_\la  \qquad \text{ where } M_\la=\st{u\in M\mid K_iu=q_i^{\langle h_i,\la \rangle}  u  }.
$$
Note that   the trivial  module $\trivial$ is contained in $\catC_\g$ and
 the tensor product $\otimes$ gives a monoidal category structure on $\catC_\g$. 
The monoidal category $\catC_\g$ is rigid.  For $M\in \catC_\g$, we denote by  $\dual M$ and $\dual^{-1} M$ the right and the left dual of $M$, respectively.
Hence we have the evaluation morphisms
$$M\tens \dual M\to\one\quad \text{and} \quad  \dual^{-1} M\tens M\to\one.$$
We extend this to $\dual^k$ for $k \in \Z$.
We set $M^{\otimes k} \seteq \underbrace{M \tens \cdots \tens M}_{k\text{-times}}$
for $k\in\Z_{\ge0}$. 
For $M,N\in \catC_\g$, we denote by $M \htens N$ the head of $M\tens N$ and by $M\stens N$ the socle of $M \tens N$. We say that $M$ and $N$ \emph{strongly commute} if $M  \tens  N$ is simple. 
A simple $U_q'(\g)$-module $L$ is \emph{real} if $L \tens L $ is simple.

A simple module $L$ in $\catC_\g$ contains a non-zero vector $u \in L$ of weight $\lambda\in \wlPc$ such that (i) $\langle h_i,\lambda \rangle \ge 0$ for all $i \in I_0$,
(ii) all the weight of $L$ are contained in $\lambda - \sum_{i \in I_0} \Z_{\ge 0} \clp (\alpha_i)$, where $\clp \colon \wl\to \wlPc$ is the canonical projection.
Such a $\la$ is unique and $u$ is unique up to a constant multiple. We call $\lambda$ the {\it dominant extremal weight} of $L$ and $u$ a {\it dominant extremal weight vector} of $L$.
 For each $i \in I_0$, we set 
$$
\varpi_i \seteq {\rm gcd}(\mathsf{c}_0,\mathsf{c}_i)^{-1}\clp (\mathsf{c}_0\Lambda_i-\mathsf{c}_i \Lambda_0),
$$
where the central element $c$ is equal to $ \sum_{ i\in I} \mathsf{c}_i h_i$.
For any $i\in I_0$, we denote by $V(\varpi_i)$ the \emph{$i$-th fundamental representation}. 
Note that the dominant extremal weight of $V(\varpi_i)$ is $\varpi_i$.

\subsection{R-matrices} \label{subsec:R}

In this subsection we review the notion of R-matrices on $U_q'(\g)$-modules and their coefficients (see \cite{D86}, \cite[Appendices A and B]{AK97} and \cite[Section 8]{Kas02} for details).

For a module $M\in \catC_\g$, we denote by $M^\aff$ the {\it affinization} of $M$ and 
by $z_M \colon M^\aff \to M^\aff$ the $U_q'(\g)$-module automorphism of weight $\delta$. Note that $M^\aff \simeq \cor[z^{\pm 1}]\otimes_{ \cor } M $ 
with the action
$$
e_i(a\tens v)=z^{\delta_{i,0}}a\tens e_iv\quad \text{for $a\in\cor[z^{\pm1}]$ and $v\in M$.}
$$ 
We sometimes write $M_z$ instead of $M^\aff$ to emphasize the endomorphism $z$. 
 For $x \in \ko^\times$, we define
$$M_x \seteq M^\aff / (z_M -x)M^\aff.$$
We call $x$ a {\it spectral parameter} (see \cite[Section 4.2]{Kas02} for details).

Take a basis $\{P_\nu\}_\nu$ of $U_q^+(\g)$ and a basis $\{Q_\nu\}_\nu$ of $U_q^-(\g)$ dual to each other with respect to a suitable coupling between
$U_q^+(\g)$ and $U_q^-(\g)$. For $U_q'(\g)$-modules $M$ and $N$, we define
\begin{align*}
\Runiv_{M,N}(u\otimes v) \seteq q^{(\wt(u),\wt(v))} \sum_\nu P_\nu v\otimes Q_\nu u \quad  \text{  for $u\in M$ and $v\in N$,}
\end{align*}
so that
$\Runiv_{M,N}$ gives a $\uqpg$-linear homomorphism $M\otimes N \rightarrow N\otimes M$, called the \emph{universal R-matrix}, 
provided that the infinite sum has a meaning.
As $\Runiv_{M,N_z}$ converges in the $z$-adic topology for 
$ M,N\in \catC_\g$, we have
a morphism of
$\ko((z))\tens\uqpg$-modules
\begin{align*} 
\Runiv_{M,N_z} \colon \ko((z))\tens_{\ko[z^{\pm1}]} (M \tens N_z) \To \ko((z))\tens_{\ko[z^{\pm1}]} (N_z\tens M).
\end{align*}
Note that  $\Runiv_{M,N_z}$ is an isomorphism.
For  non-zero modules $M, N \in \catC_\g$, 
we say that the universal R-matrix  $\Runiv_{M,N_z}$ is \emph{rationally renormalizable}
if there exists $f(z) \in \ko((z))^\times$ 
such that
$$f(z) \Runiv_{M,N_z}\bl M\tens N_z\br\subset N_z\tens M. $$
In this case, we can choose
$c_{M,N}(z) \in \ko((z))^\times$ such that,
 for any $x \in \ko^\times$, the specialization of $\Rren_{M,N_z} \seteq c_{M,N}(z)\Runiv_{M,N_z}
\col M \otimes N_z \to N_z \otimes M$ at $z=x$
$$  \Rren_{M,N_z}\big\vert_{z=x} \colon M \otimes N_x  \to N_x \otimes M$$
does not vanish.  Note that $\Rren_{M,N_z}$ and $c_{M,N}(z)$ are unique up to a multiple of $\cz^\times = \bigsqcup_{\,n \in \Z}\ko^\times z^n$. 
We call $\Rren_{M,N_z}$ {\em the renormalized R-matrix} and 
$c_{M,N}(z)$  the \emph{renormalizing coefficient}.
We denote by $\rmat{M,N}$ the specialization at $z=1$ 
\eq
\rmat{M,N} \seteq \Rren_{M,N_z}\vert_{z=1} \colon M \otimes N \to N \otimes M,
\label{eq:rmat}
\eneq
and call it the \emph{R-matrix}. The  R-matrix $\rmat{M,N}$ is well-defined up to a constant multiple whenever $\Runiv_{M,N_z}$ is \rr. By the definition, $\rmat{M,N}$ never vanishes.

Let $M$ and $N$ be simple modules in $\catC_\g$ and let $u$ and $v$ be dominant extremal weight vectors of $M$ and $N$,
respectively. 
Then there exists $a_{M,N}(z)\in\ko[[z]]^\times$
such that
$$
\Runiv_{M,N_z}\big( u \tens v_z\big)= a_{M,N}(z)\big( v_z\tens u \big).
$$
Thus we have  a unique $\ko(z)\tens\uqpg$-module isomorphism
\begin{align*}
 \Rnorm_{M,N_z}\seteq a_{M,N}(z)^{-1} &  \Runiv_{M,N_z}\big\vert_{\;\ko(z)\otimes_{\ko[z^{\pm1}]} ( M \otimes N_z) } 
\end{align*}
from $\ko(z)\otimes_{\ko[z^{\pm1}]} \big( M \otimes N_z\big)$ to $\ko(z)\otimes_{\ko[z^{\pm1}]}  \big( N_z \otimes M \big)$, which
satisfies
\begin{equation*}
\Rnorm_{M, N_z}\big( u  \otimes v_z\big) = v_z\otimes u .
\end{equation*}
We call $a_{M,N}(z)$ the {\it universal coefficient} of $M$ and $N$, and $\Rnorm_{M,N_z}$ the {\em normalized $R$-matrix}.

Let $d_{M,N}(z) \in \ko[z]$ be a monic polynomial of the smallest degree such that the image of $d_{M,N}(z)
\Rnorm_{M,N_z}(M\tens N_z)$ is contained in $N_z \otimes M$, which is called the {\em denominator of $\Rnorm_{M,N_z}$}. 
Then we have
\begin{equation*}
\Rren_{M,N_z}  =  d_{M,N}(z)\Rnorm_{M,N_z}
\col M \otimes N_z \To N_z \otimes M
\quad \text{up to a multiple of $\cz^\times$.}
\end{equation*}
Thus
\begin{align*}
 \Rren_{M,N_z} =a_{M,N}(z)^{-1}d_{M,N}(z)\Runiv_{M,N_z}
\quad \text{and} \quad  c_{M,N}(z)= \dfrac{d_{M,N}(z)}{a_{M,N}(z)}
\end{align*}
up to a multiple of $\ko[z^{\pm1}]^\times$.
In particular, $\Runiv_{M,N_z}$ is \rr whenever $M$ and $N$ are simple.

The following proposition was one of the main results of \cite{KKKO15}.
\Prop[{\cite[Theorem 3.12]{KKKO15}}] \label{Prop: r matrix hd soc}
Let $M$ and $N$ be simple modules, and assume that one of them is real.
Then $\Im(\rmat{M,N})$ is a simple module and it coincides with the head of $M\tens N$ and with the socle of $N\tens M$.
\enprop

Let $M$ and $N$ be simple modules in $\catC_\g$. Suppose that one of them is real. 
Thanks to Proposition \ref{Prop: r matrix hd soc}, the diagram
\begin{align} \label{Eq: cd for r}
	\xymatrix{
		M \otimes N \ar@{->>}[r]  \ar@/^{1.5pc}/[rrr]^{\rmat{M,N}}  & M \htens N\ar@{-}[r]^-\sim&
		N\stens M\akew\ar@{>->}[r] & N \otimes M
	}
\end{align}
commutes.  Here $\twoheadrightarrow$ denotes the natural projection and $\rightarrowtail$ denotes the embedding.

\begin{lemma} [{\cite[Corollary 3.13]{KKKO15}}] \label{Lem: MNDM}
	Let $L$ be a real simple module.
	Then for any simple module $X$, we have
	$$
	(L \htens X) \htens \dual L \simeq X, \qquad \dual^{-1} L \htens (X \htens L)  \simeq X,
	$$
	and 
	$$
	L \htens (X \htens \dual L) \simeq X, \qquad (\dual^{-1} L \htens X) \htens L  \simeq X.
	$$ 
\end{lemma}

\begin{lemma} [{\cite[Corollary 3.14]{KKKO15}}] \label{Lem: crystal for real}
	Let $X$, $Y$ and $L$ be simple modules in $\catC_\g$. Suppose that  $L$ is real. 
	\bnum
	\item $X \simeq L \htens Y$ if and only if $X \htens \dual L \simeq Y$,
	\item  $X \simeq Y \htens L$ if and only if $(\dual^{-1}L) \htens X \simeq Y$.
	\ee
\end{lemma}

In the following theorem, we refer \cite{Kas02} 
for the notion of good modules. We only note that the fundamental module
$V(\vpi_i)$ is a good module.

\begin{theorem}[{\cite{AK97,Chari02,Kas02,KKKO15}}]  \label{Thm: basic properties}
\hfill
\bnum
\item For good modules $M$ and $N$, the zeroes of $d_{M,N}(z)$ belong to
$\C[[q^{1/m}]]q^{1/m}$ for some $m\in\Z_{>0}$.
\item \label{it:comm} For simple modules $M$ and $N$ such that one of them is real, $M_x$ and $N_y$ strongly commute to each other if and only if $d_{M,N}(z)d_{N,M}(1/z)$ does not vanish at $z=y/x$.
\item  Let $M_k$ be a good module
with a dominant extremal vector $u_k$ of weight $\lambda_k$, and
$a_k\in\ko^\times$ for $k=1,\ldots, t$.
Assume that $a_j/a_i$ is not a zero of $d_{M_i, M_j}(z) $ for any
$1\le i<j\le t$. Then the following statements hold.
\bna
\item  $(M_1)_{a_1}\otimes\cdots\otimes (M_t)_{a_t}$ is generated by $u_1\otimes\cdots \otimes u_t$.
\item The head of $(M_1)_{a_1}\otimes\cdots\otimes (M_t)_{a_t}$ is simple.
\item Any non-zero submodule of $(M_t)_{a_t}\otimes\cdots\otimes (M_1)_{a_1}$ contains the vector $u_t\otimes\cdots\otimes u_1$.
\item The socle of $(M_t)_{a_t}\otimes\cdots\otimes (M_1)_{a_1}$ is simple.
\item  Let $\rmat{}\col (M_1)_{a_1}\otimes\cdots\otimes (M_t)_{a_t} \to (M_t)_{a_t}\otimes\cdots\otimes (M_1)_{a_1}$  be 
 $\rmat{ (M_1)_{a_1},\ldots, (M_t)_{a_t} }\seteq\prod\limits_{1\le j<k\le t}
\rmat{(M_j)_{a_j},\,(M_k)_{a_k}}$. 
Then the image of $\rmat{}$ is simple and it coincides with the head of $(M_1)_{a_1}\otimes\cdots\otimes (M_t)_{a_t}$
and also with the socle of $(M_t)_{a_t}\otimes\cdots\otimes (M_1)_{a_1}$.
\end{enumerate}
\item\label{Thm: bp5}
For any simple module $M\in\catC_\g$, there exists
a finite sequence $\st{(i_k,a_k)}_{1\le k\le t}$ 
in  $\sig$ \ro see $\eqref{eq:sig}$ below \rf\ 
such that
$M$ 
has $\sum_{k=1}^t \varpi_{i_k}$ as a dominant extremal weight
and it is isomorphic to a simple subquotient of
$V(\vp_{i_1})_{a_1}\tens\cdots V(\vp_{i_t})_{a_t}$.
Moreover, such a sequence $\st{(i_k,a_k)}_{1\le k\le t}$
is unique up to a permutation.

We call
$\sum_{k=1}^t(i_k,a_k)\in \awlP\seteq\Z^{\oplus  \sig}$, the {\em affine highest weight}
of $M$. 
\end{enumerate}
\end{theorem}

\subsection{Hernandez-Leclerc categories} \label{Sec: HL cat} \
For $i\in I_0$, let $m_i$ be a positive integer such that
$$
\weyl\pi_i\cap\bigl(\pi_i+\Z\delta\bigr)=\pi_i+\Z m_i\delta,
$$
where $\pi_i$ is an element of $\wl$ such that $\clp(\pi_i)=\varpi_i$.
Note that $m_i=(\al_i,\al_i)/2$ in the case when $\g$ is the dual of an untwisted affine algebra, and $m_i=1$ otherwise.
Then,
$V(\varpi_i)_x \simeq V(\varpi_i)_y  $ if and only if $x^{m_i}=y^{m_i}$ for $x,y\in \ko^\times$ (see \cite[Section 1.3]{AK97}).
We define
\eq
 \sig \seteq I_0 \times \cor^\times / \sim, 
\label{eq:sig}
\eneq
where the equivalence relation $\sim$ is given by 
\begin{align*} 
 (i,x) \sim (j,y) \Longleftrightarrow
V(\varpi_i)_x \simeq V(\varpi_j)_y
 \Longleftrightarrow\text{$i=j$ and $x^{m_i}=y^{m_j}$.}
\end{align*}
We denote by $[(i,a)]$ the equivalence class of $(i,a)$ in $ \sig$.  When no confusion arises, we simply write $(i,a)$ for the equivalence class $[(i,a)]$.
For $(i,x)$ and $ (j,y) \in \sig$, we put $d$ many arrows  from $(i,x)$ to $(j,y)$, where $d$ is the order of zeros of $d_{ V(\varpi_i), V(\varpi_j) } \allowbreak ( z_{V(\varpi_j)} / z_{V(\varpi_i)}  )$
at $  z_{V(\varpi_j)} / z_{V(\varpi_i)}  = y/x$.
Thus, $\sig$ has a quiver structure. 

We choose a connected component $\sigZ$ of $\sig$. Since 
a connected component of $\sig$ is unique up to a spectral parameter shift, 
$\sigZ$ is uniquely determined up to a quiver isomorphism.
We define $\catCO$ to be the smallest full subcategory of $\catC_\g$ such that 
\bna
\item $\catCO$ contains $V(\varpi_i)_x$ for all $(i,x) \in \sigZ$,
\item $\catCO$ is stable by taking subquotients, extensions and tensor products.
\ee
For symmetric affine types, this category was introduced in \cite{HL10}. 
Note that every simple module in $\catC_\g$ is isomorphic to a tensor product of certain spectral parameter shifts  
of some simple modules in $\catCO$ (\cite[Section 3.7]{HL10}).

\subsection{Invariants related to R-matrices}\

Let us recall the new invariants introduced in \cite{KKOP19C}.
We set
\begin{align*} 
 \varphi(z) \seteq \prod_{s=0}^\infty (1-\tp^{s}z)
 =\sum_{n=0}^\infty\hs{.3ex}\dfrac{(-1)^n\tp^{n(n-1)/2}}{\prod_{k=1}^n(1-\tp^k)}\;z^n
 \in\ko[[z]],
\end{align*}
where $p^* \seteq (-1)^{ \langle \rho^\vee, \delta \rangle} q^{\langle c, \rho \rangle}$  and  
$\tp \seteq (p^*)^2 = q^{2 \langle c, \rho \rangle}$.
 We consider the subgroup $\G$ of $\cor((z))^\times$ given by 
\begin{align*}
\G \seteq \left\{ cz^m \prod_{a \in \ko^\times} \varphi(az)^{\eta_a} \ \left|  \
\begin{matrix} \ c \in \ko^\times, \ m \in \Z , \\
\eta_a \in \Z \text{ vanishes except finitely many $a$'s. } \end{matrix} \right. \right\}.
\end{align*}

For a subset $S$ of $\Z$, let $\tp^{S} \seteq \{ \tp^k \ | \ k \in S\} $.
We define the group homomorphisms
\begin{align*}
\Deg :   \G \to  \Z \quad \text{ and } \quad \Di :   \G \to  \Z,
\end{align*} by
$$
\Deg(f(z)) = \sum_{a \in \tp^{\,\Z_{\le 0}} }\eta_a - 
\sum_{a \in \tp^{\,\Z_{> 0}} } \eta_a \quad \text{ and } \quad \Di(f(z)) = \sum_{a \in \tp^{\,\Z}} \eta_a
$$
for $f(z)=cz^m \prod_{  a\in\cor^\times } \varphi(az)^{\eta_a} \in \G$.

Note that
\eq
&&\Deg(f(z))=2\zero_{z=1}f(z)
\qt{for $f(z)\in\cor(z)^\times\subset\G$}
\eneq
(see \cite[Lemma 3.4]{KKOP19C}). 

\begin{definition} \label{def: Lams}
For non-zero $U_q'(\g)$-modules $M$ and $N$ such that $\Runiv_{M,N_z}$ is \rr, we define 
\begin{align*}
  \Lambda(M,N) & \seteq\Deg(c_{M,N}(z)),\\
\Lambda^\infty(M,N) & \seteq\Deg^\infty(c_{M,N}(z)), \\
\de(M,N) &\seteq  \frac{1}{2} (\La(M,N) + \La(N,M)).
\end{align*}
\end{definition}
Note that $\La(M,N)\equiv \Li(M,N)\mod 2$.

\begin{prop}  [\protect{\cite[Lemma 3.7, 3.8 and Corollary 3.23]{KKOP19C}}] \label{Prop: Li}
Let $M,N$ be simple modules in $\catC_\g$.
\bnum
\item  $\Li(M,N)   =-\Deg^\infty(a_{M,N}(z))$. 
\item  \label{Prop: Li: item2} $ \Li(M,N) = \Li(N,M)$.
\item $ \Li(M,N) = - \Li(\dual M, N) = - \Li(M,\dual N)$. 
\item In particular, $ \Li(M,N)= \Li( \dual M,\dual N)$.
\ee
\end{prop}

\begin{prop} [\protect{\cite[Lemma 3.7  and Proposition 3.18]{KKOP19C}}] \label{Prop: prop for La1}
Let $M,N$ be simple modules in $\catC_\g$.
\bnum
\item $ \La(M,N) =  \La(N, \dual M) =  \La( \dual^{-1}N, M) $. 
\item In particular, 
\begin{align*}
 \La(M,N) &= \La( \dual M , \dual N) = \La( \dual^{-1} M, \dual^{-1} N).
\end{align*}
\ee
\end{prop}

\begin{prop} [\protect{\cite[Proposition 3.11]{KKOP19C}}] \label{Prop: subquotients for Li}
Let $M$, $N$ and $L$ be non-zero modules in $\catC_\g$,
and let $S$ be a non-zero subquotient of $M\tens N$.
\bnum
\item
Assume that $\Runiv_{M,L_z}$ and $\Runiv_{N,L_z}$ are \rr.
Then $\Runiv_{S,L_z}$ is \rr and
$$ \Lambda(S,L)\le  \Lambda(M,L) + \Lambda(N,L)
\qtq
\Li(S,L)=\Li(M,L) + \Li(N,L).
$$
\item
Assume that $\Runiv_{L,M_z}$ and $\Runiv_{L,N_z}$ are \rr.
Then $\Runiv_{L,S_z}$ is \rr and
$$\Lambda(L,S)\le\Lambda(L,M) + \Lambda(L,N)
\qtq \Li(L,S)=\Li(L,M) + \Li(L,N).$$
\ee
\end{prop}

\begin{prop}  [\protect{\cite[Proposition 3.16]{KKOP19C}}] \label{Prop: d and d}
Let $M$ and $N$ be simple modules in $\catC_\g$. Then we have 
\bnum
\item  $ \de(M,N)=\zero_{z=1}\bl d_{M,N}(z)d_{N,M}(z^{-1})\br$,
\item  $\de(M,N) = \de(N,M)$.
\ee
In particular, we have $\de(M,N)\in\Z_{\ge0}$. 
\end{prop}

\begin{corollary}  [\protect{\cite[Corollary 3.17 and 3.20]{KKOP19C}}]  \label{Cor: comm}
Let $M$ and $N$ be simple modules in $\catC_\g$.
\bnum
\item
Suppose that one of $M$ and $N$ is real.
Then $M$ and $N$ strongly commute if and only if
$\de(M,N)=0$.
\item
In particular, if $M$ is real, then $ \La(M,M)=0.$
\ee

\end{corollary}

Proposition \ref{Prop: La and d} (i) and (ii) below were proved in \cite[Proposition 3.22]{KKOP19C} and Proposition \ref{Prop: La and d} (iii) is new.  
We add a whole proof of Proposition \ref{Prop: La and d} for the reader's convenience.

\begin{prop} \label{Prop: La and d}
For simple modules $M$ and $N$ in $\catC_\g$,  we have the followings:
\bnum
\item $\Lambda(M,N)=
 \sum_{k \in \Z} (-1)^{k+\delta(k<0)} \de(M,\D^{k}N)
=\sum_{k \in \Z} (-1)^{k+\delta(k>0)} \de(\dual^kM,N)$,
\item $\Lambda^\infty(M,N)= \sum_{k \in \Z} (-1)^{k} \de(M,\D^{k}N)$,
\item
$\zero_{z=1}c_{M,N}(z)=\displaystyle\sum_{k=0}^\infty(-1)^k\de(M,\D^kN)$. 
\ee
\end{prop}
\begin{proof}
We write $c_{M,N}(z) \equiv \prod_{a \in \cor^\times} \varphi(az)^{\eta_a} \mod \ \ko[z^{\pm1}]^\times$. Then we have
$$
 \dfrac{c_{M,N}(z)}{c_{M,N}(\tp z)} \equiv \prod_{a \in \cor^\times} (1-az)^{\eta_a},
$$
which yields that
\begin{align*}
\eta_{\tp^{-k}} & =\zero_{z= \tp^{k}}\left(\dfrac{c_{M,N}(z)}{c_{M,N}(\tp z)}\right)
=  \zero_{z= 1}\left(\dfrac{c_{M,N}(\tp^{k}z)}{c_{M,N}(\tp^{k+1} z)}\right) 
=  \zero_{z= 1}\left(\dfrac{c_{M,N_{ \tp^k }}(z)}{c_{M,N_{\tp^{k}}}(\tp z)}\right) \\
& \underset{(*)}{=}
 \zero_{z= 1} \left( \dfrac{d_{M,N_{\tp^k  }}(z)d_{N_{ \tp^k},M}(z^{-1})}
{d_{\dual^{-1}M ,N_{\tp^{k}}}(z)d_{N_{\tp^{k}}, \dual^{-1}M}(z^{-1})} \right)
 \underset{(**)}{=}
\de(M,N_{\tp^{k}}) - \de(\dual^{-1} M,N_{\tp^{k}}) \\
&  = \de(M,\D^{2k} N) - \de(M,\D^{2k+1}N),
\end{align*}
where $ (*)$ follows from \cite[Lemma 3.15]{KKOP19C} and 
$(**)$ from Proposition~\ref{Prop: d and d}.
Therefore, we have
\begin{align*}
\Lambda(M,N) & =\sum_{k \in \Z} (-1)^{\delta(k>0)} \eta_{\tp^k} 
 =\sum_{k \in \Z} (-1)^{\delta(k<0)} \eta_{\tp^{-k}} \\
&=\sum_{k \in \Z} (-1)^{\delta(k<0)} (\de(M,\D^{2k}N )- \de(M,\D^{2k+1}N)) \\
&= \sum_{k \in \Z} (-1)^{k+\delta(k<0)} \de(M,\D^{k}N ),
\end{align*}
which imply the first assertion (i).
Similarly, we obtain the second assertion (ii) as follows:
\begin{align*}
\Lambda^\infty(M,N) & = \sum_{k \in \Z} \eta_{\tp^k} = \sum_{k \in \Z}     (\de(M,\D^{-2k}N )-\de(M,\D^{-2k+1}N)) \\
& = \sum_{k \in \Z} (-1)^{k}  \de(M,\D^{k}N ). 
\end{align*}
 Finally we obtain
\begin{align*}
\zero_{z=1}c_{M,N}(z) &=\sum_{k=0}^\infty\eta_{\tp^{-k}}
=\sum_{k=0}^\infty\bl\de(M,\D^{2k} N) - \de(M,\D^{2k+1}N)\br \\
&=\sum _{k=0}^\infty(-1)^k\de(M,\D^{k} N),
\end{align*}
which gives the third assertion (iii).
\qedhere
\end{proof}

\Prop[{\cite[Corollary 4.12]{KKOP19C}}]\label{prop:tens<}
Let $L$ be a real simple module, and $M$ a simple module.
Assume that $\de(L,M)>0$. Then we have
$$\de(L,S)<\de(L,M)$$
for any simple subquotient $S$ of $L\tens M$
and also for any simple subquotient $S$ of $M\tens L$.
\enprop

The assumption in the following definition is slightly weaker than the one in \cite[Definition 4.14]{KKOP19C}.  
Under this weak assumption, the same statements as in \cite[Lemma 4.15 -- 4.18]{KKOP19C} can be proved in the same manner.

\begin{definition}
Let $L_1, L_2, \ldots, L_r$ be simple modules such that they are real except for at most one. 
The sequence $(L_1, \ldots, L_r)$ is called a \emph{normal sequence} if the composition of the R-matrices
\begin{align*}
\rmat{L_1, \ldots, L_r} \seteq&  \prod_{1 \le i < j \le r} \rmat{L_i, L_j} \\
 = & (\rmat{L_{r-1, L_r}}) \circ \cdots \circ ( \rmat{L_2, L_r}  \circ \cdots \circ \rmat{L_2, L_3} )  \circ ( \rmat{L_1, L_r}  \circ \cdots \circ \rmat{L_1, L_2} ) \\
& \cl L_1 \otimes L_2 \otimes  \cdots \otimes L_r \longrightarrow L_r \otimes \cdots \otimes L_2 \otimes L_1.
\end{align*}
does not vanishes.
\end{definition}

\Lemma[{\cite[Lemma 4.15]{KKOP19C}}] \label{Lem: normal head socle}
Let $(L_1, \ldots, L_r)$ be a normal sequence of simple modules
such that they are real except for at most one. 
Then $\Im (\rmat{L_1,\ldots, L_r})$ is simple
and it coincides with the head of $L_1\tens\cdots \tens L_r$ 
and also with the socle of $L_r\tens\cdots \tens L_1$.
\enlemma

\begin{lemma}  [\protect{\cite[Lemma 4.16]{KKOP19C}}] \label{Lem: normal}
Let $L_1, L_2, \ldots, L_r$ be simple modules such that they are real except for at most one. 
\bnum
\item If $(L_1, \ldots, L_r)$ is normal, then we have
\bna
\item
$(L_2, \ldots, L_r)$ and $(L_1, \ldots, L_{r-1})$ are normal.
\item
$
\La( L_1, \hd( L_2 \tens \cdots \tens L_r )) = \sum_{j=2}^r \La(L_1, L_j),
$ and \\
$
\La(  \hd( L_1 \tens \cdots \tens L_{r-1} ), L_r) = \sum_{j=1}^{r-1} \La(L_j, L_r).
$
\ee
\item Assume that $L_1$ is real, $(L_2, \ldots, L_r)$ is normal and
$$
\La( L_1, \hd( L_2 \tens \cdots \tens L_r )) = \sum_{j=2}^r \La(L_1, L_j),
$$
then $(L_1, \ldots, L_r)$ is normal.
\item Assume that $L_r$ is real,  $(L_1, \ldots, L_{r-1})$ is normal and 
$$
\La(  \hd( L_1 \tens \cdots \tens L_{r-1} ), L_r) = \sum_{j=1}^{r-1} \La(L_j, L_r),
$$
then  $(L_1, \ldots, L_r)$ is normal.
\ee
\end{lemma}

\begin{lemma}  [\protect{\cite[Lemma 4.3 and Lemma 4.17]{KKOP19C}}] 
\label{Lem: normal for 3}
Let $L,M,N$ be simple modules such that they are real except for at most one. 
If one of the following conditions
\bna
\item $\de(L,M)=0$ and $L$ is real, 
\item $\de(M,N)=0$ and $N$ is real,
\item $ \de(L, \dual^{-1}N) = \de(\dual L, N)=0$  and $L$ or $N$ is real 
\ee
holds,  then $(L,M,N)$ is a normal sequence, i.e., 
\begin{align*}
\La(L,M\htens N) = \La(L,M ) + \La(L,N), \quad \La(L\htens M, N) = \La(L,N ) + \La(M,N).
\end{align*}
\end{lemma}

\Lemma[{\cite[Corollary 4.18]{KKOP19C}}]\label{lem:LMN} Let $L, M,N$ be simple modules.
Assume that $L$ is real and one of $M$ and $N$ is real. 
Then $(L,M,N)$ is a normal sequence if and only if
$(M,N,\D L)$ is a normal sequence.
\enlemma
In \cite[Corollary 4.18]{KKOP19C}, we proved it in a stronger condition,
but the same proof still works without change.

\begin{lemma} \label{Lem: dual head}
	Let $L_1, L_2, \ldots, L_r$ be simple modules such that they are real except for at most one. 
	Suppose that the sequence $(L_1, \ldots, L_r)$ is normal.  For any $m\in \Z$, we have 
	$$
	\dual^{m}(  \hd ( L_1 \tens L_2 \tens  \cdots \tens L_r ) ) \simeq \hd ( \dual^{m} L_1 \tens \dual^{m} L_2 \tens  \cdots \tens \dual^{m} L_r ).
	$$ 
\end{lemma}
\begin{proof}
	It suffices to prove the case for $m= \pm1$. 
	We assume that $m = 1$. 
	Since the sequence $( \dual L_1, \ldots, \dual  L_r)$ is normal, 
	by Lemma \ref{Lem: normal head socle}, we have 
	\begin{align*}
		\dual (  \hd (  L_1 \tens \cdots \tens  L_r ) ) \simeq  \soc ( \dual L_r \tens \cdots \tens \dual L_1 ) \simeq 
		\hd ( \dual L_1 \tens \cdots \tens \dual  L_r ).
	\end{align*}
	
	The case for $m =-1$ can be proved in the same manner as above.	
\end{proof}

\begin{lemma} \label{Lem: normal for triple}
Let $L, M,N$ be simple modules.
Assume that $L$ is real and one of $M$ and $N$ is real. 
Then $\de(L, M\hconv N) = \de(L, M) + \de(L,N)$ if and only if
$(L,M,N)$ and $(M,N,L) $ are normal sequences.  
\end{lemma}
\begin{proof}
 By the assumption, we have 
\eqn
2\bl\de(L,M)+\de(L,N)-\de(L, M\hconv N)\br
&&=\bl(\La(L,M)+\La(L,N)-\La(L, M\hconv N)\br\\
&&\hs{5ex}+\bl\La(M,L)+\La(N,L)-\La(M\hconv N,L)\br.
\eneqn
Since $\La(L,M)+\La(L,N)-\La(L, M\hconv N)$ and 
$\La(M,L)+\La(N,L)-\La(M\hconv N,L)$ are non-negative by 
Proposition \ref{Prop: subquotients for Li}, we have
\eqn
&&\La(L, M\hconv N) = \La(L, M) + \La(L,N)\qtq
\La( M\hconv N, L) = \La( M, L) + \La(N, L),
\eneqn
if and only if $\de(L, M\hconv N) = \de(L, M) + \de(L,N)$.
Then the assertion follows from Lemma \ref{Lem: normal}. 
\end{proof}

\Cor\label{cor: normal for triple}
Let $L$ and $M$ be real simple modules 
and $X$ a simple module.
\bnum
\item
If $\de(L,M)=\de(\D^{-1} L,M)=0$,
then  we have
$\de(L, M\hconv X)=\de(L,X)$.
\item
If $\de(L,M)=\de(\D L,M)=0$,
then  we have
$\de(L, X\hconv M)=\de(L,X)$.
\ee
\encor
\Proof
(i)  The triples $(L,M,X)$ and $(M,X,L)$ are normal by 
Lemma~\ref{Lem: normal for 3}, and hence we have
$\de(L, M\hconv X)=\de(L,M)+\de(L,X)=\de(L,X)$ by the preceding lemma.

\snoi
(ii) can be proved similarly.
\QED

\smallskip

The following lemma can be proved similarly to
\cite[Proposition 3.2.17]{KKKO18}, and we do not repeat the proof here.
\Lemma\label{lem:MN1}
Let $M$ and $N$ be simple modules, and assume that one of them is real.
If $\de(M,N)=1$, then $M\tens N$ has length $2$ and we have an exact sequence
$$0\To N\hconv M\To M\tens N\To M\hconv N\To0.$$
\enlemma

The following lemma gives a criterion for a simple module to be real. 
\begin{lemma} \label{Lem: for real}
Let $X$ be a simple module such that $\de(X,X)=0$ and $X \otimes X$ has a simple head. Then $X$ is real.
\end{lemma}
\begin{proof}
Since $\de(X,X)=0$, we have $ \Rren_{X,X_z} \circ \Rren_{X_z,X} = f(z) \id $ for some $f(z) \in \cor(z)$ which is invertible at $z=1$.
Thus we have $\rmat{X,X}^2 \in \cor^\times \id $. By normalizing, 
we may assume that 
$ \rmat{X,X}^2 = \id$. Then we have
$$
X \otimes X = \Ker(\rmat{X,X}-\id) \oplus  \Ker( \rmat{X,X}+\id).
$$
Since  $X \otimes X$ has a simple head, we conclude that $ \rmat{X,X} $ should be $ \pm \id$, 
which implies the assertion by \cite[Corollary 3.3 and Theorem 3.12]{KKKO15}.
\end{proof}

\begin{lemma} \label{Lem: real for MN}
Let $M,N$ be real simple modules such that $ \de(M,N)=1 $. Then $M \htens N$ is real.
\end{lemma}
\begin{proof}
It follows from Proposition~\ref{prop:tens<} that
$$
\de(M, M\htens N) < \de(M,N) = 1, \quad \de(N, M\htens N) < \de(M,N) = 1,
$$
which implies that $ \de(M, M\htens N)=\de(N, M\htens N)=0$. We set $X \seteq  M\htens N$. 
Since $\de(M,X)=\de(N,X)=0$,
we have
$0\le\de(X,X)\le\de(M,X)+\de(N,X)=0$, i.e., 
$$
\de(X,X)=0.
$$
Since $N$ is real and $X\tens M$ is simple,
$(X \tens M)\tens N$ has a simple head. Thus the surjection 
$$
  (X \tens M)\tens N \twoheadrightarrow X \tens X
$$
tells us that $X \tens X$ has a simple head. Then the assertion follows from Lemma \ref{Lem: for real}.
\end{proof}

\Lemma\label{lem:simplylinked}
Let $M$ and $N$ be real simple modules such that $\de(M,N)=1$.
Then we have
\bnum
\item
$M\hconv N$ commutes with $M$ and $N$,
\item for any $m,n\in \Z_{\ge0}$, we have

$$M^{\tens m}\hconv N^{\tens n}
\simeq
\bc (M\hconv N)^{\tens m}\tens N^{\tens(n-m)}&\text{if $m\le n$,}\\
(M\hconv N)^{\tens n}\tens M^{\tens (m-n)}&\text{if $m\ge n$.}
\ec
$$
\ee
\enlemma
\Proof
(i) It follows from $\de(M, M\hconv N)\le\de(M,N)-1=0$ and
$\de(N, M\hconv N)\le\de(N,M)-1=0$.

\snoi
(ii) We shall prove only the first isomorphism.
We shall argue by induction on $m\le n$.
If $m=0$ it is obvious.
Assume that $m>0$.
Then we have
\eqn
M^{\tens m}\tens N^{\tens n}
&&\epito M^{\tens (m-1)}\tens (M\hconv N)\tens N^{\tens (n-1)}
\simeq (M\hconv N)\tens M^{\tens (m-1)}\tens N^{\tens (n-1)}\\
&&\epito (M\hconv N)\tens \bl (M\hconv N)^{\tens (m-1)}\tens N^{\tens (n-m)}\br
\simeq (M\hconv N)^{\tens m}\tens N^{\tens (n-m)}.
\eneqn
Then the assertion follows from the fact that
$(M\hconv N)^{\tens m}\tens N^{\tens (n-m)}$ is a simple quotient of
$M^{\tens m}\tens N^{\tens n}$ which has a simple head.
\QED

\Lemma
Let $M$ and $N$ be real simple modules such that $\de(M,N)=1$.
Then for any simple module $X$, we have
\bnum
\item The simple module $M\hconv (N\hconv  X)$ is isomorphic to either $(M\hconv N)\hconv  X$ or $(N\hconv M)\hconv  X$.
\item The simple module $(X\hconv M)\hconv N$ is isomorphic to either $X\hconv(M\hconv N)$ or $X\hconv (N\hconv M)$.
\ee
\enlemma

\Proof
Since the proof is similar, we prove only (i).
Let us consider a commutative diagram with an exact row:
$$\xymatrix{
0\ar[r]& (N\hconv M)\tens X\ar[r]\ar[dr]_f& M\tens N\tens X\ar[r]\ar@{->>}[d]&(M\hconv N)\tens X\ar[r]&0\\
&&M\htens (N\htens X)}
$$
The exactness follows from Lemma~\ref{lem:MN1}.
By Lemma~\ref{Lem: real for MN}, $M\hconv N$ and $N\hconv M$ are real simple modules.
If $f$ does not vanish, then we have
$(N\hconv M)\hconv X\simeq M\htens (N\htens X)$.

\snoi
If $f$ vanishes then there exists an epimorphism
$(M\hconv N)\tens X\epito M\htens (N\htens X)$ and hence
$(M\hconv N)\hconv X\simeq M\htens (N\htens X)$.
\QED

\vskip 2em 

\section{Root modules} \label{Sec: root modules}

In this section, we investigate properties of \emph{root modules}. 

\begin{definition}
A module $L\in \catC_\g$ is called a \emph{root module} if $L$ is a real simple module such that
\eq\label{eq:root}
&\phantom{aaaa}&\de\bl L, \dual^kL\br\ms{3mu} = \ms{4mu}
\delta(k=\pm 1)
\qtext{for any $k\in\Z$.}
\eneq
\end{definition}
Note that for a root module $L$, we have
$$
\Li(L,L)=-2
$$
by Proposition \ref{Prop: La and d}.

The name ``root module'' comes from Lemma~\ref{lem:root} below.

\begin{example} 
Using the denominators for fundamental modules (see \cite[Appendix A]{KKOP20} for example), 
one can easily prove that any fundamental module $V(\varpi_i)_a$ ($i\in I_0 $, $a \in \cor^\times$) is a root module.
\end{example}

\subsection{Properties of  root modules}

\begin{lemma} \label{Lem:Dk L X}
Let $L$ be a root module and let $X$ be a simple module. 

\bnum
\item For $k\in \Z$, we have 
\eqn
&&\de( \dual^k L, X  )- \delta(k=0,2)\le
\de( \dual^k L, L \htens X  ) \le \de( \dual^k L, X ) + \delta(k=\pm1).
\eneqn
In particular, 
$ \de( \dual^k L, L\htens X ) = \de(\dual^k L, X) $ for $k \ne -1,0,1,2$, and
\eqn
\de( \dual^k L, L \htens X  )-\de( \dual^k L, X )
\in\bc\st{0,1}&\text{for $k=\pm1$,}\\
\st{0,-1}&\text{for $k=0,2$.}
\ec
\eneqn
Moreover, we have
\eqn
&& \bl \de( \dual^{-1} L, L \htens X) - \de(\dual^{-1}  L,X)  \br
+ \bl \de(   L,X) -  \de(   L, L \htens X)\br \\*
&&\hs{5ex}+ \bl \de( \dual  L, L \htens X) - \de(\dual   L,X) \br 
+ \bl   \de(\dual^{2}  L,X) - \de( \dual^{2} L, L \htens X) \br =2.
\eneqn

\item For $ k\in \Z$, we have 
\eqn
\de( \dual^k L,  X  )-\delta(k=-2,0)\le
\de( \dual^k L, X \htens L  ) \le \de( \dual^k L, X ) + \delta(k=\pm1).
\eneqn
In particular, $ \de( \dual^k L, X\htens L ) = \de(\dual^k L, X) $
for $k \ne -2,-1,0,1$, and 
\eqn
\de( \dual^k L, X\htens L  )-\de( \dual^k L, X )
\in\bc\st{0,1}&\text{for $k=\pm1$,}\\
\st{0,-1}&\text{for $k=0,-2$.}
\ec
\eneqn
\ee
Moreover, we have
\eqn
&& \bl \de( \dual^{-2} L, X) - \de(\dual^{-2}  L,X\htens L)  \br
+\bl \de( \dual^{-1} L, X \htens L) - \de(\dual^{-1}   L,X) \br \\*
&&\hs{10ex}
+ \bl \de(L,X) -  \de(   L, X \htens L)\br + \bl \de( \dual L, X \htens L) - \de(\dual   L,X) \br 
=2.
\eneqn

\end{lemma}
\begin{proof}

(i) 
 By \cite[Proposition 4.2]{KKOP19C}, we have 
$$
\de( \dual^k L, L \htens X  ) \le \de( \dual^k L, X ) + \de( \dual^k L, L ) = \de( \dual^k L, X ) + \delta(k=\pm1).
$$
By the same reason, it follows from $X \simeq ( L \htens X   ) \htens \dual L $ (see Lemma \ref{Lem: MNDM}) that 
\begin{align*}
\de( \dual^k L,  X  ) &= \de( \dual^k L,  ( L \htens X   ) \htens \dual L  )  \\
&\le \de( \dual^k L, L \htens  X  ) + \de( \dual^k L, \dual  L ) \\
&= \de( \dual^k L, L \htens X) + \delta(k=0,2).
\end{align*}
Hence we obtain the first assertion.
Since 
\begin{align} \label{Eq: d2LLX2}
\Li( L, L \htens X) = \Li(L,L) + \Li(L,X) = -2 + \Li(L,X),
\end{align}
it follows from $\eqref{Eq: d2LLX2}$ and  Proposition \ref{Prop: La and d} that  
\begin{align*}
2 &=  \Li(L,X) - \Li( L, L \htens X)  \\
&=\sum_{k\in\Z}(-1)^k\bl\de(\D ^k L,X)-\de(\D^kL, L\hconv X)\br\\
&=  \bl \de( \dual^{-1} L, L \htens X) - \de(\dual^{-1}  L,X)  \br
+ \bl \de(   L,X) -  \de(   L, L \htens X)\br \\
& \ \ \ \  + \bl \de( \dual  L, L \htens X) - \de(\dual   L,X) \br 
+ \bl   \de(\dual^{2}  L,X) - \de( \dual^{2} L, L \htens X) \br,
\end{align*}
which yields the last assertion.

\snoi
(ii) Using the fact that $X \simeq  (\dual^{-1} L) \htens (X \htens L )  $, it can be proved in the same manner as above.
\end{proof}

\begin{lemma} \label{Lem: dLX>0} 
Let $L$ be a root module and let $X$ be a simple module. Suppose that  
$$
\de(L,X) > 0.
$$
Then we have the following: 
\bnum
\item $\de(L, L \htens X) = \de(L,X) -1$  and $\de(\dual^{-1} L, L \htens X) = \de(\dual^{-1}L, X)$,

\item 
$\de(L, X \htens L) = \de(L,X) -1$  and 
$\de(\dual  L, X \htens L) = \de(\dual L, X)$.
\ee 
\end{lemma}
\begin{proof}
We shall prove only (i) since the proof of (ii) is similar.

Since $\de(L,X) > 0$, $L$ does not commute with $X$. By 
Proposition~\ref{prop:tens<}, we have 
$$
\de( L, L \htens X) < \de(L, X).
$$
On the other hand, Lemma \ref{Lem:Dk L X} implies
$\de(L,X)\le\de(L,L\hconv X)+1$
which implies the first assertion.

\snoi
Let us show the second equation in (i).
By Lemma \ref{Lem:Dk L X}, we have $\de(\dual^{-1} L, L \htens X) = \de(\dual^{-1} L, X)$ or $\de( \dual^{-1} L,X )+1$.
If
\begin{align*} 
\de(\dual^{-1} L, L \htens X) = \de(\dual^{-1} L, X)+1,
\end{align*}
then Lemma \ref{Lem: normal for triple} says that $ (\dual^{-1} L,L,X)$ is a normal sequence,  and hence $(L,X,L)$ is also a normal sequence by 
Lemma~\ref{lem:LMN},
which implies that  
$$
L \htens X \simeq X \htens L.
$$
But it contradicts $\de(L,X) > 0$. 
Therefore we conclude that $\de(\dual^{-1} L, L \htens X) = \de(\dual^{-1} L, X)$. 
\end{proof}

\Lemma
Let $L$ be a root module and $X$ a simple module.
\bnum
\item
Assume one of the following conditions:
\bna
\item $\de(\D L, L \htens X)>0$,
\item $\de(\D L, X)>0$, 
\item $\de(\D^2L,X)=0$.
\ee
Then we have
$$\de(\D L, L\hconv X)=\de(\D L,X)+1.$$
\item
Assume one of the following conditions:
\bna
\item $\de(\D^{-1}L, X \htens L)>0$,
\item $\de(\D^{-1}L, X )>0$, 
\item $\de(\D^{-2}L,X)=0$.
\ee
Then we have
$$\de(\D^{-1}L, X\hconv L)=\de(\D^{-1}L,X)+1.$$
\ee
\enlemma
\Proof
We shall only prove (i) since the proof of (ii) is similar.

\snoi
(a)\ Assume first that $\de(\D L, L\htens X)>0$.
Setting $Y=L\hconv X$, we have $X\simeq Y\hconv \D L$.
Hence Lemma~\ref{Lem: dLX>0} (ii)  implies 
$\de(\D L, Y\hconv \D L)=\de(\D L, Y)-1$.

\snoi
(b) If $\de(\D L,X)>0$, then we have
$\de(\D L,L\htens X)\ge \de(\D L,X)>0$ by Lemma~\ref{Lem:Dk L X}.

\snoi
(c) Finally, assume that $\de(\D^2 L,X)=0$.
If $ \de( L, X)=0$, then  we have 
\begin{align*}
\de(\dual L, L \htens X) &=  \de ( \dual L, L) + \de(\dual L, X) \\
&= \de(\dual L, X)+1.
\end{align*}

Suppose that $ \de( L, X)>0$. 
Since $\de( \dual^2 L, X ) = 0$ and $\de( \dual^2 L, L )=0$, we have 
\begin{align} \label{Eq: d2LLX}
\de(\dual^2 L, L \htens X )  =   0.
\end{align}
Moreover, Lemma \ref{Lem: dLX>0} tells us that 
\begin{align}
\de(L, L\htens X) =\de(L,X)-1\qtq \de(\dual^{-1} L, L \htens X) = \de(\dual^{-1}L, X).\label{Eq: dLLX}
\end{align}

By Lemma~\ref{Lem:Dk L X}, \eqref{Eq: d2LLX} and \eqref{Eq: dLLX}, we have
\eqn
&& 2=\bl \de( \dual^{-1} L, L \htens X) - \de(\dual^{-1}  L,X)  \br
+ \bl \de(   L,X) -  \de(   L, L \htens X)\br \\
&&\hs{5ex}+ \bl \de( \dual  L, L \htens X) - \de(\dual   L,X) \br 
+ \bl   \de(\dual^{2}  L,X) - \de( \dual^{2} L, L \htens X) \br\\
&&\hs{3ex}=1+ \bl \de( \dual  L, L \htens X) - \de(\dual   L,X) \br,
\eneqn
which yields the desired result.
\QED

\begin{lemma} \label{Lem: DLX = LX+1} 
Let $L$ be a root module, $X$ a simple module, and $k\in \Z_{ \ge  0}$.
\bnum
\item Suppose that one  of the following conditions
\bna
\item $\de( \dual L,\, L^{\tens k} \htens X) \ge k$,
\item  $\de( \dual^2 L, X ) = 0$  
\ee
is true. Then  we have 
$$
\de(\dual L,\, L^{\tens k} \htens X) = \de( \dual L, X ) + k.
$$
\item Suppose that one of the following
\bna
\item  $\de( \dual^{-1} L, X \htens L^{\tens k}) \ge k$,
\item  $\de( \dual^{-2} L, X ) = 0$  
\ee
is true. Then we have 
$$
\de\bl\dual^{-1} L,\ms{2mu} X \htens L^{\tens k}\br = \de( \dual^{-1} L, X ) + k.
$$
\ee
\end{lemma}
\begin{proof}
This follows from the preceding lemma by induction on $k$.
Note that 
$\de(\D^2L, \allowbreak  L^{\tens k}\hconv X)=0$ as soon as $\de(\D^2L,X)=0$.
\end{proof}

\begin{prop}
Let $L$ be a root module and $X$ a simple module.
If  $ \de( L,  X)\not=1 $,  then $(L \hconv X)\hconv L$ is isomorphic to $L \hconv (X \hconv L)$.
\end{prop}
\begin{proof}
 If $\de(L,X)=0$, then it is obvious.
Hence we may assume that $\de(L,X)\ge2$.

Set $Y\seteq (L\hconv X)\hconv L$.
Then, we have
$$L\hconv X\simeq \dual^{-1}L\hconv Y\qtq X\simeq
(\dual^{-1}L\hconv Y)\hconv \D L.$$
Lemma \ref{Lem: dLX>0} says that 
$\de(L,L\hconv X)=\de(L,X)-1>0$ and
$\de(L,Y)=\de(L,L\hconv X)-1$.
Hence we obtain
$$\de(L,\D^{-1}L\hconv Y)=\de(L,Y)+1=\de(L,\D^{-1}L)+\de(L,Y).$$
Then Lemma~\ref{Lem: normal for triple} implies that
$(L,\D^{-1}L,Y)$ is a normal sequence,
and 
$(\D^{-1}L, Y, \allowbreak  \D L)$ is also a normal sequence
by Lemma~\ref{lem:LMN}.
Hence we have
$$(\D^{-1}L\hconv Y)\hconv\D L\simeq \D^{-1}L\hconv(Y\hconv\D L).$$
Since
$X\simeq \D^{-1}L\hconv(Y\hconv\D L)$, we obtain
$Y\simeq L\hconv (X\hconv L)$.\qedhere
\end{proof}

\subsection{Properties of pairs of  root modules} \

Let $L$ and $L'$ be root modules. 
Throughout this subsection, we assume that  
\begin{align} \label{Eq: assumption}
 \de( \dual^k L, L')\ms{5mu} = \ms{7mu}\delta\bl k=0\br \quad \text{ for } k\in \Z.
\end{align}
Note that, by Proposition~\ref{Prop: La and d}, 
\begin{itemize}
\item $\La(L,L')=\Li(L,L')=1$.
\end{itemize}

\Lemma\label{lem:real}
The simple module $L\hconv L'$ is a root module.
\enlemma
\Proof
Set $L''\seteq L\hconv L'$.
By Lemma~\ref{Lem: real for MN}, $L''$ is real.

It is obvious that
$\de(\dual^kL,L'')=\de(\dual^kL',L'')=0$ for $k\not=0,\pm1$.
On the other hand, Lemma~\ref{Lem: dLX>0} implies that
$\de(L,L'')=\de(L',L'')=0$.
Hence, we have
$\de(\dual^kL'',L'')=0$ for $k\not=\pm1$.
Now, we have
$$\Li(L'',L'')=\Li(L,L)+2\Li(L, L')+\Li(L',L')=(-2)+2+(-2)=-2.$$
Then Proposition~\ref{Prop: La and d} implies that
$$-2=\sum_{k\in\Z}(-1)^k\de(\dual^kL'',L'')
=-\de(\dual L'',L'')-\de(\dual^{-1}L'' , L'').$$
Since $\de(\dual L'',L'')=\de(\dual^{-1}L'' , L'')$,
we obtain
$\de(\dual^{\pm1} L'',L'')=1$. 
\QED

\begin{lemma} \label{Lem: two root modules}
We have
$$\de(\dual^k L,L\hconv L')\ms{3mu}=\ms{4mu}\delta(k=1)\qtq
\de(\dual^k L,L'\hconv L)\ms{3mu}=\ms{4mu}\delta(k=-1).
$$
\end{lemma}
\begin{proof}
Since
$\de(\dual^k L,L\hconv L')\le \de(\dual^k L,L)+\de(\dual^k L,L')$, we have
$\de(\dual^k L,L\hconv L')=0$ for $k\not=-1,0,1$.
It follows from Lemma \ref{Lem: dLX>0} that 
\begin{align*}
\de(L, L \hconv L') &= \de(L,L')-1 = 0, \\
 \de(\dual^{-1} L, L \hconv L') &= \de(\dual^{-1} L,L') =0.
\end{align*}

On the other hand, since 
\begin{align*}
\de(\dual L, L \hconv L') & \le \de(\dual L, L ) + \de(\dual L,  L')=1, 
\end{align*}
we have $\de(\dual L, L \hconv L')\in \{0,1\}$.

If $ \de(\dual L, L \hconv L')=0$, then we have
$$
L' \simeq (L \hconv L')\hconv \dual L \simeq (L \hconv L')\otimes \dual L,
$$
which implies that 
$$
 0= \de ( \dual^2 L, L' ) = \de ( \dual^2 L, (L \htens L')\tens \dual L ) = \de(\dual^2 L, L \hconv L') + \de (\dual^2 L, \dual L) \ge 1. 
$$ 
This is a contradiction. Hence, we conclude that $ \de(\dual L, L \hconv L')=1$.
Thus we obtained the first equality.

The second equality can be proved similarly.\qedhere
\end{proof}

\begin{lemma} \label{Lem:LL'X} 
Let $X$ be a simple module.
\bnum
\item If $k\ne 0,1$, then 
$$
\de(  \dual^k L, L' \hconv X) = \de(\dual^k L, X).
$$
 As for $k=0$ and $1$,  one and only one of the following  two statements is true.
\bna
\item $ \de(L, L' \hconv X) = \de(L,X) $ and $\de( \dual L, L' \hconv X) = \de (\dual L, X)-1$,
\item $ \de(L, L' \hconv X) = \de(L,X) +1 $ and $\de( \dual L, L' \hconv X) = \de (\dual L, X)$.
\ee
\item If $k\ne -1,0$, then 
$$
\de(  \dual^k L, X \hconv L') = \de(\dual^k L, X).
$$ 
 As for $k=-1$ and $0$, one and only one of the following two
statements is true.
\bna
\item $ \de(L, X \hconv L') = \de(L,X) $ and $\de( \dual^{-1} L, X \hconv L' ) = \de (\dual^{-1} L, X)-1$,
\item $ \de(L,  X\hconv L') = \de(L,X) +1 $ and $\de( \dual^{-1} L, X \hconv L') = \de (\dual^{-1} L, X)$.
\ee
\ee

\end{lemma}
\begin{proof}
(i)
By \cite[Proposition 4.2]{KKOP19C}, we have 
\begin{align*}
\de(\dual^k L, L' \hconv X) &\le \de(\dual^k L, X) + \de(k=0), \\
 \de(\dual^k L, X) &\le \de(\dual^k L, L' \hconv X) + \de(k=1),
\end{align*}
where the second inequality follows from $X  \simeq (L'\hconv X)\hconv \dual L'$.
The above inequalities give the first assertion and 
\begin{equation} \label{Eq: de LL'}
\begin{aligned} 
& \text{ $\de( L, L' \hconv X) = \de(   L, X) $ or $\de(   L, X)+1$,}\\
& 
\text{ $\de(\dual L, L' \hconv X) = \de( \dual L, X) $ or $\de(\dual L, X)-1$.} 
\end{aligned}
\end{equation}
By the assumption $\eqref{Eq: assumption}$, we have $\Li(L,L')=1$, which implies 
\begin{align*}
1  &= \Li(L, L' \hconv X) - \Li(L,X)=\sum_{k\in\Z}(-1)^k\bl
\de(\dual^k L, L' \hconv X) -\de(\dual^k L,X)\br \\
&=\bl \de(L, L' \hconv X) - \de(L, X)  \br +  \bl \de(\dual L, X) - \de( \dual L, L' \hconv X   ) \br.
\end{align*}
 Then \eqref{Eq: de LL'} implies
the second assertion. 

\snoi
(ii) can be similarly proved 
by using $ X \simeq \dual^{-1} L' \htens (X \htens L')$. 
\end{proof}

\begin{prop} \label{Prop: dLL'X}
Let $X$ be a simple module. 
\bnum
\item If $\de( \dual L, X )=0$, then we have
$$
\de(L, L' \hconv X) = \de(L,X) + 1, \qquad \de( \dual L, L' \hconv X )=0.
$$
\item If $\de( \dual^{-1} L, X )=0$, then we have
$$
\de(L, X \hconv L') = \de(L,X) + 1, \qquad \de( \dual^{-1} L, X \hconv L' )=0.
$$
\ee
\end{prop}
\begin{proof}
(i)
Since $\de( \dual L, L' \hconv X ) \le\de( \dual L, L'  ) + \de( \dual L,   X )=0  $,  Lemma \ref{Lem:LL'X} tells us that $\de(L, L' \hconv X) = \de(L,X) + 1$.

\snoi
(ii) can be proved in the same manner as above. 
\end{proof}

\begin{corollary} \label{Cor: de for LL' left}
Let $n \in \Z_{\ge0}$ and let $X$ be a simple module.
\bnum
\item If $\de(\dual L, X)=0$, then we have 
$$
\de(L, L'^{\ms{3mu} \otimes n} \hconv X ) = \de(L, X) + n \quad  \text{ and } \quad \de( \dual L, L'^{\ms{3mu}  \otimes n} \hconv X)=0.
$$ 

\item If $\de(\dual^{-1} L, X)=0$, then we have 
$$
\de(L,  X \htens L'^{\ms{3mu}  \otimes n} ) = \de(L, X) + n \quad  \text{ and } \quad \de( \dual^{-1} L,  X \htens L'^{ \ms{3mu} \otimes n} )=0.
$$ 

\ee
\end{corollary}
\begin{proof}
They easily follow from Proposition \ref{Prop: dLL'X} by
induction on $n$.
\end{proof}

\begin{prop}\label{prop:XY}
Let $m \in \Z_{\ge0}$ and let $Y$ be a simple module. Set 
$$
X \seteq L^{\tens m} \htens Y.
$$ 
Suppose that 
$$
\de(L,X)=0, \quad \de(\dual L', Y )=0, \quad \de( \dual^{t} L, Y)=0  \quad \text{ for $t=1,2$}.
$$
Then  
\bni
\item $\de(\dual L, X)=m$, 
\item for any integer $k$ such that $0\le k\le m$, we have 
$$
\de(L, L'^{\tens k} \htens X)=0, \qquad \de( \dual L, L'^{\tens k} \htens X)= m - k,
$$
\item for any integer $k \ge m$, we have 
$$
\de(L, L'^{\tens k} \htens X)=k-m, \qquad \de( \dual L, L'^{\tens k} \htens X)= 0.
$$
\ee
\end{prop}
\begin{proof}

(i) As $\de( \dual^{t} L, Y)=0$ for $t=1,2$, we have 
$$
\de(\dual L, X) = \de(\dual L, L^{\tens m} \htens Y) = \de(\dual L, Y) + m =m
$$
by Lemma \ref{Lem: DLX = LX+1} (i).

\mnoi
(ii) 
We have 
\begin{itemize}
\item  $L' \htens L$ commutes with
$L$, $L'$ and $\D L$ by Lemma \ref{Lem: two root modules},
\item $\de(\dual L', X)=0$ because $\de(\dual L', Y )=0$ and $\de(\dual L', L )=0$,
\item the triple $(L'^{\tens a}, A, X)$ is normal because $\de(\dual L', X)=0 $,
\item the triple $(L'^{\tens a}, B, L )$ is normal 
because $\de(\dual L', L )=0$,
\item the triples $(L'^{\tens a}, A, Y )$, $(L^{\tens a}, A, Y )$ and 
$\bl(L' \htens L)^{\tens a} , A, Y \br$ are normal
because $\de(\dual L', Y )=0$ and $\de(\dual L , Y )=0$,
\end{itemize}
Here $A$ is a real simple module, $B$ is a simple module, and $a$ is an arbitrary non-negative integer.
We will use these fact freely in the subsequent arguments.

For any $k$ such that $0\le k\le m$, we have 
\begin{align*}
L'^{\tens k} \htens X & \simeq L'^{\tens k} \htens ( L ^{\tens m } \htens Y) \\ 
&\simeq ( L'^{\tens k} \htens  L ^{\tens m } ) \htens Y \\
& \simeq \bl  (L' \htens  L)^{\tens k} \tens L^{\tens m-k} \br \htens Y \\
& \simeq \bl   L^{\tens m-k} \tens (L' \htens  L)^{\tens k}   \br \htens Y \\
& \simeq    L^{\tens m-k} \htens \bl (L' \htens  L)^{\tens k}    \htens Y\br,
\end{align*}
where the third isomorphism follows from 
Lemma~\ref{lem:simplylinked}.
Thus, we have, for $1\le k\le m$,
\begin{align*}
L\htens (L'^{\tens k} \htens X ) & \simeq L\htens \bl (  L^{\tens m-k} \tens (L' \htens  L)^{\tens k}   ) \htens Y \br \\
& \simeq \bl L^{\tens m-k+1}  \tens   (L' \htens  L)^{\tens k} \br    \htens Y \\ 
& \simeq  \bl (L' \htens  L) \tens  L^{\tens m-k+1}  \tens   (L' \htens  L)^{\tens k-1} \br    \htens Y \\
& \simeq   (L' \htens  L) \htens \bl( L^{\tens m-(k-1)}  \tens   (L' \htens  L)^{\tens k-1} )    \htens Y\br \\
& \simeq   (L' \htens  L) \htens  (L'^{\tens k-1} \htens X ),
\end{align*}
and 
\begin{align*}
(L'^{\tens k} \htens X ) \htens L  & \simeq L'^{\tens k} \htens ( X  \tens L )
\simeq L'^{\tens k} \htens ( L  \tens X ) \\
& \simeq (L'^{\tens k} \htens  L)  \htens X \\
&\simeq  \bl (L' \htens  L) \tens  L'^{\tens k-1} \br \htens X \\
&\simeq  (L' \htens  L) \htens  (L'^{\tens k-1} \htens X ).
\end{align*}
This tells us that
$$
L\htens (L'^{\tens k} \htens X ) \simeq (L'^{\tens k} \htens X ) \htens L,
$$
which implies that $\de( L,  L'^{\tens k} \htens X) = 0$.

\smallskip 

On the other hand, since $\de(\dual^2 L, Y)=0$ and $\de(\dual^2 L, L'\hconv L )=0$, we have 
$$
\de\bl\dual^2 L,  (L' \htens  L)^{\tens k}    \htens Y\br=0.
$$
Then, Lemma \ref{Lem: DLX = LX+1} implies that 
\begin{align*}
\de( \dual L, L'^{\tens k} \htens X ) &= \de\Bigl(\dual L,\ms{5mu} L^{\tens m-k} \htens 
\bl (L' \htens  L)^{\tens k}    \htens Y\br\Bigr) \\
&= \de\bl \dual L, (L' \htens  L)^{\tens k}    \htens Y \br +    m-k \\
&  = m-k,
\end{align*}
where the last equality follows from $\de(\dual L, L' \htens L)=0$ and $ \de( \dual L, Y)=0 $.

\mnoi
(iii) By (ii), we have 
$$
\de(L, L'^{\tens m} \htens X)=0, \qquad \de( \dual L, L'^{\tens m} \htens X)= 0.
$$
Since 
$$
L'^{\tens k} \htens X \simeq L'^{\tens k-m} \htens( L'^{\tens m} \htens X),
$$
we have the assertion by Corollary \ref{Cor: de for LL' left} (i).
\end{proof}

\vskip 2em

\section{\CWS} \label{Sec: SW}

Let $\ddD \seteq \{ \Rt_i \}_{i\in J} \subset \catC_\g$
be a family of simple modules of $\catC_\g$.
The family $\ddD $ is called a \emph{duality datum} associated with a generalized Cartan matrix $\cmC = (c_{i,j})_{i,j\in J}$ of symmetric type if it satisfies the following:
\bna
\item for each $i\in J$, $\Rt_i$ is a real simple module,
\item for any $i,j\in J$ such that $i\ne j$, $\de(\Rt_i, \Rt_j) = -c_{i,j}$.
\ee
Then one can construct a monoidal functor 
$$
\F_\ddD \col R_\cmC\gmod \longrightarrow \catC_\g
$$
 using the duality datum $\ddD$ (see \cite{KKK18A, KP18}).

The functor $\F_\ddD$ is called a \emph{\WS functor} or shortly a \emph{duality functor}.

In \S\;\ref{subsec:WS} below, we slightly modify the definition of
\WS functor in order that it commutes with the affinization.

\subsection{Affinizations} \

\subsubsection{Pro-objects}
Let $\cor$ be a base field and let $\shc$ be an essentially small $\cor$-abelian category.
Let $\Pro(\shc)$ be the category of pro-objects of $\shc$ (see \cite{KS} for details).  One can show that 
$$
\Pro(\shc)\simeq\bst{\text{left exact $\cor$-linear   functors from $\shc$
to $\cor\text{-Mod}$}}^\opp
$$
by the functor
\eqn
&&\proolim[i] M_i\To\Bigr(\shc\ni X\mapsto\indlim\Hom_{\ms{2mu}\shc}(M_i,X)\Bigr).
\eneqn
Here, $\cor\text{-Mod}$ is the category of
vector spaces over $\cor$, and $\proolim$ denotes the \emph{pro-lim} (see \cite[Section 2.6 and Proposition 6.1.7]{KS} for notations and details).
Then, $\Pro(\shc)$ is a $\cor$-abelian category which admits small 
projective limits. If no confusion arises, 
we regard $\shc$ as a full subcategory of $\Pro(\shc)$, which is stable by extensions and subquotients.
Any functor $F\cl \shc\to \shc'$ extends to
$\P F\cl \Pro(\shc)\to\Pro(\shc')$
which commutes with small filtrant projective limits:
$$
\P F(\proolim[i] M_i)\simeq \proolim[i] F(M_i).
$$

\smallskip

\subsubsection{Affinization in  quiver Hecke algebra case}
Let $R$ be a symmetric quiver Hecke algebra.
Note that
$$
R(\beta)\Mod\to \Pro\bl R(\beta) \gmod\br.
$$
Recall that $R(\beta)\Mod$ is the category of graded $R(\beta)$-modules. 
Let
$$
\Pro(R)\seteq\soplus_{\beta \in \rootl^+}\Pro(R(\beta)\gmod ),
$$
which is a monoidal category.  Let $z$ be an indeterminate of homogeneous degree $2$, and we set 
$$
R(\beta)^\aff \seteq \cor[z]\tens_\cor \one_zR(\beta),
$$
which has the graded $R(\beta)$-bimodule structure.
Here $ \one_zR(\beta)$ is a free right $R(\beta)$-module of rank one and
the left module structure is given by
$$
e(\nu)\one_z = \one_z e(\nu), \quad x_k\one_z=\one_zx_k+z\one_z\qtq\tau_k\one_z=\one_z\tau_k.
$$
Hence we have
\eq
\one_z x_k=(x_k-z)\one_z.
\label{eq:right}
\eneq
For $X\in R(\beta)\gmod$, the affinization $X^\aff$ of $X$ is isomorphic to $R(\beta)^\aff\tens_{R(\beta)}X.$
Since $X^\aff$ is not in $R(\beta) \gmod$, we set
$$
X^\Aff\seteq\proolim[m]X^\aff/z^m X^\aff\in\Pro(R(\beta) \gmod).
$$
Note that 
$$
X^\Aff\simeq\cor[[z]]\tens_\cor X
$$
as an object of $\Pro(\cor\smod)$ forgetting the action of $R(\beta)$.
Here we regard $\cor[[z]]$ as the object of $\Pro(\cor\smod)$:
$$\proolim[m]\cor[z]\;/\;\cor[z]z^m.$$
Similarly we set
\eqn
&&R(\beta)^\Aff \seteq \proolim[m] R(\beta)^\aff
/\R(\beta)^\aff(z,x_1,\ldots, x_{\height{\beta}})^m, 
\eneqn
which is an object of $\Pro(R(\beta)\gmod )$ with a right $R(\beta)$-action.
Here, $ ( z,x_1,\ldots, x_{ m})$ is the ideal of $\cor[z, x_1, \ldots, x_m]$ generated by $z, x_1, \ldots, x_m$. 
 Then we have
$$M^\Aff\simeq R(\beta)^\Aff\tens_{R(\beta)}M
\qt{for any $M\in R(\beta)\gmod$.}
$$

\smallskip
 For $M,N\in\R \gmod $, we have
$$M^\Aff \zconv  N^\Aff\simeq (M\conv N)^\Aff,$$
where
$$M^\Aff\zconv N^\Aff\seteq \Coker\bl M^\Aff\conv N^\Aff
\To[z_M-z_N]M^\Aff\conv N^\Aff\br.$$

We remark that, in this paper, we use the language of pro-objects instead of the completion in 
\cite[Section 3.1]{KKK18A} and \cite{Fu17}.

\subsubsection{Affinization in quantum affine algebra case}
Let $U_q'(\g)$ be a quantum affine algebra and let $\catC_\g$ be the category of finite-dimensional integrable $U_q'(\g)$-modules.
We embed $\catC_\g$ into $\Pro(\catC_\g)$. Note that $\Pro(\catC_\g)$ is a 
$\cor$-abelian monoidal category.
For $M\in\catC_\g$, let $M^\aff$ be the affinization of
$M$.
Recall that 
\eqn
M^\aff \simeq \cor[\zz_M^{\pm1}]\tens_\cor M
\eneqn
with the action
$$
e_i(a\tens v)=\zz_M^{\delta_{i,0}}a\tens e_iv\quad \text{for $a\in\cor[\zz_M^{\pm1}]$ and $v\in M$.}
$$
Here we use $\zz$ to distinguish from $z$ in the quiver Hecke algebra setting.
We set 
$$
M^\Aff\seteq \proolim[m] M^\aff/(\zz_M-1)^m M^\aff\in\Pro(\catC_\g).
$$
Note that there is a canonical algebra homomorphism 
$$
\cor[[\zz_M-1]]\To \End_{\Pro(\catC_\g)}(M^\Aff).
$$

\smallskip
For $M,N\in\cat$, we have
$$M^\Aff\tens_\zz N^\Aff\simeq (M \tens N)^\Aff,$$
where
$$M^\Aff\tens_\zz N^\Aff\seteq \Coker\bl M^\Aff\tens N^\Aff
\To[\zz_M-\zz_N]M^\Aff\tens N^\Aff\br.$$

For simple modules $M,N$ in $\catC_\g$, we can define
the renormalized R-matrix
$$
\Rren_{M,N}(\zz_N/\zz_M)\cl M^\Aff\tens N^\Aff \To 
N^\Aff\tens M^\Aff.
$$

\subsection{\CWS\  functor}\label{subsec:WS}
We now consider 
a duality datum $\ddD= \{ \Rt_i \}_{i\in J}$
associated with a symmetric generalized Cartan matrix $\cmC = (c_{i,j})_{i,j\in J}$.
For $i,j\in J$, we choose $\cc_{i,j}(x)\in\cor[[x]]$ such that
$$
\cc_{i,j}(x)\cc_{j,i}(-x)=1\qtq  \cc_{i,i}(0)=1. 
$$
We set 
$$
P_{ij}(u,v) \seteq \cc_{ij}(u-v) \cdot (u-v)^{d_{i,j}},
$$
where $d_{i,j}\seteq  \zero_{z=1} d_{\ms{3mu} \Rt_i, \Rt_j} (z)$.

Let $\prtl[\cmC]$ be the positive root lattice associated with $\cmC$.
For $\beta\in \prtl[\cmC] $ with $\ell=\height{\beta}$
and $\nu=(\nu_1,\ldots,\nu_\ell)\in J^\beta$, we set
$$
\hL_\nu \seteq \Rt_{\nu_1}^\Aff\tens\cdots\tens \Rt_{\nu_\ell}^\Aff
$$
and
$$
\hL(\beta)\seteq \soplus_{\nu\in J^\beta}\hL_\nu\in\Pro(\catC_\g).
$$
The algebra $R(\beta)$ acts $\hL(\beta)$ from the right as follows:
\bna
\item $e(\nu)$ is the projection to $\hL_\nu$
\item $x_k\in R(\beta)$ acts by $\log \zz_{L_{\nu_k}}$,
where $\log \zz_{L_{\nu_k}}\in\cor[[\zz_{L_{\nu_k}}-1]]\subset \End(\hL_\nu)$,
\item $e(\nu)\tau_k$ ($1\le k<\ell$) acts on $\hL_\nu$ by
$$\bc
\Rnorm_{\Rt_{\nu_k},\Rt_{\nu_{k+1}}}\circ P_{\nu_k,\nu_{k+1}}
(x_k,x_{k+1})&\text{if $\nu_k\not=\nu_{k+1}$,}\\[1ex]
(x_k-x_{k+1})^{-1}\bl \Rnorm_{\Rt_{\nu_k},\Rt_{\nu_{k+1}}}
 \circ P_{\nu_k,\nu_k}(x_k, x_{k+1}) 
-\id_{\hL_\nu}\br&\text{if $\nu_k=\nu_{k+1}$.}
\ec
$$
\ee

Note that we used $\zz-1$ instead of $\log \zz$ in \cite{KKK18A}.
 We have also relaxed the condition on $\cc_{i,i}(u)$.
We have changed the definition in order that we have Theorem~\ref{Thm: affinization} below. 

Then $\hL(\beta)$ gives the monoidal functor 
$$\hF_\ddD\cl R \gmod \to\Pro(\catC_\g)$$
defined by
$$
\hF_\ddD(M)=\hL(\beta)\tens_{R(\beta)}M \qt{for $M \in  R(\beta)  \gmod$. }
$$
It extends to
$$\hF_\ddD\cl \Pro(R)\to\Pro(\catC_\g)$$
such that $\hF_\ddD$ commutes with filtrant projective limits.

The following proposition can be proved in a similar manner to \cite{KKK18A}.
\Prop\label{Prop: Aff} 
$\hF_\ddD$ is a monoidal functor and it induces a monoidal functor
$\F_\ddD \cl \allowbreak  R \gmod\to \catC_\g$.
\enprop

Then the following theorem tells us that the functor $\hF_\ddD$ preserves affinizations.

\Th \label{Thm: affinization}
Functorially in $M\in R\gmod$, we have an isomorphism
\eqn
&&\hF_\ddD(M^\Aff)\simeq(\F_\ddD(M))^\Aff.
\eneqn
Moreover, we have
\bnum
\item the action of $z_M$ on the left term  coincides with
$\log \zz_{\F_\ddD(M)}$ on the right term,
\item for $M$, $N\in R\gmod$,
the following diagram commutes:
$$\xymatrix{
\hF_\ddD\bl M^\Aff\conv N^\Aff\br\ar[r]^-{\sim}\ar[d]&\hF_\ddD(M^\Aff)\tens\hF_\ddD(N^\Aff)
\ar[r]^-{\sim} &(\Fd(M))^\Aff\tens(\Fd(N))^\Aff\ar[d]\\
\hF\bl(M\conv N)^\Aff\br\ar[r]^-\sim&\bl\Fd(M\conv N)\br^\Aff
\ar[r]^-{\sim}&\bl\Fd(M)\tens \Fd(N)\br^\Aff\;.
}
$$
\ee
\enth
\Proof
Let us show (i).
Since $\hF_\ddD(R(\beta)^\Aff)\tens_{R(\beta)}M\simeq
\hF_\ddD(M^\Aff)$ for any $M\in R\gmod$, it is enough to show that
$$
\hF_\ddD(R(\beta)^\Aff)\simeq \bl\hF_\ddD(R(\beta))\br^\Aff
$$
compatible with the right actions of $R(\beta)$.

Set $\ell \seteq  \height{\beta}$ and $x_k \seteq \log \zz_{L_{\nu_k}}\in\End(\hL(\beta))$ for $k=1, \ldots, \ell$.
Then we have
$\hL_\nu=\cor[[x_1,\ldots,x_\ell]]\tens \Rt_\nu$.
Here we set
$$\Rt_\nu=\Rt_{\nu_1}\tens\cdots\tens \Rt_{\nu_\ell}\qtq \Rt(\beta)=\soplus_{\nu\in J^\beta}\Rt_\nu.$$
Then, $e_i$ acts on $\hL_\nu$ by
\eq \sum_{k=1}^\ell \e^{\delta_{i,0}x_k}(e_i)_k.
\label{eq:e}
\eneq
Here $\e^x$ is the exponential function and  $(e_i)_k$ denotes the action on $L_\nu$ given by
$$
\underbrace{\id\tens\cdots\tens\id}_{\txt{$(k-1)$-times}}\tens e_i\tens 
\underbrace{K_i^{-1}\tens\cdots\tens K_i^{-1}}_{\txt{$(\ell-k)$-times}}.
$$ 
Then we have
\bnum\item
$\hF_\ddD(R(\beta)^\Aff)\simeq 
\cor[[z,x_1,\ldots, x_\ell]]\tens \Rt(\beta)$.
Here $e_i$ acts by \eqref{eq:e}.
The right action of $x_k\in R(\beta)$ is given by
$x_k-z$ by \eqref{eq:right}.

\item
$\bl\hF_\ddD(R(\beta))\br^\Aff\simeq 
\cor[[z,x_1,\ldots, x_\ell]]\tens \Rt(\beta)$.
Here $e_i$ acts by
\eq  \e^{\delta_{i,0}z}\sum_{k=1}^\ell\e^{\delta_{i,0}x_k}(e_i)_k
=
 \sum_{k=1}^\ell \e^{\delta_{i,0}(x_k+z)}(e_i)_k .
\label{eq:e1}
\eneq
The right action of $x_k\in R(\beta)$ is given by
$x_k$.
\ee
Hence,
the morphism
$$f\cl \hF_\ddD(R(\beta)^\Aff)\to \bl\hF_\ddD(R(\beta))\br^\Aff$$
given by $a(z,x)\tens v\mapsto a(z, x_1+z,\ldots, x_\ell+z)\tens v$
(with $a(z,x)\in \cor[[z,x_1,\ldots, x_\ell]]$ and $v\in \Rt(\beta)$)
gives an isomorphism in $\Pro(\catC_\g)$ and the right action of $x_k\in R(\beta)$
commutes.
The compatibility of the right action of $\tau_k\in R(\beta)$ 
easily follows from the fact that $P_{i,j}(u,v)$
is a function in $u-v$.

\snoi
The second assertion (ii) is immediate. 
\QED

\subsection{\CWS \  with simply laced Cartan matrix}

{\em Hereafter, we assume that
$ \cmC= (c_{i,j})_{i,j\in J} $ is a simply laced Cartan matrix of finite type.}

Let $R_\cmC$ be the symmetric quiver Hecke algebra associated with $\cmC$. If no confusion arises, we simply write $R$ for $R_\cmC$.

Let $\ddD= \{ \Rt_i \}_{i\in J}$ be a duality datum 
associated with the Cartan matrix $\cmC $.

\Prop[\cite{KKK18A}]
We have
\bnum
\item
$\hF_\ddD$ is an exact functor and it commutes with projective limits,
\item
$\F_\ddD$ sends a simple module to a simple module or zero.
\ee
\enprop

\Lemma
Let $M\in R\gmod$ be a real simple module, and assume that $\Fd(M)$ is simple.
Then $\Fd(M)$ is also a real simple module. 
\enlemma
\Proof
Since $\Fd(M) \tens \Fd(M) \simeq \Fd(M \conv M)$ and $M\conv M$ is simple, $\Fd(M) \tens \Fd(M)$ is simple, i.e., $\Fd(M)$ is real. 
\QED

\begin{lemma} \label{Lem: same d}
Let $M,N \in R\gmod$ be simple modules such that $ \Fd(M) $ and $\Fd(N)$ are simple modules. 
Assume that one of $M$ and $N$ is real.
\bnum
\item $\de(\Fd(M), \Fd(N))\le \de(M,N)$.
\item The following conditions are equivalent:
\bna
\item
$\de(\Fd(M), \Fd(N))=\de(M,N)$.

\item
$\Fd(M\hconv N)$ and $\Fd(N\hconv M)$ are simple.
\ee
If these conditions hold, then
we have
\be[{\rm(1)}]
\item $\Fd(\rmat{M,N})\not=0$ and $\Fd(\rmat{N,M})\not=0$,
\item$\Fd(M)\hconv\Fd(N)\simeq\Fd(M\hconv N)$ and
$\Fd(N)\hconv\Fd(M)\simeq\Fd(N\hconv M)$.
\ee
\ee
\end{lemma}

\begin{proof}
Set $ z = \log \zz$  and 
$d \seteq  \de(M,N)$. 
By the definition of $\de(M,N)$, we have the following commutative diagram (up to a constant multiple): 
$$
\xymatrix@C=9ex{
M_z \conv N \ar[rr]^{\Rren_{M_z, N}} \ar@/^{2.5pc}/[rrrr]^{z^d\id}  && N \conv M_z \ar[rr]^{\Rren_{N,M_z}} && M_z \conv N .
}
$$
Applying $\hF$ to the above diagram, by Proposition \ref{Prop: Aff}
and Theorem \ref{Thm: affinization}, 
we obtain
\begin{align} \label{Eq: Rmatrix}
\xymatrix{
\Fd(M)_\zz \tens \Fd(N) \ar[rr]^{ \hF( \Rren_{M_z, N})} \ar@/^{2.5pc}/[rrrr]^{z^d\id}  && \Fd(N) \tens \Fd(M)_\zz \ar[rr]^{\hF(\Rren_{N,M_z})} && \Fd(M)_\zz \tens \Fd(N) .
}
\end{align}
Since  $ z^d \id$ is non-zero,
$ \hF(\Rren_{M_z, N})$ and $ \hF(\Rren_{N, M_z})$ are non-zero.
Note that 
\begin{align*}
&\Hom_{ \cor[\zz^{\pm1}] \tens  U_q'(\g)}(U\tens V_\zz,V_\zz\tens U)
= \cor[\zz^{\pm1}]\Rren_{U,V_\zz}, \\
& \Hom_{ \cor[\zz^{\pm1}] \tens  U_q'(\g)}(U_\zz\tens V,V\tens U_\zz)
= \cor[\zz^{\pm1}]\Rren_{U_\zz,V}.
\end{align*}
for any simple modules $U,V \in \catC_\g$ by  \cite[Proposition 9.5]{Kas02}. 
Hence, we have
$$
 \hFd(\Rren_{M_z, N}) =  z^a f(z) \Rren_{\Fd(M)_\zz, N}, \qquad \hFd(\Rren_{N, M_z}) =  z^b 
g(z) \Rren_{  N, \Fd(M)_\zz} 
$$
for some $a,b\ge 0$ and $f(z), g(z) \in \cor[[z]]^\times$. Hence it follows from $\eqref{Eq: Rmatrix}$ that 
$$
d = a+b+\de( \Fd(M), \Fd(N) ).
$$ 
Hence, we have $\de( \Fd(M), \Fd(N) )\le d$.

Moreover $d=\de( \Fd(M), \Fd(N) )$ if and only if $a=b=0$.
Since $a=b=0$ is equivalent to
$\Fd(\rmat{M,N})=\hFd(\Rren_{M_z, N})\vert_{z=0}\not=0$ and $\Fd(\rmat{N,M})=\hFd(\Rren_{N, M_z})\vert_{z=0} \not=0$.
The last two conditions conditions are equivalent to
$\Im(\Fd(\rmat{M,N}))\simeq\Fd(M\hconv N)\not\simeq0$ and 
$\Im(\Fd(\rmat{N,M}))\simeq\Fd(N\hconv M)\not\simeq0$.\qedhere
\end{proof}

\Lemma \label{Prop: LaLaLa}
Let $ \ddD= \{ \Rt_i \}_{i\in J}$ be a duality datum associated with a simply-laced finite Cartan matrix $\cmC$. 
Let $L,M,N$ be simple $R_\cmC$-modules and $S$ a simple subquotient of $M \conv N$.
Assume that $\Fd(M)$, $\Fd(N)$ and $\Fd(S)$ are simple.
\bnum
\item Assume that $\Fd(\rmat{M,L})$ and $\Fd(\rmat{N,L})$ are non-zero.
Then we have
\eqn
&&\La(\F_\ddD(M),  \F_\ddD(L) ) + \La(\F_\ddD(N),  \F_\ddD(L) )
 - \La(\F_\ddD(S),  \F_\ddD(L) )\\
&& \qquad  \qquad  \ge\La(M,L) + \La(N, L) - \La(S, L).
\eneqn
The equality holds if and only if $\Fd(\rmat{S,L})$ does not vanish.
\item 
Assume that $\Fd(\rmat{L,M})$ and $\Fd(\rmat{L,N})$ are non-zero.
\eqn
&&\La(\F_\ddD(L),  \F_\ddD(M) ) + \La(\F_\ddD(L),  \F_\ddD(N) ) 
- \La(\F_\ddD(L),  \F_\ddD(S) )\\
&&\hs{35ex}\ge\La(L, M) + \La( L, N) - \La( L, S).\nn
\eneqn
The equality holds if and only if $\Fd(\rmat{L,S})$ does not vanish.
\ee
\enlemma

\begin{proof}
Since the proof of (ii) is similar, 
we shall prove only (i).

As $S$ is a simple subquotient of  $M \conv N$, there exists a submodule $K$ of $M \conv N$ such that $S$ is a quotient of $K$.
We consider the following commutative diagram  in $R \gmod$
$$
\xymatrix{
(M \conv N)\conv L_z \ar[rr]^{\Rren_{M \circ N, L_z}} &&    L_z \conv (M \conv N) \\
K \conv L_z\akeu[1.5ex]\ar[rr]  \ar@{>->}[u] \ar@{->>}[d] && L_z \conv K \akeu[1.5ex]\ar@{>->}[u] \ar@{->>}[d] \\
S \conv L_z \ar[rr]^{z^c\Rren_{S, L_z}} && L_z \conv S
}
$$
for some $c \in \Z_{\ge 0}$. Comparing the homogeneous degrees of morphisms in the above diagram, we have 
\begin{align} \label{Eq: La in Rgmod}
2c = \La(M,L) + \La(N,L) - \La(S,L).
\end{align}
We set $ z = \log \zz $.
Applying the duality functor $\hF_\ddD$ to the above diagram, we obtain 
$$
\xymatrix@C=8ex{
( \widetilde{M} \tens \widetilde{N})\tens \widetilde{L}_{\zz} \ar[rr]^{ \hF_\ddD(\Rren_{M \circ N, L_z})} &&    \widetilde{L}_\zz \tens (\widetilde{M} \tens \widetilde{N}) \\
\widetilde{K} \tens \widetilde{L}_\zz \akeu[2ex]\ar[rr]  \ar@{>->}[u] \ar@{->>}[d] && \widetilde{L}_\zz \tens \widetilde{K}\akeu[2ex] \ar@{>->}[u] \ar@{->>}[d] \\
\widetilde{S} \tens \widetilde{L}_\zz \ar[rr]^{ z^c \hF_\ddD(\Rren_{S, L_z})} && \widetilde{L}_\zz \tens \widetilde{S}
}
$$
where $\widetilde{X} $ denotes $\F_\ddD(X)$ for a simple $R_\cmC$-module $X$.
There exist $a\in\Z_{\ge0}$ and $f(z)\in\cor[[z]]^\times$ such that 
\eqn
\hF_\ddD(\Rren_{S, L_z})=z^af(z)\Rren_{\widetilde{S},\widetilde{L}}.
\eneqn

Since $\Fd(\rmat{M,L})$ and $\Fd(\rmat{N,L})$ do not vanish, we have 
$$
\hF_\ddD( \Rren_{M, L_z} ) \equiv \Rren_{\widetilde{M}, \widetilde{L}_\zz}, \qquad 
\hF_\ddD( \Rren_{N, L_z} ) \equiv \Rren_{\widetilde{N}, \widetilde{L}_\zz} 
$$
up to a multiple of $\cor[[z]]^\times$.
The above diagram tells us that 
$$
{ {c_{\widetilde{M}, \widetilde{L}} (\zz) c_{\widetilde{N},\widetilde{L}}(\zz)}\over{ (\zz-1)^{c+a} \ c_{\widetilde{S},\widetilde{L}}(\zz)}} 
$$
is a rational function in $\zz$ which is regular and invertible at $\zz=1$.
Hence, by \cite[Lemma 3.4]{KKOP19C},  we have 
\begin{align*} 
\Deg \left( { {c_{\widetilde{M}, \widetilde{L}} (\zz) c_{\widetilde{N},\widetilde{L}}(\zz)}\over{  c_{\widetilde{S},\widetilde{L}}(\zz)}}  \right) =
2 \cdot \zero_{\zz=1} \left( { {c_{\widetilde{M}, \widetilde{L}} (\zz) c_{\widetilde{N},\widetilde{L}}(\zz)}\over{ c_{\widetilde{S},\widetilde{L}}(\zz)}}  \right)=2(c+a). 
\end{align*}
Therefore, by $\eqref{Eq: La in Rgmod}$, we conclude that 
\begin{align*}
\La(M,L) + \La(N,L) - \La(S,L)=2c
&=\Deg \left( { {c_{\widetilde{M}, \widetilde{L}} (\zz) c_{\widetilde{N},\widetilde{L}}(\zz)}\over{  c_{\widetilde{S},\widetilde{L}}(\zz)}}  \right)-2a \\*
&=\La(\widetilde{M},\widetilde{L}) + \La(\widetilde{N},\widetilde{L}) - \La(\widetilde{S},\widetilde{L})-2a. 
\end{align*}
Hence we have
\eqn
&&\La(M,L) + \La(N,L) - \La(S,L)\le
\La(\widetilde{M},\widetilde{L}) + \La(\widetilde{N},\widetilde{L}) - \La(\widetilde{S},\widetilde{L}).
\eneqn
The equality holds if and only if $a=0$
which is equivalent to 
$\Fd(\rmat{S,L})\not=0$.
\end{proof}

\subsection{Strong duality datum} \

\begin{definition} \label{Def: SDD}
A {\em strong} duality datum   $ \ddD= \{ \Rt_i \}_{i\in J}$ 
is a duality datum associated with a simply-laced finite Cartan matrix
 $\cmC=(c_{i,j})_{i,j\in J}$ 
such that all \/$\Rt_i$'s are root modules and
\begin{align*}
\de(\Rt_i, \dual^k(\Rt_j)) = - \delta(k= 0) c_{i,j}
\end{align*}
for any $k\in \Z$ and $i,j\in J$ such that $i\not=j$.
\end{definition}

In particular, we have
\eqn
&&\text{$\La(\Rt_i,\Rt_j)=-c_{i,j}$ for $i\not=j$,}\\
&&\text{$\Li(\Rt_i,\Rt_j)=-c_{i,j}$ for all $i,j\in J$.}
\eneqn 

Let $ \ddD= \{ \Rt_i \}_{i\in J}$ be a strong duality datum associated with a Cartan matrix $\cmC = (c_{i,j})_{i,j\in J}$ of finite ADE type.
Let $R_\cmC$ be the symmetric quiver Hecke algebra associated with $\cmC$. If no confusion arises, we simply write $R$ for $R_\cmC$.
We denote by 
$$
\F_\ddD \col R_\cmC\gmod \longrightarrow \catC_\g
$$
the duality functor arising from $\ddD$.
Recall that $\Fd$ sends simples to simples or zero.
However if $\ddD$ is strong, we can say more as we see below.
 
{\em Throughout this subsection, 
we assume that $\ddD$ is a strong duality datum.}

\Lemma
For $w \in \weyl$, $\La \in \wlP^+$ and $i\in J$, we have
\eqn
\eps_i\bl\dM (w\La, \La) \br&&=
\bc
- (\al_i,w\La)&\text{if $s_iw<w$,}\\
0&\text{if $s_iw>w$,}
\ec\\
\eps^*_i\bl\dM (w\La, \La) \br&&=
\bc
(\al_i,\La)&\text{if $w\ge s_i$,}\\
0&\text{otherwise,}
\ec\\
\de \bl L(i), \dM(w\La, \La)\br
&&=\bc
0&\text{if $s_i w <w$,}\\
(\al_i,w\La)&\text{if $s_iw>w$ and $w\ge s_i$,}\\
(\al_i,w\La-\La)&\text{otherwise,}
\ec
\eneqn
where $\dM(w\La, \La)$ is the determinantial module appeared in \S\,{\rm\ref{Sec: QHA}}.
\enlemma
\Proof
The equality for $\eps_i$ is proved in \cite[Proposition 10.2.4]{KKKO18}.
Let us show the equality for $\eps^*_i$.
If $w\ge s_i$, then we have
$w\La\preceq s_i\La$.
Hence $\eps_i^*\bl\dM(w\La, s_i\La)\br=0$ by the same proposition.

Since we have
$$\dM(w\La, \La)\simeq\dM(w\La, s_i\La)\hconv\dM(s_i\La, \La)
\simeq\dM(w\La, s_i\La)\hconv L(i)^{\circ (\al_i,\La)} $$
by \cite[Theorem 10.3.1]{KKKO18},
we have $\eps_i^*\bl\dM(w\La, \La)\br=(\al_i,\La)$.

Assume that $w\not\ge s_i$. Then $\La-w\la$ does not contain
$\al_i$, and hence
$\eps^*_i\bl\dM(w\La, s_i\La)\br \allowbreak =0$.

The equality for $\de$ immediately
follows from
\eq&&\de\bl L(i),M\br
=\eps_i(M)+\eps^*_i(M)+(\al_i,\wt(M))
\eneq
(\cite[Corollary 3.8]{KKOP18}).
\QED

\begin{theorem} \label{Thm: d=d}
Let $w \in \weyl$, $\La \in \wlP^+$, and set 
$$
 V_w(\La)\seteq  \Fd( \dM (w\La, \La) ).
$$ Then $V_w(\La)$ is simple and 
\begin{align*}
\de( \Rt_i, V_w(\La) ) &= \de (L(i), \dM(w\La, \La)), \\ 
\de( \dual \Rt_i, V_w(\La) ) &= \eps_i ( \dM(w\La, \La)),\\
\de(\D^2\Rt_i, V_w(\La) ) &= 0  .
\end{align*}
\end{theorem}

\begin{proof} 
First note that
\eqn
&&\text{Once we prove that $V_w(\La)$ is a simple module,
we have $\de(\D^2\Rt_i, V_w(\La) )=0$,}
\eneqn
since $\dual^2 \Rt_i$ commutes with  all $\Rt_j$'s.

\vs{.5ex}
Since $ \dM( w\La, \La ) \conv \dM( w\La', \La') \simeq \dM( w(\La+\La'), \La+\La') $ up to a grading shift for any $\La, \La \in \wlP^+$ and $ w \in \weyl$ (\cite[Proposition 4.2]{KKOP18}),
we may assume that $\La = \La_t$ for some $t\in J$. 
We may assume further that $\La$ is $w$-regular:
that is, $\ell(w)\le\ell(w')$ for any $w'\in\weyl$ such that $w'\La=w\La$.
Then, by the preceding lemma, we have
\begin{align*}
 \de (L(i), \dM(w\La, \La)) &= 
\begin{cases}
0 &  \text{ if }  s_i w < w, \\
(\alpha_i, w\La) & \text{ if } s_iw > w, 
\end{cases}
\\
 \eps_i ( \dM(w\La, \La)) &= 
\begin{cases}
- (\alpha_i, w\La) & \text{ if } s_iw < w, \\
0 &  \text{ if }  s_i w > w, 
\end{cases}
\end{align*}
if $w\La\not=\La$.

We shall argue by induction on $\ell(w)$. 

If $\ell(w)=0$, then there is noting to prove.  

If $\ell(w) = 1$, then $V_w(\La) = \Rt_t $. Then it is straightforward that the assertion is true. 

We now assume that $\ell(w) \ge 2$. 

\mnoi
\textbf{(Case 1):}\  Assume that $s_iw < w$. We set
$$
w' = s_i w, \qquad n \seteq  ( \al_i, w' \La )  \in\Z_{\ge0}.
$$
Then $\La$ is $w'$-regular and $w'\La\not=\La$.
Hence, by the induction hypothesis, we have
$$
 \de( \dual \Rt_i, V_{w'}(\La) )=0, \quad  
\de( \Rt_i, V_{w'}(\La) ) = \de(L(i), \dM(w' \La, \La))=n.
$$
Since $\dM(w\La, \La)\simeq L(i)^{\circ n}\hconv\dM(w' \La, \La)$, we have
\eq
V_w(\La) \simeq \Rt_i^{\tens n} \hconv V_{w'}(\La)
\label{eq:wn}
\eneq
by Lemma \ref{Lem: same d}. 
In particular, $V_w(\La)$ is simple.

It follows from Lemma \ref{Lem: dLX>0} that 
$$
\de(\Rt_i, V_w(\La)) = 0.
$$
Moreover, we have $\de( \dual^2 \Rt_i, V_{w'}(\La) )=0$.
Applying Lemma \ref{Lem: DLX = LX+1} (i) to the setting $L =  \Rt_i$ and $X = V_{w'}(\La)$, we obtain
$$
\de( \dual \Rt_i, V_w(\La) ) = \de( \dual \Rt_i, V_{w'}(\La) )+n =n,
$$
which gives the assertion.

\mnoi
\textbf{(Case 2):}\  Assume that $s_iw > w$. Since $\ell(w) \ge 2$, there exists $j \in J$ such that $s_j w < w$. 
We set 
$$
w' \seteq  s_j w, \qquad n \seteq  (\al_j, w'\La) \in\Z_{\ge0}.
$$
Note that $\La$ is $w'$-regular and $w'\La\not=\La$.
By \eqref{eq:wn}, we have
\begin{align*} 
V_w(\La) \simeq \Rt_j^{\tens n} \hconv V_{w'}(\La). 
\end{align*}
We set  $Z \seteq V_w(\La)$ and $Z' \seteq  V_{w'}(\La)$. 
Hence 
\eq
Z\simeq \Rt_j^{\tens n} \hconv Z'.
\label{Eq: Vw = VjVw'}
\eneq

\mnoi
(1)\  Suppose that $c_{i,j} =0 $. Then $s_i w' > w'$ since $s_is_j = s_js_i$. By the induction hypothesis, we have 
$$
 \de( \Rt_i, V_{w'}(\La)  ) = (\al_i, w'\La) = (s_j(\al_i), w \La) = (\al_i, w\La) , \qquad  \de( \dual \Rt_i, V_{w'}(\La)  ) = 0.
$$
Since $ \de( \dual^k \Rt_i, \Rt_j ) =0$ for any $k\in \Z$, it follows from $\eqref{Eq: Vw = VjVw'}$ and  Corollary~\ref{cor: normal for triple} 
that 
$$
\de(\dual^k \Rt_i , Z) = \de(\dual^k \Rt_i, Z') \qquad \text{ for any } k\in \Z.
$$
In particular, we have
\begin{align*}
\de(\Rt_i, Z) &= \de(\Rt_i, Z') = (\al_i, w\La), \\
 \de(\dual \Rt_i, Z) &= \de(\dual \Rt_i, Z') = 0.
\end{align*}

\mnoi
(2) We now assume that $c_{i,j} = -1 $. Then we have two cases: $s_i w' > w'$ and $s_i w' < w'$.
\begin{enumerate}
\item[(a)] Assume that $s_i w' > w'$. Then $\de(\dual \Rt_i, Z')=0$ by the induction hypothesis.
Hence,
by $\eqref{Eq: Vw = VjVw'}$ and Corollary \ref{Cor: de for LL' left} (i), 
we have
$$
\de(\dual \Rt_i, Z)  =0
$$ 
and 
\begin{align*}
\de(\Rt_i, Z)  &= \de(\Rt_i, Z') + n =  (\al_i, w'\La) + n \\
&= (\al_i, w'\La - n\al_j) = (\al_i, w\La). 
\end{align*}
Here the second identity follows from the induction hypothesis.

\item[(b)] Assume that $s_i w' < w'$. 
Letting $w'' \seteq  s_i w' $, we have $w = s_j s_i w'' $ and $\ell(w) = 2 + \ell(w'')$.
If $s_j w'' < w''$, then $ \ell(w) = 3+ \ell(s_j w'')  $ and $ w = s_js_is_j (s_jw'') = s_is_js_i (s_jw'')$.
This implies that $ s_i w < w$, which contradicts the assumption \textbf{(Case 2)}. Hence we have 
$ s_j w'' > w'' $, which tells us that 
$$
(\al_j, s_i w' \La) =  (\al_j, w'' \La)  \ge 0.
$$
Set $m \seteq  (\al_i, s_i w' \La ) \in\Z_{\ge0}$. 
Then, we have
\begin{align*}
(\al_j, s_i w' \La) = (\al_j, w'\La + m\alpha_i ) = n-m,
\end{align*}
which says that 
$n-m \ge 0$.

Set $Z''\seteq V_{w''}(\La)$.
By the induction hypothesis, we have
$$
\de(\Rt_i, Z' ) = 0,\quad 
\de(\D \Rt_j,Z'')=0,\quad
\de(\dual^2\Rt_j, Z'') =0.
$$
Applying Proposition~\ref{prop:XY} (iii) to the setting 
$L \seteq  \Rt_i$, $L' \seteq  \Rt_j$ and $X \seteq  Z'$,  $Y\seteq Z''$
and $k\seteq n$,
we have 
\begin{align*}
\de(\Rt_i, Z) &= \de(\Rt_i, \Rt_j^{\tens n} \hconv Z') = n-m, \\
\de(\dual \Rt_i, Z) &= 0 .
\end{align*}
Since $  (\al_i, w\La) = (\al_i, w'\La - n\al_j) = - (\al_i, s_iw'\La ) + n = n-m$, we conclude that 
$$
\de(\Rt_i, Z) = n-m =  (\al_i, w\La),
$$
which completes the proof. \qedhere
\end{enumerate}
\end{proof}

\begin{theorem}  \label{Thm: simple to simple}
Let $ \ddD= \{ \Rt_i \}_{i\in J}$ be a {\em strong} duality datum associated with a simply-laced finite Cartan matrix $\cmC$.
Then the duality functor $\F_\ddD$ sends simple modules to simple modules.
\end{theorem}
\begin{proof}
Since the duality functor $\F_\ddD$ sends a simple module to a simple module or zero, it suffices to show that 
$\F_\ddD(X)$ is non-zero for any simple module $X \in R\gmod$.

Let $w_0$ be the longest element of the Weyl group $\weyl$ of $\cmC$. Note that the category $\cC_{w_0}$ is equal to $R\gmod $.
For $i\in J$, we set $\dC_i \seteq  \dM(w_0\La_i, \La_i)$  and denote by $(\dC_i, R_{\dC_i})$ the \emph{non-degenerate braider} induced from R-matrices (\cite[Proposition 4.1]{KKOP19A}).
It is proved in \cite[Section 5]{KKOP19A} that there is a localization $\lR \seteq  R\gmod[\dC_i^{\conv -1} \mid i\in J]$ of $R\gmod$ by the braiders $\dC_i$.
Moreover $\lR$ is left rigid (\cite[Corollary 5.11]{KKOP19A}). Thus, for any simple module $X \in R\gmod$, there exists a module $Y \in R\gmod$ and $\La \in \wlP^+$ such that 
there exists a surjective homomorphism
$$
Y \conv X \twoheadrightarrow \dM(w_0\La, \La).
$$
Applying the duality functor $\F_\ddD$ to the above surjection, we have 
$$
\F_{\ddD}(Y) \tens \F_{\ddD}(X) \twoheadrightarrow \F_{\ddD}(\dM(w_0\La, \La)).
$$
Since $\F_{\ddD}(\dM(w_0\La, \La))$ is simple by Theorem \ref{Thm: d=d}, $\F_{\ddD}(X)$ does not vanish. 
\end{proof}

\Cor\label{Cor:faithful}
Let $ \ddD$ be a strong duality datum associated with a simply-laced finite Cartan matrix $\cmC$.
Then $\F_\ddD$ is faithful: i.e.,
for any non-zero morphism $f$ in $R\gmod$, $\F_\ddD(f)$ is non-zero.
\encor

\begin{theorem}  \label{Thm: invariants preserve}
Let $ \ddD= \{ \Rt_i \}_{i\in J}$ be a strong duality datum associated with a simply-laced finite Cartan matrix $\cmC= (c_{i,j})_{i,j\in J}$. 
Then, for any simple modules $M$, $N$ in $ R_\cmC\gmod$, we have 
\bnum
\item $\La(M,N) = \La( \F_\ddD(M), \F_\ddD(N) )$,
\item $\de(M,N) = \de( \F_\ddD(M), \F_\ddD(N) )$,
\item $(\wt M,\wt N) =-\Li( \F_\ddD(M), \F_\ddD(N) )$,\label{it:wt}
\item $\de\bl \dual^k\F_\ddD(M), \F_\ddD(N)\br=0$ for any $k\not=0,\pm1$,
\item \raisebox{-.7ex}{$\ba[t]{rl}\tL(M,N) =\de\bl \dual\F_\ddD(M),\F_\ddD(N)\br =\de\bl \F_\ddD(M), \dual^{-1}\F_\ddD(N)\br.\ea$}
\ee
\end{theorem}
\begin{proof}

Set $\beta \seteq  - \wt(M)$ and $\gamma \seteq  - \wt(N)$ and write $m \seteq  \height{\beta}$ and $n\seteq  \height{\gamma}$.

\snoi
(i) We shall use induction on $ m + n$. 
If $m = 0$ or $n=0$, then it is obvious. 
Hence we assume that $m$, $n\ge1$.

If $m+n=2$, then $M = L(i)$ and $N = L(j)$ for some $i,j\in J$. 
Since the assertion is obvious in the case $i=j$, 
we assume that $i\not=j$. 
Since $\F_\ddD(M) \simeq \Rt_i$ and $\F_\ddD(N) \simeq \Rt_j$, we have
\begin{align*}
\La(\Rt_i, \Rt_j) = \de(\Rt_i, \Rt_j) = -c_{i,j} = \La(L(i), L(j)). 
\end{align*} 

Suppose that $m +n\ge3$. If $m \ge 2$, then there exist simple  modules $M_1$ and $M_2$ such that
\bna
\item $\wt(M_1)\ne0 $ and $\wt(M_2)\ne0 $,
\item one of $M_1$ and $M_2$ is real,
\item $M \simeq M_1 \htens M_2$.
\ee
Hence, by Lemma~\ref{Prop: LaLaLa} together with Corollary~\ref{Cor:faithful},
we obtain
\begin{align*}
& \La(M_1, N) + \La(M_2, N) - \La(M,N) \\ 
& \quad  =  \La( \F_\ddD(M_1), \F_\ddD(N) ) + \La(\F_\ddD(M_2), \F_\ddD(N) ) - \La(\F_\ddD(M),\F_\ddD(N)).
\end{align*}
Since $\La(M_k, N) = \La( \F_\ddD(M_k), \F_\ddD(N) )$ for $k=1,2$  by the induction hypothesis, we have 
$$
\La(M, N) = \La( \F_\ddD(M), \F_\ddD(N) ).
$$

The case where $ n\ge 2$ can be similarly proved.

\mnoi
(ii) immediately follows from (i).

\snoi
(iii) There exist sequences $(i_1, \ldots, i_m)$ and $(j_1, \ldots, j_n)$ in $J$ such that 
$M$ and $N$ appear as quotients of $L(i_1)\conv \cdots \conv L(i_m)$ and  $L(j_1)\conv \cdots \conv L(j_n)$, respectively.
Note that $\beta = -\sum_{p=1}^m \al_{i_p}$ and $\gamma = -\sum_{q=1}^n \al_{j_q}$.
Since $\F_\ddD$ is exact and $ \F_\ddD(M)$ and $\F_\ddD(N)$ are simple, 
$ \F_\ddD(M)$ and $\F_\ddD(N)$ appear as quotients in $\Rt_{i_1}\tens \cdots \tens \Rt_{ i_m}$ and  $\Rt_{j_1}\tens \cdots \tens \Rt_{j_n}$, respectively.
Therefore, by \cite[Proposition 3.11]{KKOP19C}, we have 
\begin{align*}
-\Li( \F_\ddD(M), \F_\ddD(N) ) &= - \sum_{p,q} \Li( \Rt_{i_p}, \Rt_{j_q}  ) 
=   \sum_{p,q} c_{i_p,j_q}  \\
&= (\beta, \gamma).
\end{align*}

\snoi
(iv) follows from $\de(\dual^k(\Rt_i), \Rt_j)=0$ for any $i,j$ and $\vert k \vert\ge 2$.

\snoi
(v) By (i), (iii) and (iv), we have 
\begin{align*}
\La(M,N) =& \de( \F_\ddD(M), \F_\ddD(N) ) - \de( \F_\ddD(M), \dual \F_\ddD(N) ) + \de( \F_\ddD(M), \dual^{-1} \F_\ddD(N) ), \\
(\beta, \gamma) = -&\de( \F_\ddD(M), \F_\ddD(N) ) + \de( \F_\ddD(M), \dual \F_\ddD(N) ) + \de( \F_\ddD(M), \dual^{-1} \F_\ddD(N) ).
\end{align*}
Thus we have 
\begin{align*}
\tLa(M,N) &= \frac{1}{2} ( \La(M,N) + (\beta, \gamma) ) \\
&= \de( \F_\ddD(M), \dual^{-1} \F_\ddD(N) ) = \de( \dual \F_\ddD(M),  \F_\ddD(N) ).
\qedhere \end{align*}
\end{proof}

\begin{corollary} \label{Cor: ep under F}
Let $ \ddD= \{ \Rt_i \}_{i\in J}$ be a {\em strong} duality datum.
For any $i\in J$ and any simple module $M \in R_\cmC\gmod$, we have 
\bnum
\item $\eps_i (M) = \de(\dual\, \Rt_i, \F_\ddD(M) ) $,
\item $\eps_i^* (M) = \de(\dual^{-1} \Rt_i, \F_\ddD(M) ) $.
\ee
\end{corollary}
\begin{proof}
It follows from \cite[Corollary 3.8]{KKOP18} and Theorem \ref{Thm: invariants preserve} (v).
\end{proof}

\begin{corollary} \label{Cor: mono [F]}
Let $ \ddD= \{ \Rt_i \}_{i\in J}$ be a {\em strong} duality datum associated with a simply-laced finite Cartan matrix $\cmC$.
Then the duality functor $\F_\ddD$ induces an injective ring homomorphism 
$$\KRc\monoto K(\catC_\g  ),$$
where $ \KRc$ is the specialization of the $K( \RC\gmod)$ at $q=1$.
\end{corollary}
\begin{proof}
Thanks to Theorem \ref{Thm: simple to simple}, it is enough to show that $\F_\ddD(M) \not\simeq \F_\ddD(N)$ for any non-isomorphic pair of simple $R$-modules $M$ and $N$. 
Let $M$ and $N$ be simple $R$-modules such that  
$$
\F_\ddD(M)\simeq \F_\ddD(N).
$$
We set $\beta\seteq- \wt(M)$ and $\gamma\seteq - \wt(N)$. 
We shall show $M\simeq N$. 

We first assume that $  \F_\ddD(M) = \F_\ddD(N) = \one $. 
Then we have 
$(\beta,\beta)=  -\Li(M,M) \allowbreak =0$, which implies  $\beta=0$. Hence $M\simeq\one$.
Similarly, we have $N\simeq\one$.

\bigskip
We now assume that $ \F_\ddD(M) \simeq \F_\ddD(N) \not\simeq \one $. Since $M \not\simeq \one $, there exists $i\in J$ such that $\eps_i(M) > 0$.
By Corollary \ref{Cor: ep under F}, we have 
\begin{align*}
\eps_i(M) & =  \La(\dual \Rt_i, \F_\ddD(M))   = \La(\dual \Rt_i, \F_\ddD(N))  = \eps_i(N), 
\end{align*}
which tells us that 
$ \te_i(M) \ne 0 $ and $\te_i(N) \ne 0$. 
Setting $M' \seteq  \te_i(M)$ and $N' \seteq  \te_i(N)$, we have 
$$
\Rt_i \hconv \F_\ddD(M') \simeq \F_\ddD(M) \simeq \F_\ddD(N) \simeq \Rt_i \hconv \F_\ddD(N'),
$$
which implies that $\F_\ddD(M') \simeq \F_\ddD(N') $ by Lemma \ref{Lem: crystal for real}.
Thus, by the standard induction argument, we conclude that 
$$
\te_i(M) = M' \simeq N' = \tei(N),
$$ 
which yields that $M \simeq N$. 
\end{proof}

\Lemma\label{lem:root}
Let $M$ be a real simple module in $R\gmod$.
Then $\F_\ddD(M)$ is a root module if and only if $\wt(M)$
is a root of $\gf$. 
\enlemma
\Proof
Set $V=\F_\ddD(M)$. Then we have
$\de(\D^k V,V)=0$ for $k\not=\pm1$.
Hence we have
$\bl\wt(M),\wt(M)\br=
-\Li(V,V)=2\de(\D V,V)$.
Therefore we have
$$\text{$V$ is a root module}\Leftrightarrow\de(\D V,V)=1
\Leftrightarrow
\bl\wt(M),\wt(M)\br=2\Leftrightarrow\text{$\wt(M)$ is a root}.
$$
\QED

\vskip 2em

\section{Strong duality datum and affine cuspidal modules}  \label{Sec: SD and ACM} \

\subsection{Unmixed pairs} \

The notion of an unmixed pair of modules over quiver Hecke algebras
has an analogue for modules over quantum affine algebras.
\begin{definition} \label{Def: unmixed}
Let $(M,N)$ be an ordered pair of simple modules in $\catC_\g$. We call it \emph{unmixed} if 
$$
\de(\dual M, N)=0,
$$
and \emph{strongly unmixed} if 
$$
\de(\dual^k M, N)=0 \qquad \text{ for any } k\in \Z_{\ge 1}.
$$
\end{definition}

\begin{lemma} \label{Lem: unmixed Li=La}
Let $M$ and $N$ be simple modules in $\catC_\g$. If $(M,N)$ is strongly unmixed, then 
$$
\Li(M,N) = \La(M,N).
$$
\end{lemma}
\begin{proof}
It follows from Definition \ref{Def: unmixed} and Proposition \ref{Prop: La and d} that 
\eqn
\La(M,N)&&= \sum_{k \in \Z} (-1)^{k+\delta(k>0)} \de(\dual^{k}M,N) 
 = \sum_{k \in \Z} (-1)^{k} \de(\dual^{k}M,N) = \Li(M,N).
\eneqn
\end{proof}

\begin{lemma} \label{Lem: unmixed and normal}
Let $ L_1, \ldots, L_r$ be real simple modules in $\catC_\g$ for $r\in \Z_{>1}$. If $(L_a, L_b)$ is  unmixed for any $a < b$, then 
$(L_1, \ldots, L_r)$ is normal.
\end{lemma}
\begin{proof}
We shall argue by induction on $r$. Since it is obvious when $r=2$, we assume that $r>2$.
By the induction hypothesis, $(L_1,\ldots, L_{r-1})$ is normal.
Set $X=\hd(L_2\tens\cdots\tens L_{r-1})$.
Then Lemma~\ref{Lem: normal} implies that 
$$\La(L_1,X)=\sum_{k=2}^{r-1}\La(L_1,L_k).$$
Since $(L_1,L_r)$ is unmixed, Lemma~\ref{Lem: normal for 3} implies that $(L_1,X,L_r)$ is normal.
Hence we have
$$\La\bl L_1, \hd(L_2\tens\cdots \tens L_r)\br=\La(L_1,X\hconv L_r)=\La(L_1,X)+\La(L_1,L_r),$$
which implies that
$$\La\bl L_1, \hd(L_2\tens\cdots \tens L_r)\br=\sum_{k=2}^{r}\La(L_1,L_k).$$
Since $(L_{2}, \ldots, L_r)$  is normal, Lemma~\ref{Lem: normal} implies that $(L_{1}, \ldots, L_r)$  is normal. 
\end{proof}

\smallskip

\subsection{Affine cuspidal modules}\label{subsec:cusp}
Let  $ \ddD= \{ \Rt_i \}_{i\in J}$ be a strong duality datum 
in $\catCO$ associated with a simply-laced finite Cartan matrix 
$\cmC =(c_{i,j})_{i,j\in J}$. 
Let $R_\cmC$ be the symmetric quiver Hecke algebra associated with $\cmC$.

The category $\catCD$ is defined to be  the smallest full subcategory of $\catCO$ such that
\bna
\item   it contains $\F_\ddD( L )$ for any simple  $R_{\cmC}$-module $L$, 
\item  it is stable by taking subquotients, extensions, and tensor products. 
\ee

Since $ \de(\dual^k \Rt_i, \Rt_j)=0 $ for any $ i,j\in J $ and $k \ge 2$, it follows from Theorem \ref{Thm: invariants preserve} that   
\begin{align} \label{Eq: d=0 for Dk}
\de( \dual^k M, N) = 0\qquad \text{for any simple module $M,N \in \catCD$ and $ k \ge 2.$}
\end{align}

For $k\in \Z$, let $\dual^k (\catCD)$  be the full subcategory of $\catCO$ whose objects are $\dual^kM$ for all $M\in \catCD$.

\begin{prop}
Let $k\in \Z$ with $k \ne 0$. If a simple module $M$ is contained in $\catCD \cap \dual^k(\catCD)$, then $M \simeq \trivial$.
\end{prop}
\begin{proof}
We may assume $k>0$ without loss of generality. 
Let $M $ be a simple module in $ \catCD \cap \dual^k (\catCD)$.
By Theorem \ref{Thm: simple to simple}, there exists a simple module $V \in R_\cmC\gmod$ such that $\F_\ddD(V) \simeq M$.
By Corollary \ref{Cor: ep under F} and Theorem \ref{Thm: invariants preserve} (iv), for any $i\in J$, we have
\begin{align*}
\eps_i^* (V) = \de(\dual^{-1} \Rt_i, M ) = \de(\Rt_i, \dual M) = 0.
\end{align*}
Thus $V$ should be in $R_\cmC(0)\gmod$, which says that $V\simeq \trivial$.
\end{proof}

\begin{lemma} \label{Lem: unmixed to unmixed}
Let $M,N$ be simple modules in $\RC\gmod$. If $(M,N)$ is unmixed, then $\bl \F_\ddD(M), \F_\ddD(N)\br$ is strongly unmixed. 
\end{lemma}
\begin{proof}
By $\eqref{Eq: d=0 for Dk}$, we know that $\de( \dual^k  \F_\ddD(M), \F_\ddD(N) )=0$ for $k \ge 2$.
It follows from \cite[Proposition 2.12]{KKOP18} that $\La( M, N )= -(\wt(M), \wt(N)) $, i.e., $\tLa( M, N )=0$. 
Thus, by Theorem \ref{Thm: invariants preserve} (v), we obtain 
$$
\de( \dual \F_\ddD(M), \F_\ddD(N) ) = \tLa( M, N )=0,
$$
which completes the proof.
\end{proof}

Let $\g_\cmC$ be the simple Lie algebra associated with 
$\cmC$.
Let $\prDC$ be the set of positive roots of $\g_\cmC$ and let $\weylC $ be the Weyl group associated with $\g_\cmC$. 
Let $w_0$ be the longest element of $\weylC$, and $\ell$ denotes the length of $w_0$.
We choose an arbitrary reduced expression $\rxw = s_{i_1} s_{i_2} \cdots s_{i_\ell}$ of the longest element $w_0$ of $\weylC$.
We extend $\st{i_k}_{1\le k\le \ell}$ to $\st{i_k}_{k\in\Z}$ by 
\eq
&&i_{k+\ell}=(i_k)^* \qt{for any $k\in\Z$.}\label{def:ik}
\eneq
(Recall that, for $i\in J$, $i^*$ is a unique element of $J$ such that
$\al_{i^*}=-w_0\al_i$.)

We can easily see that $s_{i_{a+1}}\cdots s_{i_{a+\ell}}$ is also a reduced expression
of $w_0$ for any $a\in\Z$.
Let 
$$
\{ \Cp_k \}_{k=1, \ldots, \ell} \subset \RC\gmod
$$ 
be the cuspidal modules associated with the reduced expression $\rxw$.
Under the categorification, the cuspidal module $\Cp_k$ corresponds to the dual PBW vector $E^*(\beta_k)$ corresponding to $\beta_k\seteq  s_{i_1} \cdots s_{i_{k-1}} (\alpha_{i_k}) \in \prDC $ for $k=1, \ldots, \ell$
(see Section \ref{Sec: QHA}).

\smallskip

We now introduce the notion of affine cuspidal modules for quantum affine algebras.

\begin{definition}\label{def:cusp}
We define a sequence of simple $U_q'(\g)$-modules $\{ \cuspS_k \}_{ k\in \Z }$ in $\catC_\g $ as follows:
\bna
\item $\cuspS_k = \F_\ddD(\Cp_k)$ for any $k=1, \ldots, \ell$, 
and we extend its definition to all $k\in\Z$ by 
\item $\cuspS_{k+\ell} = \dual( \cuspS_k )$ for any $k\in \Z$.
\ee
The modules $\cuspS_k$ $(k\in \Z)$ are called the \emph{affine cuspidal modules} corresponding to $\ddD$ and $\rxw$.
\end{definition}

\begin{prop} \label{Prop: cusp} 
The affine cuspidal modules satisfy the following properties.
\bnum
\item $\cuspS_a$ is a root module for any $a\in\Z$.
\item
For any $a,b \in \Z $ with $a > b$, the pair $( \cuspS_a, \cuspS_b )$ is strongly unmixed.
\item Let $ k_1 > \cdots > k_t $ be decreasing integers and $ (a_1, \ldots, a_t) \in \Z_{ \ge0 }^t  $. Then 
\bna
\item the sequence  
$ ( \cuspS_{k_1}^{\tens a_1}, \ldots, \cuspS_{k_t}^{\tens a_t} ) $ is normal,
\item the head of the tensor product $ \cuspS_{k_1}^{\tens a_1} \tens \cdots \tens \cuspS_{k_t}^{\tens a_t} $  is simple.
\ee
 
\ee
\end{prop}

\begin{proof}

(i) follows immediately from Lemma~\ref{lem:root}.

\snoi
(ii) Without loss of generality, we may assume that $ 1 \le b \le \ell$.
 We write $a = \ell \cdot t + r $ for some $t \in \Z_{\ge 0}$ and $1 \le r \le \ell$. By the definition, we have $\cuspS_a = \dual^t \cuspS_r$.
If $t \ge 1$, then we have
$$
\de(\dual^{k} \cuspS_a, \cuspS_b   ) = \de(\dual^{k+t} \cuspS_r, \cuspS_b   ) = 0 \qquad \text{ for any $k\ge 1$,}
$$
by $\eqref{Eq: d=0 for Dk}$.

Suppose that  $t=0$. As $\ell \ge a > b \ge 1$,  the pair $ ( \Cp_a, \Cp_b)$ is unmixed. Thus Lemma \ref{Lem: unmixed to unmixed} says that $(\cuspS_a, \cuspS_b)$ is strongly unmixed.

\snoi
(iii) follows from Lemma \ref{Lem: unmixed and normal} and
Lemma~\ref{Lem: normal head socle}.  
\end{proof}

\begin{example} \label{Ex: ex1} 

Let $U_q'(\g)$ be the quantum affine algebra of affine type $A_2^{(1)}$, and 
let $\catCO$ be the Hernandez-Leclerc category corresponding to $\sigZ = \{ (1, (-q)^{2k}), (2, (-q)^{2k+1}) \mid k\in \Z  \}$.   
For $i\in I_0$ and $m\in \Z_{>0}$, we denote the \emph{Kirillov-Reshetikhin module} by 
$$ 
V(i^m) := \hd \left( V(\varpi_i)_{ {(-q)}^{m-1}} \tens V(\varpi_i)_{ {(-q)}^{m-3}} \tens \cdots \tens V(\varpi_i)_{{(-q)}^{-m+1}}  \right).
$$
We simply write $V(i)$ instead of $V(i^1)$, which is the $i$th fundamental module $V(\varpi_i)$.

Let $\Rt_1 := V(1) $ and $\Rt_2 := V(1)_{(-q)^2} $, and define $\ddD := \{ \Rt_1, \Rt_2 \} \subset \catCO$. Then $\ddD$ is a strong duality datum (see \cite[Section 4.1]{KKK18A}).
Let $\cmC$ be the Cartan matrix of finite type $A_2$.
Then we have the duality functor $\F_\ddD \col R_\cmC\gmod \longrightarrow \catCO$.

\bni
\item We choose a reduced expression $\rxw = s_1s_2s_1$. Then we have 
$$
\beta_1 := \al_1, \qquad  \beta_2 := s_1(\al_2)= \al_1 + \al_2, \qquad \beta_3 = s_1s_2(\al_1) = \al_2,
$$
and the affine cuspidal modules corresponding to $\ddD$ and $\rxw$ are given as follows: 
\begin{align*}
	\cuspS_1 &= \F_\ddD(L(1)) = \Rt_1 = V(1), \\
	\cuspS_2 &= \F_\ddD(L(1) \hconv L(2)) = \Rt_1 \htens \Rt_2 = V(1) \htens V(1)_{(-q)^2} = V(2)_{-q}, \\
	\cuspS_3 &= \F_\ddD(L(2)) = \Rt_2 = V(1)_{(-q)^2},
\end{align*}
and $\cuspS_{ k+3} = \dual (\cuspS_k)$ for $k\in \Z$.
Here $L(i)$ be the self-dual 1-dimensional  simple $R(\al_i)$-module.
It is easy to see that the set $\{ \cuspS_k \mid k\in \Z \}$ of all affine cuspidal modules is equal to the set of all fundamental modules in $\catCO$.

\item We choose another reduced expression $\rxw' = s_2s_1s_2$. Then we have 
$$
\beta_1' := \al_2, \qquad  \beta_2' := s_2(\al_1)= \al_1 + \al_2, \qquad \beta_3' = s_2s_1(\al_2) = \al_1,
$$
and the affine cuspidal modules corresponding to $\ddD$ and $\rxw'$ are given as follows: 
\begin{align*}
	\cuspS_1' &= \F_\ddD(L(2)) = \Rt_2 = V(1)_{(-q)^2}, \\
	\cuspS_2' &= \F_\ddD(L(2) \hconv L(1)) = \Rt_2 \htens \Rt_1 = V(1)_{(-q)^2} \htens V(1) = V(1^2)_{-q}, \\
	\cuspS_3' &= \F_\ddD(L(1)) = \Rt_1 = V(1),
\end{align*}
and $\cuspS_{ k+3}' = \dual (\cuspS_k')$ for $k\in \Z$.  
Note that the affine cuspidal modules $\cuspS_{2+3t}'$ ($t\in \Z$) are not fundamental modules. 
\ee

\end{example}

\subsection{Reflections} \

For any $k \in J$, we set 
\begin{align} \label{Def: refl}
 \Refl_k (\ddD) \seteq  \{ \Refl_k (\Rt_i)  \}_{i\in J}\qtq
 \Refl_k^{-1} (\ddD) \seteq  \{ \Refl^{-1}_k (\Rt_i)  \}_{i\in J},
\end{align}
where 
$$
\Refl_k (\Rt_i) \seteq  
\begin{cases}
 \dual \Rt_i  & \text{ if } i=k, \\
\Rt_k \htens \Rt_i &  \text{ if }  c_{i,k}=-1, \\
\Rt_i &  \text{ if } c_{i,k}=0,
\end{cases}
\qtq
\Refl^{-1}_k (\Rt_i) \seteq  
\begin{cases}
 \dual^{-1} \Rt_i  & \text{ if } i=k, \\
\Rt_i \htens \Rt_k &  \text{ if }  c_{i,k}=-1, \\
\Rt_i &  \text{ if } c_{i,k}=0.
\end{cases}
$$
It is easy to see that $\Refl_k \circ \Refl^{-1}_k (\ddD) = \ddD $ and $ \Refl^{-1}_k \circ \Refl_k (\ddD) = \ddD  $ for any $k\in J$.

\begin{prop} \label{Prop: strong to strong}
Let $k\in J$.
\bnum
\item For any $i\in J$, $\Refl_k(\Rt_i)$ and $\Refl_k^{-1}(\Rt_i)$ are root modules.
\item $\Refl_k(\ddD)$ and $\Refl_k^{-1}(\ddD)$ are strong duality data
associated with the Cartan matrix $\cmC$.  
\ee
\end{prop}
\begin{proof}
We shall focus on proving the case for $\Refl_k$ since the case for $\Refl^{-1}_k$ can be proved in a similar manner.

\smallskip

Set $\Rt_i' \seteq  \Refl_k(\Rt_i)$ for $i\in J$. For $i,j \in J$, we write $i \sim j$ if $c_{i,j}=-1$ and $i \not\sim j$ if $c_{i,j}=0$.
Note that,  for real simple modules $L$, $M$ and $N$, Lemma \ref{Lem: normal for 3} says that if one of the following conditions
\begin{itemize}
\item $\de(L,M)=0$,
\item $\de(M,N)=0$,
\item $ \de(L, \dual^{-1}N) = \de(\dual L, N)=0$
\end{itemize}
holds, then  
$$
\La(L, M\htens N) = \La(L, M) + \La(L, N), \quad \La(L \htens M, N) = \La(L, N) + \La(M, N),
$$
 which will be used several times in the proof.

\mnoi
(i) follows from Lemma~\ref{lem:real}.

\mnoi
(ii) Thanks to (i), it suffices to prove that 
$$
\de( \dual^t \Rt_i', \Rt_j') = - \delta(t=0) c_{i,j} \qquad \text{ for $t\in \Z$ and $ i\ne j$.} 
$$ 
Let $i,j\in J$ with $i \ne j$. We shall prove it case by case.

\snoi
\textbf{(Case 1): }  If $i \not \sim k$ and $j \not \sim k$, then 
$$
\de(\dual^t\Rt_i', \Rt_j') = \de(\dual^t\Rt_i, \Rt_j) = -\delta(t=0) c_{i,j} \qquad \text{ for }t\in \Z.
$$

\snoi
\textbf{(Case 2): } If $i \not \sim k$ and $j  = k$, then $c_{i,j}=0$ and 
\begin{align*}
\de( \dual^t \Rt_i', \Rt_j' ) = \de( \dual^t \Rt_i, \dual \Rt_j ) = \de( \dual^{t-1} \Rt_i,  \Rt_j ) = 0 = -\delta(t=0) c_{i,j}.
\end{align*}

\snoi
\textbf{(Case 3): }  Suppose that $i \not \sim k$ and $j  \sim k$. Then we have
$$
\de(\dual^t \Rt_i', \Rt_j' ) = \de(\dual^t \Rt_i , \Rt_k \htens \Rt_j ).
$$
Then we have
\begin{align*}
\La( \dual^t \Rt_i, \Rt_k \htens \Rt_j) &= \La( \dual^t \Rt_i, \Rt_k) + \La( \dual^t \Rt_i,  \Rt_j), \\
\La( \Rt_k \htens \Rt_j, \dual^t \Rt_i) &= \La( \Rt_k, \dual^t \Rt_i) + \La(  \Rt_j, \dual^t \Rt_i),
\end{align*}
where the first equality follows from $\de( \dual^t \Rt_i, \Rt_k )=0 $
and the second  from
$\de(\dual\Rt_k, \dual^t \Rt_i)=0 $. 
Hence we obtain
\begin{align*}
\de( \dual^t \Rt_i', \Rt_j') &= \de( \dual^t \Rt_i, \Rt_k) + \de( \dual^t \Rt_i,  \Rt_j) = \de( \dual^t \Rt_i,  \Rt_j) \\ 
& = - \delta(t=0)c_{i,j}.
\end{align*}

\mnoi
\textbf{(Case 4): }  Suppose that $i \sim k$ and $j  \sim k$. Then, by Lemma \ref{Lem: dual head}, 
$$
\de(\dual^t \Rt_i', \Rt_j' ) = \de( \dual^t( \Rt_k \htens \Rt_i ), \Rt_k \htens \Rt_j ) = \de( \dual^t \Rt_k \htens \dual^t \Rt_i , \Rt_k \htens \Rt_j ).
$$
Since $\cmC$ is of finite type, we have $ c_{i,j}= 0$, i.e., $\de( \dual^t \Rt_i, \Rt_j )=0$ for any $t \in \Z$.  

\bna
\item If $t\ne 0,\pm1$, then 
$$
\de(\dual^t \Rt_i', \Rt_j' ) =  \de( \dual^t \Rt_k \htens \dual^t \Rt_i, \Rt_k \htens \Rt_j ) = 0
$$
since $\de( \dual^t (\Rt_a), \Rt_b )=0$ for $a,b=i,j,k$ by Theorem \ref{Thm: invariants preserve} (iv).

\item Suppose that $t=0$. Then 
\begin{align*}
\La( \Rt_k \htens \Rt_i, \Rt_k \htens \Rt_j ) &= \La( \Rt_k \htens \Rt_i, \Rt_k ) + \La( \Rt_k \htens \Rt_i, \Rt_j ) \\
 &= - \La( \Rt_k, \Rt_k \htens \Rt_i ) + \La( \Rt_k \htens \Rt_i, \Rt_j ) \\
&= - \La( \Rt_k, \Rt_i ) + \La( \Rt_k , \Rt_j ) + \La( \Rt_i, \Rt_j ).
\end{align*}
Here the first and second identities follow from $\de( \Rt_k, \Rt_k \htens \Rt_i )=0$ by Lemma \ref{Lem: two root modules} and the third follows from $\de(\Rt_i, \Rt_j)=0$. 
Exchanging $i$ and $j$, we have 
\begin{align*}
\La( \Rt_k \htens \Rt_j, \Rt_k \htens \Rt_i ) = - \La( \Rt_k, \Rt_j ) + \La( \Rt_k , \Rt_i ) + \La( \Rt_j, \Rt_i ),
\end{align*}
which tells us that 
$$
\de( \Rt_i', \Rt_j' ) = \de( \Rt_k \htens \Rt_j, \Rt_k \htens \Rt_i ) =  \de( \Rt_i, \Rt_j)  = -c_{i,j}.
$$

\item  Suppose that $t= \pm 1$. 
We have 
\eqn
&&\Li (\Rt'_i,\Rt'_j)=\Li(\Rt_k \htens \Rt_i, \Rt_k \htens \Rt_j)\\
&&\hs{5ex}=\Li(\Rt_k,\Rt_k)+\Li(\Rt_k,\Rt_j)+\Li(\Rt_i,\Rt_k)+\Li(\Rt_i,\Rt_j)
=(-2)+1+1+0=0.\eneqn
Hence we have
$$
0=\sum_{t\in\Z}(-1)^t\de(\dual^t\Rt'_i,\Rt'_j)
=-\de(\dual\Rt'_i,\Rt'_j)-\de(\dual^{-1}\Rt'_i,\Rt'_j).$$
Hence $\de(\dual\Rt'_i,\Rt'_j)=\de(\dual^{-1}\Rt'_i,\Rt'_j)=0$.

\ee

\snoi
\textbf{(Case 5): } Suppose that $i \sim k$ and $j  = k$. Then, we have
$$
\de(  \dual^t \Rt_j', \Rt_i' ) = \de(  \dual^{t+1} \Rt_k, \Rt_k\htens \Rt_i ) \quad \text{ for } t\in \Z,
$$
which is equal to $\delta(t=0)$ by Lemma \ref{Lem: two root modules}. 
\end{proof}

\Prop\label{prop:ReflCusp}
Let $\st{\cuspS_k}_{k\in\Z}$ be the sequence of the affine cuspidal modules 
corresponding to $\ddD$
and a reduced expression $\rxw=s_{i_1}\cdots s_{i_{\ell}}$ of $w_0$.
Set $\cuspS'_k=\cuspS_{k+1}$ for $k\in \Z$.
Then $\st{\cuspS'_k}_{k\in\Z}$ is the affine cuspidal modules corresponding to
$\Refl_{i_1}\ddD$
and the reduced expression $\rxw'=s_{i_2}\cdots s_{i_{\ell+1}}$ \ro see \eqref{def:ik}\rf. 
\enprop
\Proof

Set $i=i_1$.
We denote by $\prec$, $ \{ \beta_k \}_{k=1, \ldots, \ell}$ and $\{ \Cp_k \}_{k=1, \ldots, \ell}$ the convex order, the ordered set of positive roots and the cuspidal modules in $R_\cmC\gmod$ corresponding to $\rxw$ as in Section \ref{Sec: QHA}.
Similarly, we write $\prec'$, $ \{ \beta_k' \}_{k=1, \ldots, \ell} $ and $\{ \Cp_k' \}_{k=1, \ldots, \ell} $ for the ones corresponding to $\rxw'$. 
It is enough to show that
\eqn
\F_{\Refl_i(\ddD)}(\Cp_k')\simeq\cuspS_{k+1}\qt{for $1\le k\le \ell$.}
\eneqn

It is easy to see that 
\begin{align}
&\text{$\bullet$ $ \beta_{k+1} = s_i \beta_k'$ for $k=1, \ldots, \ell-1$,} \label{Eq: b=sib'} \\
&\text{$\bullet$ $ \Cp_{k+1} \simeq \rT_i ( \Cp_k')$ for $k=1, \ldots, \ell-1$,} \label{Eq: S=siS'} \\
&\text{$\bullet$ $\al_i$ is smallest (resp.\ largest) with respect to $\prec$ (resp.\ $\prec'$).}  \label{Eq: al in w0 w0'}
\end{align}

It follows from $\eqref{Eq: al in w0 w0'}$ that $ \Cp_1 \simeq L(i)\simeq \Cp_\ell' $. Thus we have
\begin{align} \label{Eq: refl of cusp i}
\F_{\Refl_i(\ddD)} ( \Cp_\ell' ) \simeq  \dual \Rt_i \simeq \dual ( \F_{ \ddD} ( \Cp_1 ) )
=\cuspS_{\ell+1}.
\end{align}

It remains to prove that
\begin{align} \label{Eq: refl on cusp}
\text{ $ \F_{\Refl_i(\ddD)} ( \Cp_k' ) \simeq   \F_{ \ddD} ( \Cp_{k+1} )  $ for $k=1, \ldots, \ell-1$.}
\end{align}
We shall use induction on $\height{\beta_k'}$.

\smallskip 

If $\height{\beta_k'}=1$, then $\beta_k' = \al_j$ for some $j\in J$. Note that $j \ne i$ because $k < \ell$. 
Thus, we have $\beta_{k+1} = s_i ( \beta_k' ) = s_i(\al_j)$ and 
$$
\Cp_{k+1} = 
\begin{cases}
L(i) \hconv L(j)& \text{if $c_{i,j}=-1$,} \\
L(j) & \text{otherwise.}
\end{cases}
$$
By the definition of $\Refl_i$, we have 
$$
\F_{\Refl_i(\ddD)}( \Cp_k') \simeq \F_{\Refl_i(\ddD)}( L(j)) \simeq\Refl_i(\Rt_j)\simeq \F_{\ddD}(\Cp_{k+1}).
$$

\smallskip 
Suppose that $\height{\beta_k'}>1$. We take a minimal pair $(\beta_a', \beta_b')$ of $\beta_k'$ with respect to $\prec'$. 
It follows from $\eqref{Eq: minimal for QHA}$ that
\begin{align} \label{Eq: minimal pair CP'}
\Cp_a' \hconv \Cp_b' \simeq \Cp_k'.
\end{align}

\mnoi
 \textbf{ (Case 1)}: Suppose that $ b \ne \ell$. Applying $\rT_i$ to $\eqref{Eq: minimal pair CP'}$, it follows from $\eqref{Eq: S=siS'}$ that 
\begin{align} \label{Eq: minimal pair CP k+1}
\Cp_{a+1} \hconv \Cp_{b+1} \simeq \Cp_{k+1}.
\end{align}
Applying $\F_\ddD$ to the above isomorphism $\eqref{Eq: minimal pair CP k+1}$
and using the induction hypothesis, we have
\begin{align*}
\F_\ddD( \Cp_{k+1} ) &\simeq \F_\ddD( \Cp_{a+1} ) \htens \F_\ddD( \Cp_{b+1} ) \\
& \simeq \F_{\Refl_i(\ddD)} ( \Cp_a' ) \htens\F_{\Refl_i(\ddD)} ( \Cp_b' )    \\
& \simeq \F_{\Refl_i(\ddD)} ( \Cp_k' ).
\end{align*}

\mnoi
 \textbf{ (Case 2)}: Suppose that $ b=\ell$. Since $\Cp_b' = L(i)$, by applying $\F_{\Refl_i(\ddD)}$ to $\eqref{Eq: minimal pair CP'}$, we have 
\begin{align} \label{Eq: b=ell 1}
\F_{\Refl_i(\ddD)}(\Cp_a') \htens \dual \Rt_i  \simeq \F_{\Refl_i(\ddD)}(\Cp_k').
\end{align}

On the other hand, it follows from $\tf_i^*( \Cp_a') =  \Cp_a' \hconv L(i) \simeq \Cp_k'$ that 
$$
\eps_i^*(\Cp_a') + 1 = \eps_i^*(\Cp_k'), \qquad \vphi_i^*(\Cp_a')  = \vphi_i^*(\Cp_k')+1.
$$ 
Thus, by $\eqref{Eq: Saito refl}$ and $\eqref{Eq: S=siS'}$, we have  
\eqn
 \Cp_{a+1}  =  \rT_i(\Cp_a') &&\simeq
\tf_i^{\vphi_i^*(\Cp_k')+1}\;\te^*_i{}^{\eps_i^*(\Cp_a')+1}\tf^*_i\Cp'_a\\
&&\simeq
L(i)\hconv\bl \tf_i^{\vphi_i^*(\Cp_k')}\;\te^*_i{}^{\eps_i^*(\Cp_k')}\Cp'_k\br
\simeq
 L(i)\hconv \rT_i (\Cp_k') = L(i) \hconv \Cp_{k+1},
\eneqn  
which implies that 
\begin{align}
\F_\ddD (\Cp_{a+1}) \simeq  \Rt_i \htens \F_\ddD (\Cp_{k+1}).
\label{eq:VDD}
\end{align}
Then, we have 
\eqn
&&\ba{rll}
\F_\ddD (\Cp_{k+1}) & \simeq \F_\ddD (\Cp_{a+1}) \htens \dual  \Rt_i 
&\text{by \eqref{eq:VDD} and Lemma~\ref{Lem: MNDM}}
\\
& \simeq \F_{\Refl_i(\ddD)} (\Cp_{a}') \htens \dual  \Rt_i
\akew[5ex]&\text{by the induction hypothesis}
 \\
& \simeq  \F_{\Refl_i(\ddD)}(\Cp_k')&\text{by \eqref{Eq: b=ell 1},}
\ea
\eneqn
which completes the proof for $\eqref{Eq: refl on cusp}$.
\QED

\begin{prop} \label{Prop: CsIDCD}
Let $i\in J$, and let $S$ be a simple module in $\catC_\g$. 
\bnum
\item
 The following conditions are equivalent:
\bna
\item $S\in\F_\ddD( \cC_{*, s_i}  )$,
\item $S\in\catC_\ddD$ and $\de(\D^{-1} \Rt_i, S)=0$,
\item $S\in\catC_{\Refl_i(\ddD)} \cap \catCD$.
\ee
\item
 The following conditions are equivalent:
\bna
\item $S\in\F_\ddD( \cC_{s_iw_0}  )$,
\item $S\in\catC_\ddD$ and $\de(\D\Rt_i, S)=0$,
\item $S\in\catC_{\Refl_i^{-1}(\ddD)} \cap \catCD$.
\ee
\ee
Here, $\cC_{*, s_i} $ and $\cC_{ s_iw_0} $ are the subcategories of $R_\cmC\gmod$ appeared in {\rm\S\;\ref{Sec: QHA}}.
\end{prop}
\begin{proof}
We shall focus on proving the assertion (i) since the assertion (ii) can be proved in a similar manner.

\smallskip
Let us take a reduced expression
$\rxw=s_{i_1}\cdots s_{i_{\ell}}$ of $w_0$ such that $i_1=i$.
Let $\{ \Cp_k \}_{k=1, \ldots, \ell}$ be the cuspidal modules in $R_\cmC\gmod$ 
corresponding to $\rxw$.
Let $\st{\cuspS_k}_{k\in\Z}$ be the affine cuspidal modules 
corresponding to $\ddD$
and $\rxw$.
Set $\cuspS'_k =\cuspS_{k+1}$ or $k\in \Z$.
 Then $\st{\cuspS'_k}_{k\in\Z}$ is the cuspidal modules corresponding to
$\Refl_{i_1}\ddD$
and $\rxw'=s_{i_2}\cdots s_{i_{\ell+1}}$ by Proposition \ref{prop:ReflCusp}. 

\smallskip

Now we shall prove (i). It is known that  
\be[(1)]
\item any simple module in $ \cC_{*, s_i}$ is isomorphic to the head of a convolution product of copies of $\Cp_2, \ldots, \Cp_\ell $, 
\item for any simple module $M \in \RC\gmod$, $M \in  \cC_{*, s_i} $ if and only if $\eps_i^*(M)=0$ 
\ee
(see \cite[Proposition 2.18 and Theorem 2.20]{KKOP18}).
Hence (a)$\Leftrightarrow$ (b) follows from (1) and Corollary  \ref{Cor: ep under F}.

Note that
$$\F_\ddD(\Cp_k)=\cuspS_k=\cuspS'_{k-1}\in  \catC_{\Refl_i(\ddD)}
\qt{for $2\le k\le\ell$.}$$
Hence, by (1), we have  
$$
\F_\ddD( \cC_{*, s_i}  ) \subset \catC_{\Refl_i(\ddD)} \cap \catCD.
$$
that is (a)$\Rightarrow$(c).

Let $S$ be a simple module in $\catC_{\Refl_i(\ddD)} \cap \catCD$. Since 
$\D^{-1}\Rt_i\in\D^{-2}\catC_{\Refl_i(\ddD)}$, 
 we have 
$$
\de( \dual^{-1} \Rt_i, S )= 0.
$$
Thus we obtain (c)$\Rightarrow$(b).
\end{proof}

\begin{example} \label{Ex: ex2}   
	
We use the same notations given in Example \ref{Ex: ex1}.
\bni
\item We shall apply $\Refl_1$ to the duality datum $\ddD = \{ \Rt_1, \Rt_2 \}$. Let 
\begin{align*}
	\widetilde{ \Rt}_1 &:= \Refl_1 (\Rt_1) = \dual \Rt_1 = V(2)_{(-q)^3}, \\
	\widetilde{ \Rt}_2 &:= \Refl_1 (\Rt_2) = \Rt_1 \htens \Rt_2 = V(2)_{-q}.
\end{align*}
Then we have $ \Refl_1 (\ddD) = \{ \widetilde{ \Rt}_1, \widetilde{ \Rt}_2 \}$. 
The affine cuspidal modules $\widetilde{\cuspS}_k$ corresponding to $ \Refl_1 (\ddD)$ and the reduced expression $s_2s_1s_2$ are given as follows:
\begin{align*}
	\widetilde{ \cuspS}_1 &=\F_{\Refl_1 \ddD} (L(2)) = \widetilde{ \Rt}_2 = V(2)_{-q}, \\
	\widetilde{\cuspS}_2 &=\F_{ \Refl_1 \ddD} (L(2) \hconv L(1)) = \widetilde{ \Rt}_2 \htens \widetilde{ \Rt}_1 = V(2)_{-q} \htens V(2)_{(-q)^3} = V(1)_{( -q)^2}, \\
	\widetilde{\cuspS}_3 &=\F_{ \Refl_1 \ddD} (L(1)) = \widetilde{ \Rt}_1 = V(2)_{(-q)^3},
\end{align*}
and $\widetilde{ \cuspS}_{ k+3} = \dual ( \widetilde{ \cuspS}_k)$ for $k\in \Z$. Note that $ \widetilde{\cuspS}_k = \cuspS_{ k+1}$ for any $k\in \Z$ (see Proposition \ref{prop:ReflCusp}).

\item We shall apply $\Refl_2$ to the duality datum $\ddD = \{ \Rt_1, \Rt_2 \}$. Let 
\begin{align*}
	\widehat{ \Rt}_1 &:= \Refl_2 (\Rt_1) = \Rt_2 \htens \Rt_1 =  V(1)_{(-q)^2} \htens V(1) = V(1^2)_{-q}, \\
	\widehat{ \Rt}_2 &:= \Refl_2 (\Rt_2) = \dual \Rt_2 = V(2)_{( -q)^5}.
\end{align*}
Then we have $ \Refl_2 (\ddD) = \{ \widehat{ \Rt}_1, \widehat{ \Rt}_2 \}$. As you observe, the duality datum $ \Refl_2 (\ddD) $ has a root module which is not fundamental.

The affine cuspidal modules $\widehat{\cuspS}_k$ corresponding to $ \Refl_2 (\ddD)$ and the reduced expression $s_1s_2s_1$ are given as follows:
\begin{align*}
	\widehat{ \cuspS}_1 &=\F_{ \Refl_2 \ddD} (L(1)) = \widehat{ \Rt}_1 = V(1^2)_{-q}, \\
	\widehat{\cuspS}_2 &=\F_{ \Refl_2 \ddD} (L(1) \hconv L(2)) = \widehat{ \Rt}_1 \htens \widehat{ \Rt}_2 = V(1^2)_{-q} \htens V(2)_{(-q)^5} \\ 
	& =  ( V(1)_{( -q)^2} \htens V(1) ) \htens V(2)_{(-q)^5}  = V(1) \qquad \qquad \text{by Lemma \ref{Lem: MNDM},} \\
	\widehat{\cuspS}_3 &=\F_{ \Refl_2 \ddD}(L(2)) = \widehat{ \Rt}_2 = V(2)_{(-q)^5},
\end{align*}
and $\widehat{ \cuspS}_{ k+3} = \dual ( \widehat{ \cuspS}_k)$ for $k\in \Z$.  
Note that $ \widehat{\cuspS}_k = \cuspS'_{ k+1}$ for any $k\in \Z$ (see Proposition \ref{prop:ReflCusp}).
	
	\ee  
\end{example}

\vskip 2em 

\section{PBW theoretic approach} \label{Sec: PBW} \

\subsection{\Principal duality datum} \label{Sec: CDD}

\begin{definition} \label{Def: principal}
A duality datum $\ddD$ is called \emph{\principal} if it is strong and, for any simple module $M \in \catCO$, there exists simple modules $M_k \in \catCD$ $(k\in \Z)$ such that 
\bna
\item $M_k \simeq \one $ for all but finitely many $k$,
\item $M \simeq \hd ( \cdots \tens \dual^2 M_2 \tens \dual M_1 \tens \ M_0 \tens \dual^{-1} M_{-1} \tens \cdots  ).$
\ee
\end{definition}

In \cite{KKOP20}, we associate 
to the category $\catCO$ a simply laced finite type root system
in a canonical way.
For a simple  module $M \in \catC_\g$, set $ \rE(M) \in \Hom(\sig, \Z)$ by 
\begin{align*}
\rE(M)(i,a)\seteq \Li(M, V(\varpi_i)_a) \qquad \text{ for } (i,a) \in \sig.
\end{align*}
Let
\begin{equation*}
\begin{aligned}
\rW_0 \seteq \{  \rE(M) \mid M \text{ is simple in } \catC_\g^0 \}
\qtq\Delta_0 \seteq \{ \rs_{i,a} \mid (i,a) \in \sigZ  \} \subset \rW_0,
\end{aligned}
\end{equation*}
where we set $\rs_{i,a} \seteq  \rE(V(\varpi_i)_a) $. Then $\rsP\seteq  ( \rW_0, \Delta_0 )$ forms a root system, and the type of $\rsP$ is given as follows (\cite[Theorem 3.6]{KKOP20}): 
\renewcommand{\arraystretch}{1.5}
\begin{align} \label{Table: root system} \small
\begin{array}{|c||c|c|c|c|c|c|c|} 
\hline
\text{Type of $\g$} & A_n^{(1)}  & B_n^{(1)} & C_n^{(1)} & D_n^{(1)} & A_{2n}^{(2)} & A_{2n-1}^{(2)} & D_{n+1}^{(2)}  \\
&(n\ge1)&(n\ge2)&(n\ge3)&(n\ge4)&(n\ge1)&(n\ge2)&(n\ge3)\\
\hline
\text{Type of $\rsP$} & A_n & A_{2n-1}    & D_{n+1}   &  D_n & A_{2n} & A_{2n-1} & D_{n+1}  \\
\hline
\hline
\text{Type of $\g$} & E_6^{(1)}  & E_7^{(1)} & E_8^{(1)} & F_4^{(1)} & G_{2}^{(1)} & E_{6}^{(2)} & D_{4}^{(3)}  \\
\hline
\text{Type of $\rsP$} & E_6 & E_{7}    & E_{8}   & E_6 & D_{4} & E_{6} & D_{4}  \\
\hline
\end{array}
\end{align} 
We denote by $\rsX$ the type of $\rsP$.

We define a symmetric bilinear form $(\ ,\ )$ on $\rW_0$ by
$( \rE(M), \rE(N) )=- \Li(M,N)$ for simple modules $M$ and $N$.
Then $(\ ,\ )$ is a Weyl group invariant positive definite bilinear form and
$\Delta_0=\set{\al\in\rW_0}{(\al,\al)=2}$.

\begin{prop} \label{Prop: cd Xg}
Let $ \ddD\seteq  \{ \Rt_i \}_{i\in J} \subset \catCO$ be a \principal duality datum  associated with a simply-laced finite Cartan matrix $\cmC$. 
Then $\cmC$ is of type $\rsX$.
\end{prop}
\begin{proof}
We denote by $\rlQ_\cmC$ and $\rDC$  the root lattice and the set of roots associated with $\cmC$.

It follows from Proposition \ref{Prop: Li}, Proposition \ref{Prop: subquotients for Li}, Theorem \ref{Thm: simple to simple} and Definition \ref{Def: principal}, 
the abelian group $\rW_0$ is generated by $\rE(M)$ for $M\in \catCD$.
Moreover $\rE(\F_\ddD(M))$ depends only on $\wt(M)$ by Theorem~\ref{Thm: invariants preserve}\,\eqref{it:wt}.
Hence the functor $\F_\ddD$ induces the surjective additive map 
$$
[\F_\ddD]\col \rlQ_\cmC \epito\rW_0 
$$
given by $[\F_\ddD](\al_i) = \rE(\Rt_i)$ for $i\in J$. 
Moreover, $[\F_\ddD]$ preserves the positive definite
pairing $(-,-)$. Hence $[\F_\ddD]$ is bijective. 
Since both of  $\rDC$  and $\Delta_0$ are characterized by the condition $(X,X)=2$ (\cite[Corollary 3.8]{KKOP20} and \cite[Proposition 5.10]{Kac}),
the set $\{  \rE(\Rt_i) \}_{i\in J}$ becomes 
a  basis of the root system $\rsP$. 
Since $c_{i,j} = (\al_i, \al_j) = (\rE(\Rt_i), \rE(\Rt_j))$ for any $i,j \in J$ by Theorem~\ref{Thm: invariants preserve}\,\eqref{it:wt},
we conclude that the Cartan matrix $\cmC = (c_{i,j})_{i,j\in J}$ is of type $\rsX$.
\end{proof}

\begin{theorem} \label{Thm: refl D}
Let $\ddD \seteq  \{ \Rt_i \}_{i\in J}$ be a \principal duality datum.  For any $i\in J$,  $\Refl_i(\ddD)$ and $\Refl_i^{-1}(\ddD)$ are \principal.
\end{theorem}
\begin{proof} 
We shall focus on proving the case for $\Refl_i$ since the other case for $\Refl_i^{-1}$ can be proved in a similar manner.
Since $\Refl_i(\ddD)$ is strong by  Proposition \ref{Prop: strong to strong}, it suffices to show that $\Refl_i(\ddD)$  satisfies the conditions of Definition \ref{Def: principal}.

\smallskip

Let $i\in J$ and choose a reduced expression $\rxw = s_{i_1} s_{i_2} \cdots s_{i_\ell}$ of the longest element $w_0$ of $\weylC$ with $i_1 = i$.
Define $\st{i_k}_{k\in\Z}$ and the cuspidal modules $\st{\cuspS_k}_{k\in\Z}$
corresponding to $\ddD$ and $\rxw$ as in \S\;\ref{subsec:cusp}.
Let $M$ be a simple module in $\catCO$. 
As $\ddD$ is \principal, there exist simple modules $M_k \in \catCD$ ($k\in \Z$) such that 
$M_k \simeq \one $ for all but finitely many $k$ and 
\begin{align} \label{Eq: M=hd Mks}
M \simeq \hd \left( \cdots \tens \dual^2 M_2 \tens \dual M_1 \tens \ M_0 \tens \dual^{-1} M_{-1} \tens \cdots  \right).
\end{align}
For each $k\in \Z$, there exist $a_{k,1}, \ldots, a_{k, \ell} \in \Z_{\ge 0}$ such that
$$
M_k\simeq\hd\bl\cuspS_{\ell}^{\tens a_{k,\ell}}\tens\cdots\tens\cuspS_1^{\tens a_{k,1}}\br.
$$
Set $c_{s+k\ell}=a_{k,s}$ for $1\le s\le \ell$ and $k\in\Z$.
Then, by Lemma \ref{Lem: dual head}, we have
$$\D^kM_k\simeq\hd\bl
\cuspS_{k\ell+\ell}^{\tens c_{k\ell+\ell}}\tens\cdots\cuspS_{k\ell+1}^{\tens c_{k\ell+1}}\br. 
$$
Hence we have
$$M\simeq\hd(\cdots \tens\cuspS_{1}^{\tens c_1}\tens
\cuspS_0^{\tens c_0}\tens \cuspS_{-1}^{\tens c_{-1}}\tens\cdots\br.$$

Set
$$N_k=\hd\bl\cuspS_{\ell+1}^{\tens c_{k\ell+\ell+1}}\tens\cdots\tens\cuspS_{2}^{\tens c_{k\ell+2}}\br.$$
Then 
$N_k\in\catC_{\Refl_i\ddD}$ by Proposition~\ref{prop:ReflCusp}, and we have
$$\D^k N_k\simeq\hd\bl\cuspS_{k\ell+\ell+1}^{\tens c_{k\ell+\ell+1}}\tens\cdots\tens\cuspS_{k\ell+2}^{\tens c_{k\ell+2}}\br.$$
Hence we obtain
$$M\simeq\hd\bl\cdots\tens\D^1N_1\tens \D^0 N_0\tens \D^{-1} N_{-1}\tens\cdots\br.
$$
\end{proof}

\subsection{Duality datum arising from $\mathrm{Q}$-datum} \label{Sec: Q-data} \

The subcategory $\catCQ$ of $\catCO$ was introduced in \cite{HL15} for simply-laced affine type ADE, in \cite{KKKO16D} for twisted affine type $A^{(2)}$ and $D^{(2)}$, 
in  \cite{KO18, OhSuh19} for untwisted affine type $B^{(1)}$ and $C^{(1)}$, and in \cite{OS19} for exceptional affine type.
Let $\gf$ be the simple Lie algebra of type $\rsX$ defined in $\eqref {Table: root system}$ and $\If$ the index set of $\gf$.
 The category $\catCQ$ categorifies the coordinate ring $\C[N]$ of the maximal unipotent group $N$
associated with $\gf$. This category is defined by a Q-datum.  
A Q-datum is a triple $ \qQ \seteq (\Dynkin, \Dat, \hf)$  consisting of the Dynkin diagram $ \Dynkin$ of $\gf$, an automorphism $\Dat$ on $\Dynkin$ and a height function $\hf$, which satisfy certain conditions (see \cite{FO20} for details,
and also \cite[\S\,6]{KKOP20B}).
When $\g$ is of untwisted affine type ADE, $\Dat$ is the identity and $\qQ$ is equal to a Dynkin quiver with a height function.
To a Q-datum $\qQ$, 
we can associate a subset $\sigQ$ of $\sigZ$.
This set $\sigQ$ is in a 1-1 correspondence to the set $\prDf$ 
of positive roots of $\gf$, 
which is denoted by 
\begin{align} \label{Eq: phi}
\phiQ\col \prDf \isoto\sigQ.
\end{align} 
Set 
$$
\ddD_\qQ \seteq  \{ \Rt_i \}_{i\in \If},
$$
where $\Rt_i$ is the fundamental module corresponding to $\phiQ(\alpha_i)$  for $i\in \If$. Then $\ddD_\qQ$ becomes a strong duality datum 
(\cite{KKK15B, KKKO16D, KO18, OS19, Fu18, FO20,KKOP20B}), 
which gives the duality functor $\F_{\ddD_\qQ}$.
By the definition, we have $\catCQ = \catC_{\ddD_\qQ}$. We simply write $\F_\qQ$ for $\F_{\ddD_\qQ}$:
$$
\F_\qQ\cl \Rgf \gmod\To\catC_\g.
$$

We refer the reader to \cite{FO20, OhSuh19,KKOP20B} 
for the notion of (twisted) $\qQ$-adapted reduced expressions 
of the longest element $w_0$ of the Weyl group of
$\gf$.

Let $\weylf$ be the Weyl group of $\gf$. 
For a $\qQ$-adapted reduced expression $\rxw= s_{  i_1 } s_{  i_2 } \cdots s_{  i_\ell }$
of the longest element $w_0$ of $\weylf$, 
we define $\beta_k\in\prDf$ ($1\le k\le \ell)$
by \eqref{Eq: beta_k}.   
Then there exist a sequence $ \{(i_k, a_k)\}_{k\in \Z} \subset \If \times \cor^\times$ and $\pi\colon \If\to I_0$ such that 
$(\pi(i_k), a_k)=  \phiQ  (\beta_k) \in \sigQ $ for $ k=1, \ldots, \ell$ and 
$$
(\pi( i_{s+m\ell} ), a_{s+m\ell})=\delta^m\bl( \pi(i_s),a_s)\br
\qt{for $1\le s\le \ell$ and $m\in\Z$.}
$$
Here we set 
$$
\delta^m\bl(i,a)\br\seteq  
\begin{cases}
 \bl i, (p^*)^ma\br & \text{ if $m$ is even}, \\
 \bl i^*, (p^*)^ma\br & \text{ if $m$ is odd}.
\end{cases}
$$
(See \cite[\S\,6]{KKOP20B}.)

We define the affine cuspidal modules $\st{\cuspS_k}_{k\in\Z}$ as in Definition~\ref{def:cusp}.

Collecting results in \cite{FO20, HL15,KKK15B,  KKKO16D,KO18,OS19,OhSuh19}, we have Proposition \ref{lem:adapted} below. 
In Proposition \ref{lem:adapted},  the symmetric cases follow from \cite{ HL15, KKK15B}, the untwisted $B^{(1)}$ and $C^{(1)}$ cases follow from \cite{OhSuh19, KO18}, the twisted $A^{(2)}$ and $D^{(2)}$ cases follow from \cite{KKKO16D}, 
and the exceptional cases follow from \cite{OS19}.   The uniform approach is given in \cite{FO20}. See also \cite[\S\,6]{KKOP20B}.

\Prop[\cite{FO20, HL15, KKK15B, KKKO16D,KO18,OS19,OhSuh19}]\label{lem:adapted}
 Let $\qQ$ be a Q-datum.
\bnum
\item $\sigma_0(\g)=\bigsqcup_{m\in\Z}\;\delta^m\,\sigQ$ {\rm(see \cite[Proposition 4.21]{FO20} for example)}. 
\item
There exists a $\qQ$-adapted reduced expression of $w_0$ {\rm(see \cite[Section 3]{OhSuh19} for example)}. 
\item For a $\qQ$-adapted reduced expression 
 $\rxw=s_{i_1} s_{i_2} \cdots s_{i_\ell}$ of $w_0$,  let
$\st{(i_k,a_k)}_{k\in\Z}$ be the sequence as above, and let $\st{\cuspS_k}_{k\in\Z}$ be the affine cuspidal modules corresponding
$\ddD_\qQ$ and $\rxw$.
Then we have 
\bna
\item $\cuspS_k\simeq V(\vpi_{\pi( i_k)} )_{a_k}$,
\item $d_{V(\vpi_{ \pi(i_s) }), V(\vpi_{\pi (i_t) })}(a_t/a_s)\not=0$
for $t,s\in\Z$ such that $s>t$.
Here, $d$ is the denominator of the R-matrix.
\ee
{\rm (See \cite[Theorem 4.3.4]{KKK15B}, \cite[Theorem 5.1 and Lemma 5.2]{KKKO16D} , \cite[Theorem 6.3, 6.4]{KO18} and \cite[Section 6]{OS19}).}
\ee
\enprop

\begin{prop} \label{Prop: DQ complete}
 The duality datum  $\ddDQ$ is a \principal duality datum.
\end{prop}
\begin{proof}
Recall that $\sigma_0(\g)=\set{(\pi(i_k),a_k)}{k\in\Z}$.
For a simple module $M$ in $\catCO$,
let $\la=\sum_{s=1}^r(\pi( i_{k_s} ),a_{k_s})$ be the affine highest weight of $M$
(see Theorem~\ref{Thm: basic properties}\;\eqref{Thm: bp5}). 
We may assume that $\st{k_s}_{s\in[1,r]}$ is a decreasing sequence.
Then, by Proposition~\ref{lem:adapted} and 
Theorem~\ref{Thm: basic properties}, we have
$M\simeq\hd\bl\cuspS_{k_1}\tens\cdots\tens\cuspS_{k_r}\br$. 
\end{proof}
Thanks to Theorem \ref{Thm: refl D}, we have the following.

\begin{corollary} \label{Cor: DQ complete}
The duality datum obtained from $ \ddDQ$ by applying a finite sequence of $\Refl_i$ and  $\Refl_i^{-1}$ $(i\in \If)$  is a \principal duality datum.
\end{corollary}

\begin{example}  
We use the same notations given in Example \ref{Ex: ex1} and Example \ref{Ex: ex2}. Let $ \Dynkin$ be the Dynkin diagram of finite type $A_2$.

\bni
\item Let $\xi$ be the height function on $\Dynkin$ defined by $\xi(1) = 0$ and $\xi(2)=1$, and let $\qQ$ be the Q-datum consisting of $\Dynkin$ and $\xi$. 
Then $\ddD$ is equal to the duality datum arises from  the Q-datum $\qQ$, which says that $\ddD$ is complete. The reduced expression $s_1s_2s_1$ is $\qQ$-adapted, but $s_2s_1s_2$ is not $\qQ$-adapted.

\item By Corollary \ref{Cor: DQ complete}, $\Refl_1 (\ddD)$ and $\Refl_2 (\ddD)$ are complete duality data. The duality datum $\Refl_1 (\ddD)$ arises from the Q-datum consisting of $\Dynkin$ and the height function $\xi'$ defined by $\xi'(1) = 2$ and $\xi'(2)=1$, but  $\Refl_2 (\ddD)$ does not come from any Q-datum.	
\ee
\end{example}

\subsection{PBW for quantum affine algebras} \

In this subsection, we develop the PBW theory for $\catCO$ using a complete duality datum. 
This generalizes the ordinary standard modules and related results (\cite{H04, Kas02, Nak01, Nak04, VV02}). 
Note that the ordinary standard modules are cyclic tensor products of fundamental modules.

Let  $\cmC=(c_{i,j})_{i.j\in J}$ be a simply-laced finite Cartan matrix. 
Throughout this subsection, we assume that  
$$
\text{$ \ddD= \{ \Rt_i \}_{i\in J}$ is a \emph{\principal} duality datum associated with  $\cmC$. }
$$
Proposition \ref{Prop: cd Xg} says that $\cmC$ is of type $\rsX$ and $J=\If$. Let $\weylC$ be the Weyl group associated with $\cmC$. 
We fix a reduced expression $\rxw = s_{i_1} s_{i_2} \cdots s_{i_\ell}$ of the longest element $w_0$ of $\weylC$, and 
let $\cuspS_k$ $(k\in \Z)$ be the affine cuspidal modules corresponding to $\ddD$ and $\rxw$. 
We define 
\begin{align} \label{Eq: ZZ}
 \ZZ \seteq \Z_{\ge0}^{\oplus \Z}
= \bigl\{ (a_k)_{k\in \Z} \in   \Z_{\ge0}^{\Z} \mid
  \text{  $a_k=0$ except finitely many $k$'s} 
\bigr\}.
\end{align}
We denote by $\prec$ the bi-lexicographic order on $\ZZ$, i.e.,   
 for any $ \bfa = (a_k)_{k\in \Z}$ and $\bfa' = (a_k')_{k\in \Z}$  in $\ZZ$, $ \bfa \prec \bfa' $ if and only if the following conditions hold:
\eq\label{eq:bilexico}
&\phantom{aa}&\left\{\parbox{68ex}{\bna
\item there exists $r \in \Z$ such that $ a_k = a_k' $ for any $k < r$ and $ a_r < a_r'$,
\item   there exists $s \in \Z$ such that $ a_k = a_k' $ for any $k > s$ and $ a_s < a_s'$.
\ee
}\right.
\eneq
Similarly, we set $\prec_r$ (resp.\ $\prec_l$) to the right (resp.\ left) lexicographic order on $\ZZ$, i.e.,   
 for any $ \bfa ,\bfa'  \in \ZZ$, $ \bfa \prec_r \bfa' $ (resp.\ $ \bfa \prec_l \bfa' $) if and only if  the condition (a) (resp.\ (b)) in $\eqref{eq:bilexico}$ holds.
Hence, we have 
\begin{align} \label{Eq: l r lexico}
\bfa \prec \bfa' \quad \Longleftrightarrow \quad  \bfa \prec_l \bfa' \ \text{ and } \    \bfa \prec_r \bfa'.
\end{align}

For $ \bfa = (a_k)_{k\in \Z} \in \ZZ$, we define
\eqn
&&\sP_{\ddD, \rxw} (\bfa)\seteq
 \bigotimes_{k =+\infty}^{-\infty}  \cuspS_k^{\otimes a_k}
=\cdots\tens\cuspS_2^{\otimes a_2}\tens\cuspS_1^{\otimes a_1} \otimes  \cuspS_{0}^{\otimes a_{0}}
 \otimes \cuspS_{-1}^{\otimes a_{-1}}
\otimes\cuspS_{-2}^{\otimes a_{-2}}\otimes \cdots.
\eneqn
Here $\sP_{\ddD, \rxw} (0)$ should be understood as 
the trivial module $\one$.
We call the modules $\sP_{\ddD, \rxw} (\bfa)$ \emph{standard modules} with respect to the cuspidal modules $\{ \cuspS_k \}_{k\in \Z}$.  

\begin{lemma}\label{Lem: d(DS, SM)}
Let $k\in \Z$, $a \in \Z_{> 0}$ and let $M$ be a simple module in $\catCO$. 
\bnum
\item If $ \de( \dual^t \cuspS_k, M )=0$ for $t=1,2$, then $ a = \de(\dual \cuspS_k, \cuspS_k^{\tens a} \htens M) $.
\item If $ \de( \dual^{t} \cuspS_k, M )=0$ for $t=-1,-2$, then $ a = \de(\dual^{-1} \cuspS_k,  M \htens  \cuspS_k^{\tens a})$.
\ee
\end{lemma}
\begin{proof}
(i) Note that $\cuspS_k$ is a root module by Proposition \ref{Prop: cusp}.
Applying Lemma \ref{Lem: DLX = LX+1} (i) to the setting  $ L\seteq\cuspS_k$ and $X \seteq M$, 
we have 
$$
\de(\dual \cuspS_k, \cuspS_k^{\tens a} \htens M) = a + \de(\dual \cuspS_k,  M) = a.
$$

\snoi
(ii) can be proved in the same manner as above.
\end{proof}

\begin{lemma} \label{Lem: uniqueness}
Let  $ m, l \in \Z $ with $m\ge l$ and $ a_m, a_{m-1} \ldots, a_l \in \Z_{\ge 0}   $. We set 
$$
M \seteq  \hd \left( \cuspS_{m}^{\tens  a_{ m}} \tens \cuspS_{ m-1}^{\tens  a_{m-1}} \tens \cdots \tens \cuspS_{l}^{\tens a_{l}}  \right).
$$
\bnum
\item $\de( \dual \cuspS_k, M)=0$ for any $k > m$.
\item Set $M_m \seteq  M$ 
and define inductively 
$$
d_{k} \seteq  \de( \dual \cuspS_{k}, M_{k}  )\qtq
M_{k-1} \seteq  M_k \htens \dual ( \cuspS_{k}^{\tens d_k }),  \qquad 
$$
for $k=m, \ldots,l$. 
Then we have
$$
d_k = a_k \qtq M_k\simeq
\hd \left( \cuspS_{k}^{\tens  a_{k}} \tens \cuspS_{k-1}^{\tens  a_{k-1}} \tens \cdots \tens \cuspS_{l}^{\tens a_{l}}  \right) 
\qt{for $k=m,\ldots, l$.} 
$$
\item $\de( \dual^{-1} \cuspS_k, M)=0$ for any $k < l$.
\item Set $N_l \seteq  M$ 
and define inductively 
$$
e_{k} \seteq  \de( \dual^{-1} \cuspS_{k}, N_{k}  )\qtq
N_{k+1} \seteq   \dual^{-1} ( \cuspS_{k}^{\tens e_k }) \htens N_k,  \qquad 
$$
for $k=l, \ldots,m$.  Then we have 
$$
e_k = a_k \qtq
N_k\simeq \hd \left( \cuspS_{m}^{\tens  a_{ m}} \tens \cdots\tens
\cuspS_{ k+1}^{\tens  a_{k+1}}\tens \cuspS_{k}^{\tens a_{k}}  \right)
\qt{for $k=m, \ldots, l$. }
$$
\ee
\end{lemma}
\begin{proof}
(i) By Proposition \ref{Prop: cusp} (ii), $( \cuspS_k, \cuspS_t )$ is strongly unmixed  for any $k > m$ and $t=m, \ldots, l$.  Thus we have $\de( \dual\cuspS_k, \cuspS_t)=0$ for $t=m, \ldots , l$, which implies that
$ \de(\dual \cuspS_k, M)=0 $.

\mnoi
(ii) 
By induction on $k$, we may assume that $k=m$.
We set 
$ N \seteq  \hd \left( \cuspS_{m-1}^{\tens  a_{ m-1}} \tens \cdots \tens \cuspS_{l}^{\tens a_{l}}  \right)$.
By (i), we have $\de( \dual^t \cuspS_m, N )=0$ for $t=1,2$.
Proposition \ref{Prop: cusp} (iii) tells us that $M \simeq\cuspS_m^{\tens a_m} \htens N $. 
Thus, by  Lemma \ref{Lem: MNDM} and Lemma \ref{Lem: d(DS, SM)}, we have  
\begin{align*}
d_m=\de( \dual \cuspS_m, M) &= \de( \dual \cuspS_m, \cuspS_m^{\tens a_m} \htens N) = a_m, \\
 M \htens \dual (\cuspS_k^{\tens a_m}) &\simeq  (\cuspS_m^{\tens a_m} \htens N) \htens \dual (\cuspS_k^{\tens a_m}) \simeq N.
\end{align*}
\snoi
The assertions (iii) and  (iv) can be proved in the same manner as above.
\end{proof}

\begin{theorem} \label{Thm: main1} \label{Thm: main pbw} \

\bnum
\item For any $\bfa \in \ZZ$, the head of \/ $\sP_{\ddD, \rxw} (\bfa)$ is simple. We denote the head by 
$$
 \sV _{\ddD, \rxw} (\bfa) \seteq  \hd \bl\sP_{\ddD, \rxw} (\bfa)\br.
$$

\item For any simple module $M \in \catCO$, there exists a unique $\bfa \in \ZZ$ such that 
$$
M \simeq  \sV _{\ddD, \rxw} (\bfa).
$$
\ee
Therefore, the set $\{   \sV _{\ddD, \rxw} (\bfa) \mid \bfa \in \ZZ \}$ is a complete and irredundant set of simple modules of $\catCO$ up to isomorphisms.
\end{theorem}
\begin{proof}
(i) follows from Proposition \ref{Prop: cusp}.

\mnoi
(ii) Let $M$ be a simple module in $\catCO$. Since $\ddD$ is \principal, there exist simple module $M_k \in \catCD$ ($k\in \Z$) such that 
$ M_k \simeq \one$ for all but finitely many $k$ and 
$$
M \simeq \hd \left( \cdots \tens \dual^2 M_2 \tens \dual M_1 \tens  M_0 \tens \dual^{-1} M_{-1} \tens \cdots  \right).
$$
Since $M_k \in \catCD$, there exist $b_{1}^k , \ldots, b_\ell^k \in \Z_{\ge0}$ such that 
$M_k \simeq \hd \left(  \cuspS_{\ell}^{b_\ell^k} \tens \cdots \tens  \cuspS_{1}^{b_1^k}  \right)$, which yields that 
$$
\dual^k M_k \simeq \hd \left(  \cuspS_{ k \cdot \ell + \ell}^{b_\ell^k} \tens \cdots \tens  \cuspS_{k \cdot \ell + 1}^{b_1^k}  \right)
$$
by Lemma \ref{Lem: dual head}. 
For $t\in \Z$, we define $a_t \seteq  b_r^k$, where $t = k \cdot \ell + r$ for some $k\in \Z$ and $r=1, \ldots, \ell$, and set $ \bfa \seteq  (a_t)_{t\in \Z}$.
By Proposition \ref{Prop: cusp}, we have 
$$
M \simeq  \sV _{\ddD, \rxw} (\bfa).
$$
The uniqueness for $\bfa$ follows from Lemma \ref{Lem: uniqueness}, which completes the proof.

\end{proof}

The element $\bfa \in \ZZ$ associated with a simple module $M$
in Theorem~\ref{Thm: main1}~(ii) is called the \emph{cuspidal decomposition} of $M$ with respect to the cuspidal modules $\{ \cuspS_k \}_{k\in \Z}$,
and it is denoted by
 \begin{align} \label{Eq: cusp decomp}
\bfa_{\ddD, \rxw}(M) \seteq \bfa.
\end{align}

\begin{lemma} \label{Lem: LMLN MN}
Let $L,M,N$ be simple modules in $\cat$ and assume that $L$ is real.
\bnum
\item 
If $(L,M)$ and $(L,N)$ are strongly unmixed and $L \hconv N$ appears in $L \tens M$ as a subquotient, then we have $M \simeq N$.

\item If $(M,L)$ and $(N,L)$ are strongly unmixed and $N \hconv L$ appears in $M \tens L$ as a subquotient, then we have $M \simeq N$.
\ee
\end{lemma}
\begin{proof}
(i) Since $(L,M)$ and $(L,N)$ are strongly unmixed, 
$$
\La(L,M) = \Li(L,M) \quad \text{ and }\quad \La(L,N) = \Li(L,N)
$$
by Lemma \ref{Lem: unmixed Li=La}. Since $L \htens N$ appears in $L\tens M$, Proposition \ref{Prop: subquotients for Li} tells us that 
$$
\La(L,M) = \Li(L,M) = \Li(L,N) = \La(L,N) =\La(L,L\hconv N).
$$
Thus it follows from \cite[Theorem 4.11]{KKOP19C} that 
$$
L \htens M \simeq L \htens N,
$$
which implies that $M \simeq N$ by Lemma \ref{Lem: MNDM}.

\snoi
(ii) can be proved in the same manner as above.
\end{proof}

For $\bfc = (c_k)_{k\in \Z} \in \ZZ$, we set $l(\bfc)$ (resp.\ $r(\bfc)$) to be the integer $t$ such that
\begin{align} \label{Eq: def lr}
c_t \ne 0, \quad \quad \text{ $c_k=0$ for any $k > t$ (resp.\ $k < t$).}
\end{align}

\begin{theorem} \label{Thm: main tri}
Let $\bfa$ be an element of $ \ZZ$. Then we have the following.
\bnum
\item The simple module $ \sV _{\ddD, \rxw} (\bfa) $ appears only once in $\sP_{\ddD, \rxw} (\bfa)$.
\item If $V$ is a simple subquotient of $\sP_{\ddD, \rxw} (\bfa)$ which is not isomorphic to $ \sV _{\ddD, \rxw} (\bfa) $, then we have 
$$
\bfa_{\ddD, \rxw}(V) \prec \bfa.
$$
\item In the Grothendieck ring, we have 
$$
[\sP_{\ddD, \rxw} (\bfa)] = [\sV_{\ddD, \rxw} (\bfa)] + \sum_{\bfa' \prec \bfa} c(\bfa') [\sV_{\ddD, \rxw} (\bfa')],
$$
for some $ c(\bfa') \in \Z_{\ge 0}$.
\ee
\end{theorem}

\begin{proof}
We focus on proving (ii) because (i) and (iii) follow from (ii).

Let $\bfa = ( a_k )_{k\in \Z}$ and set 
$$
l \seteq  l(\bfa) \quad \text{ and } \quad r\seteq  r(\bfa).
$$
Let $V$ be a simple subquotient of $\sP_{\ddD, \rxw} (\bfa)$ which is not isomorphic to $ \sV _{\ddD, \rxw} (\bfa) $. We set 
$$
\bfb  = (b_k)_{k\in \Z}\seteq  \bfa_{\ddD, \rxw}(V). 
$$
For $ k > l$  and $r > t$, 
since  $( \cuspS_{k}, \cuspS_l )$ and $( \cuspS_{r}, \cuspS_t )$ is strongly unmixed by Proposition \ref{Prop: cusp}, 
we have 
$$
\de( \dual \cuspS_{k}, \sP_{\ddD, \rxw} (\bfa)  ) = 0 ,\qquad \de( \dual^{-1} \cuspS_{t}, \sP_{\ddD, \rxw} (\bfa)  ) = 0,
$$
which implies that 
$\de( \dual \cuspS_{k}, V  ) = 0 $ and $ \de( \dual^{-1} \cuspS_{t}, V) = 0$ by \cite[Proposition 4.2]{KKOP19C}.
Thus, Lemma \ref{Lem: uniqueness} tells us that 
$$
l \ge l(\bfb) \quad \text{ and } \quad  r(\bfb) \ge r.
$$

\smallskip

We now shall prove $\bfb \prec_l  \bfa$, where $\prec_l$ is the left lexicographical order on $\ZZ$.
Note that, by Lemma \ref{Lem: uniqueness}, Proposition \ref{Prop: cusp} and \cite[Proposition 4.2]{KKOP19C}, we have 
$$
b_l = \de (\dual \cuspS_l, V) \le \de (\dual \cuspS_l,  \sP_{\ddD, \rxw} (\bfa) ) = \de (\dual \cuspS_l,  \cuspS_l^{\tens a_l} ) = a_l.
$$ 
When either  $ l > l(\bfb) $ or $l= l(\bfb)$ and $ b_l < a_l$,
it is obvious that $\bfb \prec_l  \bfa$ by the definition.
We assume that $l= l(\bfb)$ and $ b_l = a_l$. Set 
$$
c \seteq  b_l = a_l, \qquad  \bfa^- = (a_k^-)_{k\in \Z}, \text{ where }
a_k^-  \seteq  
\bc
0 & \text{ if } k=l, \\
a_k & \text{ otherwise, }
\ec
$$
and
\begin{align*}
P^- \seteq   \cuspS_{l-1}^{\otimes a_{l-1}}\tens \cdots \tens  \cuspS_{r}^{\otimes a_{r}}, \quad 
V^- \seteq  \hd \left(  \cuspS_{l-1}^{\otimes b_{l-1}}\tens \cdots \tens  \cuspS_{r(\bfb)}^{\otimes b_{r (\bfb)}} \right). 
\end{align*}
Note that 
\begin{align} \label{Eq: P- V-}
P^- = \sP_{\ddD, \rxw} (\bfa^-), \quad  \sP_{\ddD, \rxw} (\bfa) = \cuspS_l^{\tens c} \tens P^-, \qquad  V \simeq (\cuspS_l^{\tens c}) \htens V^- ,
\end{align}
where the third follows from proposition \ref{Prop: cusp} (iii).
As $V$ appears in $\cuspS_l^{\tens c} \tens P^-$ as a simple subquotient, there exist a simple subquotient $L$ of $P^-$ such that 
$$
\text{ $V$ appears in $\cuspS_l^{\tens c} \tens L$ as a simple subquotient.}
$$
By Proposition  \ref{Prop: cusp} (ii), we know that $(\cuspS_l, V^-)$ and $ (\cuspS_l, L)$ are strongly unmixed. 
Hence, by Lemma \ref{Lem: LMLN MN}, we conclude that 
$$
V^- \simeq L.
$$
If $V^- $ is isomorphic to $ \hd (P^-)$, then we have $ V \simeq \hd \bl \sP_{\ddD, \rxw} (\bfa)\br  $ by $\eqref{Eq: P- V-}$, which 
contradicts the assumption. Hence  $V^- $ is not isomorphic to $ \hd( P^-)$. 
Applying the standard induction argument to the setting $V^-$ and $P^-$, we obtain 
$$
\bfa_{\ddD, \rxw}(V^-) \prec_l \bfa^-,
$$
which implies that 
$
\bfb \prec_l  \bfa.
$

In the same manner as above, one can prove that $\bfb \prec_r \bfa$.
Therefore it follows from $\eqref{Eq: l r lexico}$ that 
$$
\bfb \prec  \bfa.
$$
\end{proof}

\begin{remark} \label{Rem: block} 
Let $V$ be a simple subquotient of $\sP_{\ddD, \rxw} (\bfa)$. Theorem \ref{Thm: main tri} says that $ \bfa_{\ddD, \rxw}(V) \prec \bfa. $
There is another condition which $V$ should satisfy.
By Proposition \ref{Prop: subquotients for Li}, we have 
\begin{align} \label{eq: same block}
\rE(V) = \rE(  \sV_{\ddD, \rxw} (\bfa) ),
\end{align}
where $\rE$ is given in Section \ref{Sec: CDD}.
 Thus they are in the same block of $\catC_\g$.
	
\end{remark}

\begin{remark} 
There is a well-known partial ordering, called \emph{Nakajima partial ordering}, in the $q$-character theory.
For simplicity,	we assume that $U_q'(\g)$ is of untwisted affine ADE type. Let $Y_{i,a}$ be an indeterminate for $ i\in I_0$ and $a\in \cor^\times$.
For $i\in I_0$ and $a\in \cor^\times $, set
 $A_{i,a} := Y_{i,aq^{-1}} Y_{i,aq} \prod_{ (\alpha_i, \alpha_j)=-1} Y_{j,a}^{-1} $. 
Then one can define a partial ordering $\le$ on the set of monomials in $\Z[Y_{i,a}^{\pm}\mid i\in I_0, \ a\in \cor^\times  ]$ as follows: for monomials $m$ and $m'$,
$m \le m'$ if and only if $m^{-1} m'$ is a product of elements of $\{ A_{i,a} \mid i\in I_0,\ a\in \cor^\times \}$ (\cite{FM21, Nak04}).
The simple modules and ordinary standard modules in $\catC_\g$ are parameterized by dominant monomials, which are denoted by $L(m)$ and $M(m)$ respectively for a dominant monomial $m$.
 Note that the fundamental module $V(\varpi_i)_a$ corresponds to $Y_{i,a}$.
>From the viewpoint of $(q,t)$-characters, it was shown in \cite{Nak01,Nak04} that 
\begin{align} \label{Eq: Nak}
[M(m)] = [L(m)] + \sum_{m' < m} P_{m,m'} [L(m')]
\end{align}
in the Grothendieck ring $K(\catC_\g)$ and the multiplicity $P_{m,m'}$ can be understood as the specialization at $t=1$ of an analogue $P_{m,m'}(t)$ of Kazhdan-Lusztig polynomial.

Let $\qQ$ be a Q-datum and let $\rxw$ be a $\qQ$-adapted reduced expression. In this case, the affine cuspidal modules $\cuspS_k$ are all fundamental modules in $\catCO$ and $\sP_{\ddD_\qQ, \rxw} (\bfa)$ are ordinary standard modules (see Example \ref{Ex: ex1} (i) for instance).
Let $m$ and $m'$ be dominant monomials and set $ \bfa := \bfa_{\ddD_\qQ, \rxw}(L(m))  $ and 
$ \bfa' := \bfa_{\ddD_\qQ, \rxw}(L(m'))  $.
 Considering the definition of $A_{i,a}$ and \cite[Proposition 6.11]{KKOP20B}, 
one can show that if $m \le m'$ in the partial ordering, then $\bfa \preceq \bfa'$ in the ordering \eqref{eq:bilexico}.
>From this observation about two orders $\le$ and $\preceq$, 
Theorem \ref{Thm: main tri} is compatible with \eqref{Eq: Nak}.
Since affine cuspidal modules are not necessary to be fundamental in general (see Example \ref{Ex: ex1} (ii) for instance), Theorem \ref{Thm: main tri} can be viewed as a generalization of \eqref{Eq: Nak}.

Remark \ref{Rem: block} says that  the condition \eqref{eq: same block} holds when $V$ is a simple subquotient of $\sP_{\ddD, \rxw} (\bfa)$.
Thus it is interesting to ask under what conditions the ordering \eqref{eq:bilexico} is equal to the ordering $\le$.

\end{remark}

For $a,b \in \Z \sqcup \{ \pm \infty \}$,  an \emph{interval}
$[a,b]$ is the set of integers between $a$ and $b$:
\begin{align*}  
 [a,b] \seteq \{ s \in \Z \ | \  a \le s \le b \}.
\end{align*}
If $a>b$, we understand $[a,b]=\emptyset$.

For an interval $[a,b]$, we define 
$ \catCab{[a,b], \ddD, \rxw}$ to be the full subcategory of $\catC_\g$ whose objects have all their composition factors $V$ satisfying the following condition:
\begin{align} \label{Def: C[a,b]} 
b \ge l( \bfa_{\ddD, \rxw}(V) ) \quad \text{ and }\quad r( \bfa_{\ddD, \rxw}(V) )\ge a.
\end{align}
Thanks to Theorem \ref{Thm: main tri},  we have the following proposition.  
\begin{prop}
The category $ \catCab{[a,b], \ddD, \rxw}$ is stable by taking subquotients, extensions, and tensor products.
\end{prop}
It is easy to show that the category $\catCab{[a,b], \ddD, \rxw}$ is equal to the smallest full subcategory of $\catC_0$ satisfying the following conditions$\colon$
\bni
\item  it is stable under taking subquotients, extensions, tensor products and
\item  it contains $\cuspS_{s}$ for all $a \le s \le b$ and the trivial module $\mathbf{1}$.
\ee
If there is no confusion arises, then we simply write $\catCab{[a,b]}$ instead of $\catCab{[a,b], \ddD, \rxw}$.

\smallskip

For an interval $[a,b]$, we set 
$$
\ZZ^{[a,b]} := \{\bfa =(a_{k})_{k\in \Z} \in \ZZ \mid a_k= 0 \text{ for either $k > b$ or $a > k$}  \}
$$
Then the theorem below follows from Lemma \ref{Lem: uniqueness},  Theorem \ref{Thm: main pbw} and Theorem \ref{Thm: main tri} directly.

\begin{theorem} \label{Thm: unitry for C[a,b]}
Let $[a,b]$ be an interval.
\bni
\item  The set $\{   \sV _{\ddD, \rxw} (\bfa) \mid \bfa \in \ZZ^{[a,b]} \}$ is a complete and irredundant set of simple modules of $\catCab{[a,b]}$ up to isomorphisms.
\item 
Let $M$ a simple modules in $\catCO$.
Then, $M$ belongs to $\catCab{[a,b]}$
if and only if
$$\text{$\de(\D \cuspS_k,M)=0$ for $k>b$ and
$\de(\D^{-1}\cuspS_k,M)=0$ for $k<a$.}
$$
\item For $\bfa \in \ZZ^{[a,b]}$, the standard module $\sP_{\ddD, \rxw} (\bfa)$ is contained in $\catCab{[a,b]}$ and,
in the Grothendieck ring, we have 
$$
[\sP_{\ddD, \rxw} (\bfa)] = [\sV_{\ddD, \rxw} (\bfa)] + \sum_{\bfa' \prec \bfa} c(\bfa') [\sV_{\ddD, \rxw} (\bfa')],
$$
for some $ c(\bfa') \in \Z_{\ge 0}$.
\ee
\end{theorem}

\begin{example} 
We use the same notations given in Example \ref{Ex: ex1}. 
\bni
\item We consider the affine cuspidal modules $\cuspS_k$ given in Example \ref{Ex: ex1} (i).
Let $l \in \Z_{\ge0}$. The category $\catCab{[1,2( l+1)]}$ is determined by $\cuspS_k$ for $k\in [1,2( l+1)]$.
It follows from
\begin{align*}
	\{  \cuspS_k \mid k \in [1, 2 ( l+1)] \} = \{ V(1)_{(-q)^{2t}}, V(2)_{(-q)^{2t+1}}  \mid t \in [0,l]  \} 	
\end{align*}
that the category $\catCab{[1,2( l+1)]}$ is equal to the Hernandez-Leclerc category $\catC_l$ defined in \cite[Section 3.8]{HL10}.

\item 
Let us take the affine cuspidal modules $\cuspS_k'$ given in Example \ref{Ex: ex1} (ii).  In this case, the category $\catCab{[a,b]}$ is not equal to $\catC_l$ in general. 
>From this viewpoint, the category $ \catCab{[a,b], \ddD, \rxw}$ is a generalization of the category $\catC_l$.
\ee	
\end{example}

\end{document}